%% file: final_arxiv_main.tex
\def\opre{}

\documentclass[opre]{informs3}

\OneAndAHalfSpacedXI % Current default line spacing

\input{preamble}

\usepackage{amsmath,amssymb,amsfonts}

\usepackage[ruled]{algorithm2e} % For algorithms

\SetAlFnt{\small}
\SetAlCapFnt{\small}
\SetAlCapNameFnt{\small}
\SetAlCapHSkip{0pt}
\IncMargin{-\parindent}

\usepackage{natbib}
 \bibpunct[, ]{(}{)}{,}{a}{}{,}%
 \def\bibfont{\small}%
  \usepackage{multibib}
\newcites{EC}{References}

\usepackage{rotating}
\usepackage{fancyvrb}

\TheoremsNumberedThrough     % Preferred (Theorem 1, Lemma 1, Theorem 2)
\ECRepeatTheorems
\JOURNAL{Operations Research}

\EquationsNumberedThrough    % Default: (1), (2), ...
\MANUSCRIPTNO{OPRE-2024-02-0847}

\begin{document}
%\setcounter{page}{1}
%\DOI{10.1287/xxxx.0000.0000}%

\RUNAUTHOR{Banerjee, Hssaine, and Sinclair}

\RUNTITLE{Online Fair Allocation of Perishable Resources}

\TITLE{Online Fair Allocation of Perishable Resources}

\ARTICLEAUTHORS{%
\AUTHOR{Sid Banerjee}
%,\textsuperscript{a} Second Author,\textsuperscript{b} Third Author,\textsuperscript{c} Fourth Author,\textsuperscript{c}

\AFF{School of Operations Research and Information Engineering, Cornell University, \EMAIL{sbanerjee@cornell.edu}}
%\textsuperscript{b}School of Industrial Engineering, Good College, Collegeville, Maine 01234 \EMAIL{secauth@goodcoll.edu}; 
%\textsuperscript{c}Their Common Affiliation \EMAIL{thauth@anywhere.edu, fourauth@anywhere.edu}

%mirko.janc@informs.org
\AUTHOR{Chamsi Hssaine}

\AFF{Department of Data Sciences and Operations, Marshall School of Business, University of Southern California, \EMAIL{hssaine@usc.edu}}

\AUTHOR{Sean R. Sinclair}

\AFF{Department of Industrial Engineering and Management Sciences, Northwestern University, \EMAIL{sean.sinclair@northwestern.edu}}
}

\ABSTRACT{%
\input{final_arxiv_parts/abstract}
}

% \SUBJECTCLASS{Please confirm subject classifications.Fourier-Motzkin elimination; adjustable robust optimization; linear decision rules; redundant constraint identification.}

% \AREAOFREVIEW{Stochastic Models}.}

%{\CQ{Kindly provide the keywords.}}

\maketitle

%%% UNCOMMENT THIS TO GO BACK TO ORIGINAL
	% \input{original_parts/intro}
	% \input{original_parts/related_work}
	% \input{original_parts/preliminaries}
 %        \input{original_parts/perishing_discussion}
 % 	\input{original_parts/results}
	% \input{original_parts/experiments}
	% \input{original_parts/conclusion}

	\input{final_arxiv_parts/intro}
	\input{final_arxiv_parts/related_work}

	\input{final_arxiv_parts/preliminaries}
    \input{final_arxiv_parts/perishing_discussion}
 	\input{final_arxiv_parts/results}

    \input{final_arxiv_parts/special-cases}
	\input{final_arxiv_parts/experiments}
	\input{final_arxiv_parts/conclusion}

\ACKNOWLEDGMENT{The authors thank the area editor, associate editor, and two anonymous referees whose comments and suggestions were tremendously helpful in improving the paper.  Part of this work was done while Sean R. Sinclair was a postdoctoral associate at Massachusetts Institute of Technology.  Part of this work was done while Sean Sinclair and Sid Banerjee were visiting the Simons Institute for the Theory of Computing for the semester on Data-Driven Decision Processes. We gratefully acknowledge funding from the National Science Foundation under grants ECCS-1847393, DMS1839346, CCF-1948256, CNS-195599, and CNS-1955997, the Air Force Office of Scientific Research under grant FA9550-23-1-0068, and the Army Research Laboratory under grants W911NF-19-10217 and W911NF-17-1-0094.
}

\bibliographystyle{informs2014} % outcomment this and next line in Case 1
% \SingleSpacedXI
\bibliography{references} % if more than one, comma separated

% \smallskip

% \textbf{Sid Banerjee} is an associate professor in the school of operations research and information engineering at Cornell, working on topics at the intersection of data-driven decision-making, network algorithms and market design. He received his PhD in 2013 from the ECE Department at UT Austin, and was a postdoctoral researcher in the Social Algorithms Lab at Stanford from 2013-2015.

% \textbf{Chamsi Hssaine} is an assistant professor of data sciences and operations at the University of Southern California, Marshall School of Business. Her research interests are in the design and analysis of algorithms for data-driven decision-making, with applications in revenue and inventory management, as well as fair resource allocation. 

% \textbf{Sean R. Sinclair} is an assistant professor at Northwestern University in the department of industrial engineering and management sciences.  His research interests are in reinforcement learning and data-driven decision-making for operations systems.

\newpage

% \ECSwitch
% \ECHead{Electronic Companion for Online Fair Allocation of Perishable Resources}

\AtBeginEnvironment{APPENDICES}{%
  % Make appendix anchors unique for hyperref
  \renewcommand{\theHsection}{appendix.\Alph{section}}%
  \renewcommand{\theHsubsection}{\theHsection.\arabic{subsection}}%
  \renewcommand{\theHsubsubsection}{\theHsubsection.\arabic{subsubsection}}%
  \renewcommand{\theHfigure}{\theHsection.\arabic{figure}}%
  \renewcommand{\theHtable}{\theHsection.\arabic{table}}%
  \renewcommand{\theHequation}{\theHsection.\arabic{equation}}%
  \renewcommand{\theHtheorem}{\theHsection.\arabic{theorem}}%

  % --- cleveref aliases so refs say "Appendix" ---
  \crefalias{section}{appendix}%
  \crefalias{subsection}{appendix}%
  \crefalias{subsubsection}{appendix}%
}

\crefalias{section}{appendix}
\begin{APPENDICES}
    \OneAndAHalfSpacedXI
    \input{final_arxiv_parts/appendix_perishing}

    \input{final_arxiv_parts/special-cases-proofs}

\input{final_arxiv_parts/experiments_full}
    \input{final_arxiv_parts/useful_lemmas}
        % \bibliographystyleEC{informs2014}
    % \bibliographyEC{references}
\end{APPENDICES}

\end{document}

%% file: preamble.tex
\usepackage{amsmath,amsfonts,cancel}
\usepackage{thmtools}
\usepackage{adjustbox}
\usepackage[shortlabels]{enumitem}
\usepackage{mathtools}
\usepackage{dsfont}
\usepackage{tcolorbox}
\usepackage{booktabs}
\usepackage{bbm}
\usepackage{xspace}
\usepackage{subcaption}
\usepackage{tensor}
\usepackage{bm}
\usepackage{tikz}
\usetikzlibrary{patterns}
\usepackage{natbib}

\usepackage{placeins}
\usepackage[colorlinks=true,breaklinks=true,bookmarks=true,urlcolor=magenta,citecolor=magenta,linkcolor=magenta,bookmarksopen=false,draft=false]{hyperref}
\usepackage[capitalize]{cleveref}
\crefname{algocf}{Algorithm}{Algorithms}
\Crefname{algocf}{Algorithm}{Algorithms}

\newenvironment{rproofof}[1]{ \ifdefined\opre \begin{proof}{\it Proof of #1.}
\else \begin{proof}[Proof of #1] \fi }{ \ifdefined\opre 
\Halmos \end{proof} \else \end{proof} \fi  }

\newenvironment{rproof}{ \ifdefined\opre \proof{\it Proof.}
\else \begin{proof} \fi }{ \ifdefined\opre 
\Halmos \endproof \else \end{proof} \fi  }

\usepackage[suppress]{color-edits}
\addauthor{min}{red}

\newcommand{\Deff}{\Delta_{\text{\it efficiency}}}
\newcommand{\Denv}{\Delta_{\text{\it EF}}}

\newcommand{\Denvhind}{\textsc{Envy}}

\newcommand{\A}{\mathcal{A}}

\newcommand{\Ind}[1]{\mathds{1}\left\{ #1 \right\}}

\newcommand{\Exp}[1]{\mathbb E \left[ #1 \right]} % Variance
\newcommand{\Var}[1]{\mathrm{Var} \left[ #1 \right]}
\newcommand{\STD}[1]{\mathrm{Std} \left[ #1 \right]}

\renewcommand{\Pr}{\mathbb{P}}

\newcommand{\pfail}{\delta}

\newcommand{\E}{\mathcal{E}}

\newcommand{\Bupper}[1]{\overline{\mathcal{B}}_{#1}}

\newcommand{\HopeGuardrail}{\textsc{Vanilla-Guardrail}\xspace}
\usepackage{changepage}
\newcommand{\HopeGuardrailPerish}{\textsc{Perishing-Guardrail}\xspace}

\newcommand{\StaticXlower}{\textsc{Static-}$\underline{X}$\xspace}
\newcommand{\StaticProp}{\textsc{Static-}$\frac{B}{\overline{N}}$\xspace}
\newcommand{\Stockout}{\textsc{Stockout}\xspace}
\newcommand{\Spoilage}{\textsc{Spoilage}\xspace}

\newcommand{\sLCB}{\sigma^{\textsf{LCB}}}
\newcommand{\sCV}{\sigma^{\textsf{CV}}}
\newcommand{\sMean}{\sigma^{\textsf{Mean}}}

\newcommand{\conf}{\textsc{Conf}}

\newcommand{\Bav}{\beta_{avg}}

\newcommand{\PUA}{\textsc{PUA}}

\newcommand{\tstar}{t^*}
\newcommand{\dperish}{\mathcal{L}^{\textsf{perish}}}

\newcommand{\taub}[1]{\tau_b(#1)}
\newcommand{\conftaub}[1]{{\tau}_b(#1)}
\newcommand{\Pupper}{\overline{P}}

\newcommand{\confn}[2]{\sqrt{2({#2}-{#1})|\Theta|\rho_{max}^2 \log(2T^2 / \delta)}}

\newcommand{\stdpt}{\STD{P_{<t}}}

\newcommand{\Balg}[1]{\mathcal{B}^{alg}_{#1}}
\newcommand{\setB}[1]{\mathcal{B}_{#1}}

%% file: final_arxiv_parts/abstract.tex
We consider a practically motivated variant of the canonical online fair allocation problem: a decision-maker has a budget of {\em perishable} resources to allocate over a fixed number of rounds. Each round sees a random number of arrivals, and the decision-maker must commit to an allocation for these individuals before moving on to the next round. The goal is to construct a sequence of allocations that is {\it envy-free} and {\it efficient}. Our work makes two important contributions toward this problem: we first derive strong lower bounds on the optimal envy-efficiency trade-off, {demonstrating that} a decision-maker is fundamentally limited in what she can hope to achieve relative to the no-perishing setting; we then design an algorithm achieving these lower bounds which takes as input $(i)$ a prediction of the perishing order, and $(ii)$ a desired bound on envy. Given the remaining budget in each period, the algorithm uses forecasts of future demand and perishing to adaptively choose {from} one of two carefully constructed {guardrail quantities}. We demonstrate our algorithm's strong numerical performance --- and state-of-the-art, perishing-agnostic algorithms' inefficacy --- on simulations calibrated to a real-world dataset.

%% file: final_arxiv_parts/intro.tex
\section{Introduction}\label{sec:intro}

Despite a consistent decline in food insecurity in the United States over the past decade, 2022 saw a marked increase in individuals struggling to access enough food to fulfill basic needs. A recent report by the U.S. Department of Agriculture found that over 44 million individuals faced some form of hunger in 2022 --- 45\% more than the previous year~\citep{npr_2023}. Due in part to rising food prices and the rolling back of pandemic-era social security measures, this disturbing statistic has further underscored the important role of local food banks; for instance, the Feeding America network of food banks, food agencies, and local food programs distributed over 5.2 billion meals that same year~\citep{feeding_america}.

In distributing food throughout their operating horizon, food banks have two competing objectives: distributing {\it as much} food as possible to communities in need, and ensuring {\it equitable} access to donations. This tension has attracted much attention in the operations literature, with recent work characterizing the fundamental trade-offs between fairness and efficiency in sequential allocation problems~\citep{bertsimas2011price,donahue2020fairness,lien2014sequential,manshadi2021fair,sinclair2021sequential}. Understanding such trade-offs in theory is useful, as they allow a system designer to recognize and choose their desired operating point, balancing the loss in efficiency and equity. Despite the useful insights derived in prior work, to the best of our knowledge an important reality of food bank operations remains overlooked: waste due to {spoilage} of perishable goods. 

Perishable goods constitute a substantial portion of food donations. For instance, 35\% of the food distributed annually by the Food Bank of Central and Eastern North Carolina is perishable \citep{orgut2023}. Similarly, in 2024 alone, the Los Angeles Regional Food Bank distributed over 38 million pounds of fresh produce, representing 29\% of the total amount of food distributed throughout the year \citep{LAFoodBank2024Annual}. At the same time, the effective distribution of perishable goods has long been a significant challenge for these organizations. For example, in 2020 Feeding America reported that 3\% of the food received by its partner agencies was disposed of due to spoilage. Similar percentages were reported by other food rescue organizations in the United States, with Boston- and Rochester-based organizations respectively reporting that 2.4\% and 6\% of all received food was wasted due to spoilage \citep{epa2020wastedfood}.

These figures highlight a persistent operational gap between the supply of perishable goods and their successful distribution. The reason for this gap is two-fold. First, food banks lack control over the kinds of donated items, meaning that the supply of goods that are about to perish can far exceed the demand for these goods. \cite{ohls2002emergencyfood} highlighted this as a significant concern: in interviews with over 300 food banks, they found that around 20\% of food banks distribute less than 90\% of the food they receive, partly due to this issue. Second, even when this is not the case, true expiration dates are frequently unknown. Rather, product packaging commonly indicates ``best-by,'' ``use-by,'' or ``sell-by'' dates, after which products are still safe to consume~\citep{usda_food_product_dating}. So, while these dates may provide an informative signal as to when goods will perish (e.g., one week after the ``best-by'' date), the actual time at which this occurs is not known ex-ante. This problem is further exacerbated within the context of food banks, where individuals are encouraged to donate perishables that are {past} their ``best-by'' dates \citep{cityharvest2019donor,wqad_donate_food_2025}. Not only does this reduce the amount of time food banks have to allocate these items before they spoil, but it also creates significant uncertainty in the quality of the donations \citep{zou2025reducing}.

Intuitively, the equity-efficiency trade-off {should be} exacerbated in the presence of perishables: while equity requires a decision-maker to allocate {conservatively} to preserve resources for future demand \citep{sinclair2021sequential}, perishability starts a ``race against time.'' As goods perish due to a slow allocation rate, not only is efficiency harmed, but so may be equity, since spoilage runs the risk of a decision-maker running out of items, with nothing left to give out by the end of the operating horizon. Thus motivated, this paper seeks to answer the following questions:
\begin{center}{\it Do established equity-efficiency trade-offs in dynamic environments persist in the presence of perishable goods? If not, what limits do they impose on fair and efficient allocation? Can we design policies that perform well under these limits?}
\end{center}

\subsection{Our Contributions}

We consider a model in which a decision-maker has a fixed budget of items (also referred to as {\it goods}, or {\it resources}) to be allocated over $T$ discrete rounds. An a priori unknown number of individuals arrives in each period, each seeking a share of goods. Each agent is characterized by an observable type (drawn from a known, potentially time-varying distribution), associated with a linear utility function over the received allocation. Moreover, each unit of good has a stochastic perishing time (similarly drawn from a known distribution, independent of incoming arrivals), after which the good spoils and can no longer be allocated. The decision-maker allocates goods according to a fixed ordering (also referred to as perishing {\it prediction}, or {\it allocation schedule}). For instance, if each good has an associated ``best-by'' date, in practice the decision-maker may allocate the goods in increasing order of these dates. The goal is to find a policy that trades off between three ex-post metrics:
\begin{enumerate}[1.]
\item \emph{Hindsight Envy} ($\Denvhind$) -- Maximum difference in utility obtained by any two agents.
\item \emph{Counterfactual Envy} $(\Denv)$ -- Maximum difference between the utility obtained by any agent, and their utility under the static, proportional allocation (optimal in the no-perishing setting).
\item \emph{Inefficiency} $(\Deff)$ -- Amount of unallocated (including spoiled) goods at the end of the horizon.
\end{enumerate}

For this setting, we first characterize the fundamental limits that perishability places on any online algorithm. 
%\chdelete{We argue that --- contrary to the setting without perishable resources ---{counterfactual envy ceases to be a meaningful metric for a large class of perishing processes}.} 
To see the considerable impact that perishability can have on counterfactual envy, consider an extreme scenario in which all items perish at the end of the first round. Clearly, there is no hope of achieving low envy in such a setting, since future demand can never be satisfied.  Our first main contribution identifies a necessary and sufficient condition on the joint perishing and arrival distribution for low counterfactual envy to be an achievable desideratum (\cref{thm:offset-expiry-iff}). From a managerial perspective, this characterization --- which at a high level states that the cumulative perishing must lag behind cumulative arrivals --- can be viewed as guidance on the composition of perishable donations that the decision-maker should accept. It moreover underscores the importance of leveraging {\it joint} information over perishing and demand: if demand is {back-loaded}, perishing times should be concentrated late in the horizon; if most demand arrives early, however, early perishing times are acceptable.

For this class of processes, which we term {\it offset-expiring}, zero spoilage occurs if goods are allocated in the exact order in which they perish. Due to stochasticity in goods' perishing times, however, the decision-maker's allocation schedule may be misaligned with the realized perishing times of goods. For instance, while one good may have an earlier ``best-by'' date than another, the latter good may spoil earlier. Alternatively, if the decision-maker allocates goods in increasing order of expected perishing time, the variance of the perishing times may cause goods with later expected perishing times to perish earlier than others. We show that the inherent aggressiveness of the perishing process (i.e., in cases where the process is not offset-expiring) and inaccuracies in the allocation schedule pose an insurmountable barrier to any algorithm's performance by deriving lower bounds on any algorithm's envy and inefficiency. These lower bounds are parameterized by an unavoidable loss due to spoilage, also referred to as the {\it perishing-induced loss} and denoted by $\dperish$ (\Cref{def:sigma-loss}), implicitly characterized as the solution to a fixed-point equation (\cref{thm:lower_bound}). In contrast to the no-perishing setting, in which the only source of loss is {\it exogenous} uncertainty in demand, in our setting the loss due to spoilage is {\it endogenous}: it crucially depends on the rate at which the algorithm allocates items. This endogeneity poses significant challenges in the analysis of any adaptive algorithm; designing a tractable approach to analyzing ex-post spoilage is the main technical contribution of our work. Additionally, the lower bounds we derive give rise to the key insight that, contrary to the ``race against time'' intuition under which a decision-maker must {\it increase} the allocation rate to prevent {\it avoidable} spoilage, achieving low hindsight envy in the presence of {\it unavoidable} spoilage requires a decision-maker to potentially allocate significantly less than the proportional allocation. Hence, perishability throws a wrench into the well-studied envy-efficiency trade-off: while hindsight and counterfactual envy are aligned in the no-perishing setting, these two may be at odds when goods spoil, since only high counterfactual-envy solutions may yield low hindsight-envy. The tension between efficiency and equity is exacerbated for the same reason, relative to the classical setting.

In our final technical contribution, we leverage these insights to construct an adaptive threshold algorithm (\Cref{alg:brief-perishing}) that achieves these lower bounds (\cref{thm:upper_bound_perishing}). Our algorithm takes as input $(i)$ the fixed allocation schedule, $(ii)$ a desired upper bound on hindsight envy $L_T$, and $(iii)$ a high-probability parameter $\delta$. Given these inputs, it computes a high-probability lower bound on a budget-respecting zero-hindsight-envy solution, and an ``aggressive'' efficiency-improving allocation that is $L_T$ away. In each round the algorithm chooses which of the two quantities to allocate to each individual, cautiously doing so by constructing pessimistic forecasts of future arrivals and perishing. While this algorithm is similar in flavor to state-of-the-art algorithms for the no-perishing setting \citep{sinclair2021sequential}, the main challenge it contends with is forecasting endogenous future spoilage. We overcome this challenge by leveraging an {\it exogenous} ``slow'' consumption process that tractably decouples past allocations from future perishing. Crucially, we show that future spoilage under this slow process upper bounds our algorithm's future spoilage.

\begin{figure}[t]
\centering
\scalebox{.85}{
\tikzset{every picture/.style={line width=0.75pt}} %set default line width to 0.75pt    
\input{figures/intro_tikz/perishing_trade_off}}
\caption{Graphical representation of \cref{thm:upper_bound_perishing,thm:lower_bound}, illustrating the envy-efficiency trade-off ($\Deff$ vs. $\Denv$) achieved by \HopeGuardrailPerish (\Cref{alg:brief-perishing}) in the ``high-envy, low-perishing'' regime. Here, $\dperish$ is used to denote the unavoidable perishing-induced loss (\Cref{def:sigma-loss}).  
The vertical and horizontal dotted lines represent the impossibility results due to demand and perishing uncertainty, respectively. The region below the solid line represents the impossibility due to the envy-efficiency trade-off; the region above is the achievable region for \HopeGuardrailPerish.}
\label{fig:uncertainty_principle}
\end{figure}

Our bounds give rise to three salient regimes that dictate the extent to which a decision-maker is restricted in leveraging inequity to improve efficiency (and vice versa). These three regimes are a function of both $L_T$ and the unavoidable perishing-induced loss dictating our lower bounds, and at a high level are as follows:
\begin{enumerate}
\item {\it Low-envy}: There are no efficiency gains from deviating from the equitable solution, where all individuals receive the same amount.
\item {\it High-envy, high-perishing}: Similarly here, inefficiency is invariant to $L_T$.
\item {\it High-envy, low-perishing}: Inefficiency decreases as $1/L_T$, until reaching the unavoidable perishing-induced loss throughout the horizon, $T\dperish$. {This regime is represented pictorially in \Cref{fig:uncertainty_principle}.}
\end{enumerate}

We conclude our theoretical results by providing {\it explicit} characterizations of the unavoidable perishing-induced loss for important special cases (\Cref{prop:nec-cond-offset,prop:basic-conditions,thm:geometric_perishing,prop:gen_distribution_example}). These characterizations shed light on the salient drivers of this loss: $(i)$ how aggressively goods perish relative to the rate at which they are demanded, and $(ii)$ whether the allocation schedule approximately preserves the order induced by the true perishing times. Our results reveal that the {variability} of the perishing distribution plays a crucial role in both cases.

We complement our theoretical bounds by testing the practical performance of our algorithm on both synthetic and real-world datasets. Our experiments show that the unfairness required to achieve efficiency gains is order-wise larger than in settings without perishable resources. Additionally, they underscore the weakness of {\it perishing-agnostic} online algorithms. We observe that these latter algorithms are incapable of leveraging unfairness to improve efficiency across a variety of perishing regimes. In contrast to these, our algorithm's construction of a {\it perishing-aware} baseline allocation is necessary to mitigate stockouts across all --- rather than simply worst-case --- instances. These include instances where offset-expiry fails to hold with high probability, as is the case in the real-world dataset we use to calibrate our experiments \citep{keskin2022data}. Perhaps most surprisingly, despite our baseline allocation being significantly lower than that of algorithms that don't take into account perishability, our algorithm is {\it more efficient} than these more aggressive algorithms, in addition to being more fair. This observation contradicts the ``race against time'' intuition that aggressive allocations are necessarily more efficient than cautious ones. Finally, we numerically explore the question of selecting practical allocation schedules, with the goal of minimizing the unavoidable loss due to perishing. Our main managerial insight is that ordering items according to a high-probability lower bound on their perishing time is a robust, interpretable choice for decision-makers; it naturally takes into account items' expected perishing times, all the while hedging against the variability of the perishing process.

Finally, we note that while our work is motivated by food bank operations, perishability is also an important consideration in medical applications, such as blood bank operations~\citep{bar2017blood} and vaccine distribution~\citep{manshadi2021fair}. For these applications, our modeling framework is applicable to settings in which allocation occurs in rural or remote areas, where uncertainty around storage conditions (e.g., temperature and humidity) may affect the lifetime of these goods.

\subsubsection*{Paper organization.} We next survey the related literature. We present the model in \cref{sec:preliminary}, and formalize the limits of perishability in \cref{sec:construction_x_lower}. In \cref{sec:alg} we design and analyze an algorithm that achieves the derived lower bounds, illustrating the impact of the quality of the allocation schedule and the aggressiveness of the perishing process for important special cases in \Cref{sec:special-cases}. We conclude with numerical experiments in \cref{sec:experiments}. 

%% file: figures/intro_tikz/perishing_trade_off.tex
\tikzset{every picture/.style={line width=0.75pt}} %set default line width to 0.75pt        

\begin{tikzpicture}[x=0.75pt,y=0.75pt,yscale=-1,xscale=1]

%Shape: Rectangle [id:dp559707459456056] 
\draw  [draw opacity=0][fill={rgb, 255:red, 184; green, 233; blue, 134 }  ,fill opacity=1 ] (220,2.68) -- (489.4,2.68) -- (489.4,100.4) -- (220,100.4) -- cycle ;
%Shape: Rectangle [id:dp5555924507039581] 
\draw  [draw opacity=0][fill={rgb, 255:red, 240; green, 180; blue, 190 }  ,fill opacity=1 ] (220,2.68) -- (251.68,2.68) -- (251.68,260.06) -- (220,260.06) -- cycle ;
%Shape: Rectangle [id:dp29237468353515483] 
\draw  [draw opacity=0][fill={rgb, 255:red, 240; green, 180; blue, 190 }  ,fill opacity=1 ] (220,2.68) -- (488.5,2.68) -- (488.5,260.06) -- (220,260.06) -- cycle ;
%Straight Lines [id:da42754641059606724] 
\draw  [dash pattern={on 4.5pt off 4.5pt}]  (252.63,100.58) -- (252.63,260.06) ;
%Shape: Arc [id:dp15076133545120962] 
\draw  [draw opacity=0][fill={rgb, 255:red, 184; green, 233; blue, 134 }  ,fill opacity=1 ] (490.57,246.11) .. controls (470.61,244.32) and (449.99,240.68) .. (429.16,235.04) .. controls (339.59,210.81) and (271.88,156.67) .. (252.63,100.58) -- (464.42,104.72) -- cycle ; \draw   (490.57,246.11) .. controls (470.61,244.32) and (449.99,240.68) .. (429.16,235.04) .. controls (339.59,210.81) and (271.88,156.67) .. (252.63,100.58) ;  
%Shape: Rectangle [id:dp14588434123826655] 
\draw  [draw opacity=0][fill={rgb, 255:red, 184; green, 233; blue, 134 }  ,fill opacity=1 ] (252.63,2.45) -- (489.43,2.45) -- (489.43,100.58) -- (252.63,100.58) -- cycle ;
%Shape: Rectangle [id:dp5990218380932004] 
\draw  [draw opacity=0][fill={rgb, 255:red, 184; green, 233; blue, 134 }  ,fill opacity=1 ] (461.4,37.4) -- (489.4,37.4) -- (489.4,235.6) -- (461.4,235.6) -- cycle ;
%Shape: Rectangle [id:dp1990239878816188] 
\draw  [draw opacity=0][fill={rgb, 255:red, 184; green, 233; blue, 134 }  ,fill opacity=1 ] (258.6,87.6) -- (489.6,87.6) -- (489.6,106.6) -- (258.6,106.6) -- cycle ;
%Straight Lines [id:da09977604831446596] 
\draw    (252.63,2.68) -- (252.63,100.58) ;
%Straight Lines [id:da8967771510446751] 
\draw  [dash pattern={on 4.5pt off 4.5pt}]  (490.57,246.11) -- (254,246.11) ;

% Text Node
\draw (194.1,81) node [anchor=north west][inner sep=0.75pt]  [font=\large,rotate=-269.93] [align=left] {$\Deff$};
% Text Node
\draw (470,262.06) node [anchor=north west][inner sep=0.75pt]  [font=\large] [align=left] {$\Denv$};
% Text Node
\draw (260.05,208.73) node [anchor=north west][inner sep=0.75pt]  [font=\footnotesize] [align=left] {Envy-Efficiency \\ Tradeoff};
% Text Node
\draw (246.65,263.46) node [anchor=north west][inner sep=0.75pt]    {$T^{-1/2}$};
% Text Node
\draw (310.22,100.89) node [anchor=north west][inner sep=0.75pt]   [align=left] {Achievable with \\ \HopeGuardrailPerish};
% Text Node
\draw (226,252.06) node [anchor=north west][inner sep=0.75pt]  [font=\footnotesize,rotate=-270] [align=left] {Impossible due to Demand Uncertainty};
% Text Node
\draw (256,249.11) node [anchor=north west][inner sep=0.75pt]  [font=\footnotesize] [align=left] {Impossible due to Perishing Uncertainty};
% Text Node
\draw (491.4,230) node [anchor=north west][inner sep=0.75pt]  [font=\large] [align=left] {$T \dperish$};

\end{tikzpicture}

%% file: final_arxiv_parts/related_work.tex
\subsection{Related Work}
\label{sec:related_work}

Fairness in resource allocation has a long history in the economics and computation literature, beginning with Varian's seminal work \citep{varian1973equity,varian1976two}. More recently, there has been work studying fairness in operations, including {resource-level fairness through load balancing}~\citep{chen2022fair,hssaine2024target,balseiro2025regularized,freund2023end}, pricing~\citep{cohen2022price,den2022waste}, incentive design~\citep{freund2021fair,hssaine2025learning}, algorithmic hiring and online selection~\citep{salem2023secretary,hu2025diversity}, and societal systems more generally \citep{gupta2021individual,lyu2025price}.  We highlight the most closely related works below, especially as they relate to {\it online} fair allocation; see \cite{aleksandrov2019online} for a survey.

\paragraph{Fair allocation without perishable resources.} There exists a long line of work in which {non-perishable} resources become available to the decision-maker online, whereas agents are fixed \citep{benade2018make,aleksandrov2015online,mattei2017mechanisms,mattei2018fairness, aleksandrov2019monotone,gorokh2020fair,bansal2020online,bogomolnaia2022fair,ijcai2019-49,aziz2016control,zeng2019fairness}. 
These models lie in contrast to the one we consider, wherein resources are fixed and {\it individuals} arrive online. Papers that consider this latter setting include \citet{kalinowski2013social}, which studies welfare maximization with {indivisible} goods; \citet{gerding2019fair} considers a scheduling setting wherein agents have fixed and known arrival and departure times, as well as demand for the resource; \citet{hassanzadeh2023sequential} allows individuals to arrive in multiple timesteps.  
A series of papers also consider the problem of fair allocation with minimal disruptions relative to previous allocations, as measured by a {\it fairness ratio}, a competitive ratio analog of counterfactual envy in our setting \citep{friedman_2017,cole_2013,10.1145/2764468.2764495}. Other works design algorithms with attractive competitive ratios with respect to Nash Social Welfare~\citep{azar2010allocate,gorokh2020fair} and the max-min objective~\citep{lien2014sequential,manshadi2021fair}. \citet{ma2025forward} recently showed how forward-backwards contention resolution schemes can be leveraged to obtain competitive ratios for online rationing in non-perishable settings.

The above papers situate themselves within the {adversarial}, or worst-case, tradition. {A separate line of work considers fair resource allocation in stochastic settings~\citep{donahue2020fairness,elzayn2019fair}, as we do. The algorithms developed in these papers, however, are {\it non-adaptive}: they decide on the entire allocation upfront, {\it before} observing any of the realized demand. In contrast, we consider a model where the decision-maker makes the allocation decision in each round {\it after} observing the number of arrivals. \citet{freeman2017fair} considers a problem in which agents' utilities are realized from an unknown distribution, and the budget resets in each round. They present algorithms for Nash social welfare maximization and discuss some of their properties. 
Our work is most closely related to (and indeed, builds upon) \citet{sinclair2021sequential}, which first introduced the envy-efficiency tradeoff we are interested in. Our work considers the more challenging setting of perishable goods, which none of the aforementioned works consider.

Finally, we note the existence of a rich line of work that takes a mechanism design approach to fair division, dating back to the 1940s \citep{steinhaus1948problem}. We refer the reader to \citet{moulin2004fair}, Chapters 11-13 of \citet{brandt2016handbook}, and Chapter 11 of \citet{kroer2026games} for textbook treatments of the topic. Generally speaking, the main objective of these works is the design of {\it nonmonetary} mechanisms that fairly allocate a set of resources amongst {\it strategic} agents with private preferences for these goods. Much of the work has focused on one-shot (i.e., static) settings, for which there exist strong impossibility results precluding the existence of fair, efficient, {\it and} incentive compatible mechanisms \citep{budish2011combinatorial}. Notable works that derive {\it approximately} fair and incentive-compatible mechanisms for these settings include \citet{budish2011combinatorial} and \citet{branzei2017nash}. In dynamic settings, models of fair allocation under incentives involve the {\it same} set of agents competing for resources across time, which differs from the arrival model we consider here. Recent works that consider this setting include \citet{onyeze2025allocating}, who study equilibrium existence, and \citet{fikioris2024incentives}, who propose a generalization of the popular dynamic max-min fair allocation mechanism to achieve strong fairness and approximate incentive compatibility guarantees.

We conclude by highlighting the closely related work of \citet{prendergast2022allocation}, which empirically studies the virtual currency-based auction mechanism used by Feeding America to allocate food to its members. This work finds that, compared to a previous queuing mechanism under which food banks were being allocated an equal amount of random food through time, the auction mechanism allows food banks to express preferences across different types of food. As it relates to our work, it would be interesting and practically relevant to study how $(i)$ the {timing} of allocations, and $(ii)$ varying preferences across resources induce strategic behavior. We leave this for future work.

\paragraph{Dynamic allocation of perishable resources.} Though dynamic resource allocation of perishable goods has a long history in the operations research literature (see, e.g.,  \citet{nahmias2011perishable} for a comprehensive survey of earlier literature), to the best of our knowledge, the question of {\it fairly} allocating perishable goods has attracted relatively little attention. We highlight the few relevant papers below.
\citet{perry99} and \citet{hanukov2020service} analyze FIFO-style policies for efficiency maximization in inventory models with Poisson demand and deterministic or Poisson perishing times.
 Motivated by the problem of electric vehicle charging, \citet{gerding2019fair} consider an online scheduling problem where agents arrive and compete for a perishable resource that spoils at the end of every period, and as a result must be allocated at every time step. They consider a range of objectives, including: maximum amount of resources allocated, maximum number of satisfied agents, as well as envy-freeness. \citet{bateni2022fair} similarly consider a setting wherein an arriving stream of goods perish immediately. {Recent empirical work by \citet{orgut2023} considers the problem of a food bank equitably and efficiently allocating perishable goods under {complete information}. Their case study on data from a partnering food bank numerically validates our theoretical results: in low-budget settings, there is little or no benefit to increasing inequity to counteract the inefficiency due (in part) to spoilage.} 
In contrast to these latter papers, the model we consider locates itself within the smaller category of inventory models in which products have {\it random} lifetimes. The majority of these assume that items have exponentially distributed or geometric lifetimes \citep{bakker2012review}.

%% file: final_arxiv_parts/preliminaries.tex
\section{Preliminaries}
\label{sec:preliminary}

We consider a decision-maker who, over $T$ rounds, must allocate $B$ divisible units (also referred to as {\it items}) of a single type of resource among a population of individuals. Let $\mathcal{B}$ denote the set of these $B$ units.

\subsection{Basic Setup}

\paragraph{Demand model.}
At the start of each round $t \in [T]$, a random number of individuals arrives, each requesting a share of units. Each individual is characterized by her {\it type} $\theta \in \Theta$, with $|\Theta| < \infty$. Each type $\theta$ is associated with a known utility function $u_\theta(x) = w_\theta \cdot x$ for a given allocation $x \in \mathbb{R}_+$ of the resource, with $w_\theta > 0$.  We let $N_{t,\theta}$ denote the number of type $\theta$ arrivals in round $t$; $N_{t,\theta}$ is drawn independently from a known distribution, with $N_{t,\theta} \geq 1$ almost surely for all $t\in [T], \theta \in \Theta$.  This latter assumption is for ease of exposition; our results continue to hold (up to constants) as long as $\Pr(N_{t, \theta} = 0)$ does not scale with $T$. For ease of notation we define $N_t = \sum_{\theta \in \Theta} N_{t, \theta}$ and $N = \sum_{t \in [T], \theta \in \Theta} N_{t, \theta}$. We assume $\mathbb{E}[N] = \Theta(T)$, and define $\Bav = B/\mathbb{E}[N]$ to be the average number of units per individual, with $\Bav = \Theta(1)$.

\paragraph{Perishing model.} 
Each unit of resource $b \in \mathcal{B}$ is associated with a perishing time $T_b \in \mathbb{N}^+$ drawn from a {known} distribution. Items' perishing times are independent of one another and of the arrival process, and perishing occurs at the {end} of each round, after items have been allocated to individuals.  For $t \in [T]$, we let $P_t = \sum_{b\in\mathcal{B}}\mathds{1}\left\{T_b = t\right\}$ denote the number of units of resource perishing in period $t$. 

The decision-maker has access to a {\em predicted ordering} according to which items perish; we will often refer to this ordering as the {\it allocation schedule}. For example, the decision-maker may order items in increasing order of their expected perishing times, or according to a lower confidence bound on their perishing times. As alluded to in the introduction, within the context of food bank operations, the decision-maker may have access to ``best-by'' dates specified on food packaging and choose to allocate items with earlier ``best-by'' dates first. We use  $\sigma: \mathcal{B} \rightarrow [B]$ to denote this ordering, i.e., $\sigma(b)$ is the {rank} of $b$ in this ordering. For $b \in [B]$, $\sigma^{-1}(b)$ is used to denote the {\it identity} of the $b$th ranked item in $\sigma$, with $\sigma^{-1}(1)$ being the item that comes first in the allocation schedule. 

\begin{example}\label{ex:sigma}
Consider an instance for which $\mathcal{B} = \{\texttt{a},\texttt{b},\texttt{c}\}$, and let $\sigma$ be such that \mbox{$\sigma(\texttt{a}) = 2$}, $\sigma(\texttt{b}) = 3$, and $\sigma(\texttt{c}) = 1$. Then, $\sigma^{-1}(1) = \texttt{c}$, $\sigma^{-1}(2) = \texttt{a}$, and $\sigma^{-1}(3) = \texttt{b}$.  
\end{example}

\begin{remark}
In this paper we restrict our attention to static (rather than time-varying and sample path-dependent) allocation schedules, given their practical relevance to the motivating real-world applications described in \cref{sec:intro}. We leave the nonstationary extension to future work.
\end{remark}

\paragraph{Additional notation.}  For any time-dependent quantity $Y_t$, we define $Y_{\leq t} = \sum_{t' \leq t} Y_{t'}$, $Y_{\geq t} = \sum_{t' \geq t}Y_{t'}$, along with their strict analogs. We let $w_{max} = \max_\theta w_\theta$, $\sigma_{t, \theta}^2 = \text{Var}(N_{t, \theta}) < \infty$, and assume $\rho_{t, \theta} = |N_{t, \theta} - \mathbb{E}[N_{t, \theta}]| < \infty$ almost surely. Finally, let $\mu_{\max} = \max_{t} \Exp{N_{t}}, \sigma^2_{\min} = \min_{t, \theta} \sigma^2_{t, \theta}, \sigma^2_{\max} = \max_{t,\theta} \sigma^2_{t, \theta}$, and $\rho_{\max} = \max_{t,\theta} \rho_{t,\theta}$. We use $\lesssim$ and $\gtrsim$ to denote the fact that inequalities hold up to polynomial factors of $\Bav, |\Theta|, w_{max}, \mu_{\max}, \sigma^2_{\min}, \sigma^2_{\max}$, $\log T$, and $\log(1 / \delta)$.

\subsection{Notions of Fairness and Efficiency}

The goal is to design a {\it fair} and {\it efficient} online algorithm that determines the amount to allocate to all $N_{t,\theta}$ individuals in each round $t$, for all $\theta\in\Theta$, given the remaining budget in each round. We assume this amount is allocated uniformly across all $N_{t, \theta}$ individuals of type $\theta$. We use $X_{t, \theta}^{alg} \in \mathbb{R}_+$ to denote the per-individual amount distributed in period $t$, with $X^{alg} = (X_{t, \theta}^{alg})_{t \in [T]}$. Finally, we let $B_t^{alg}$ denote the remaining budget at the beginning of period $t$, defined inductively as follows:
\begin{align}
\begin{cases}
B_1^{alg} &= B,  \\
B_{t+1}^{alg}
&= B_t^{alg}
   - \sum_{\theta \in \Theta} N_{t,\theta} X_{t,\theta}^{alg}- \PUA_t^{alg},
   \qquad t \in [T].
\end{cases}
\end{align}
where $\PUA_t^{alg}$ denotes the quantity
of unallocated items that perish at the end of round $t$ (i.e., the {\it perished and unallocated} items). We moreover let $\mathcal{B}_t^{alg}$ be the {\it set} of items remaining at the beginning of period $t$.

Our notions of online fairness and efficiency are motivated by the offline notion of {\it Varian Fairness} \citep{varian1973equity}, and are the same as those considered in past works \citep{sinclair2021sequential}.

\begin{definition}[Counterfactual Envy, Hindsight Envy, and Efficiency]
\label{def:distance}
Given budget $B$, realized demands $(N_{t, \theta})_{t \in [T], \theta \in \Theta}$, perishing times $(T_b)_{b \in [B]}$, and allocation schedule $\sigma$, for any online allocation defined by $X^{alg}$ we define:
\begin{itemize}
    \item\emph{Counterfactual Envy}: %The counterfactual distance of  $X^{alg}$ to envy-freeness as 
    \begin{align}
    \label{eq:cenv}
    \Denv\triangleq \max_{t \in [T], \theta \in \Theta} \left|w_\theta \left(X_{t, \theta}^{alg} -\frac{B}{N}\right)\right|.\end{align}
    \item \emph{Hindsight Envy}: %The hindsight distance of $X^{alg}$ to envy-freeness as 
    \begin{align}
    \label{eq:env_hind}
        \Denvhind \triangleq \max_{t,t' \in [T]^2, \theta, \theta' \in \Theta^2} w_\theta (X_{t', \theta'}^{alg} - X_{t, \theta}^{alg}).
    \end{align}
    \item\emph{Inefficiency}: %The distance to efficiency as
    \begin{align}
    \label{eq:efficiency}
    \Deff \triangleq B - \sum_{t \in [T], \theta \in \Theta} N_{t, \theta} X_{t,\theta}^{alg}.
    \end{align}
\end{itemize}

\end{definition}

In the offline setting {\it without} perishable goods, \citet{varian1973equity} established that $X^{opt}_{t,\theta} = B/N$ (referred to as the {\it proportional allocation}) is the optimal fair and efficient per-individual allocation. Hence, counterfactual envy $\Delta_{EF}$ can be interpreted as a form of regret with respect to this strong no-perishing benchmark, and can be used to characterize the impact of perishability on our algorithm's performance. Hindsight envy, on the other hand, measures how differently the online algorithm treats any two individuals across time. Finally, the efficiency of the online algorithm, $\Deff$, measures how {wasteful} the algorithm was in hindsight. This could happen in two ways: either through spoilage, or because the decision-maker allocated too conservatively throughout the horizon, thus leaving a large number of unspoiled goods unallocated by $T$.

Even in simple settings without perishability, it is known that these metrics are at odds with each other in online settings. To see this, consider the following two scenarios. On the one hand, an algorithm can trivially achieve a hindsight envy of zero by allocating nothing to individuals in any round; this, however, would result in both high {\it counterfactual} envy, in addition to maximal inefficiency. On the other hand, a distance to efficiency of zero can trivially be satisfied by exhausting the budget in the first round, at a cost of maximal hindsight envy as individuals arriving at later rounds receive nothing. {\citet{sinclair2021sequential} formalized this tension for the additive utility setting without perishable resources via the following lower bounds.}
%, which we re-state here for completeness.
\begin{theorem}[Theorems 1 and 2, \citet{sinclair2021sequential}]
\label{thm:lower_bounds}
Under any arrival distribution satisfying {mild regularity conditions}, there exists a problem instance without perishing, such that any algorithm must incur $\Denv \gtrsim \frac{1}{\sqrt{T}}$. Moreover, any algorithm that achieves $\Denv = L_T = o(1)$ or $\Denvhind = L_T = o(1)$ must also incur $\Deff \gtrsim \min\{\sqrt{T}, 1/L_T\}.$
\end{theorem}

Since settings without perishable resources are a special case of our setting (e.g., a perishing process with $T_b > T$ a.s., for all $b\in \mathcal{B}$), this lower bound holds in our case; the goal then is to design algorithms that achieve this lower bound with high-probability. However, as we will see in \cref{sec:construction_x_lower}, perishing stochasticity is fundamentally distinct from, and more challenging than, demand stochasticity. This difference is particularly salient in regards to the envy-efficiency trade-off.  

\subsection{Discussion of Modeling Assumptions}\label{ssec:discuss-model}
We conclude the section with a discussion of our modeling assumptions.

\paragraph{Age-dependent quality considerations.} While our model does not explicitly account for freshness considerations (i.e., allocating an item that perishes tomorrow is considered to be the same ``good'' outcome as allocating an item that perishes in 10 days), one natural way to take such considerations into account would be to define a pseudo-perishing time $\widetilde{T}_b = T_b-a_b$, where $a_b > 0$ is an item-dependent constant that guarantees a minimum amount of time that an item lasts in the household once it is allocated. (In practice, the decision-maker may want to set $a_b$ in a category-specific way, such that the allocated item meets a minimum quality threshold at the time of allocation. For instance, milk can last 5-7 days after its ``sell-by'' date if stored properly in a refrigerator \citep{usdaire_milk_sit_out_2025}. In this case, it may be desirable to set $a_b = 3$ in order to ensure that a household can consume the milk for at least 3 days.) Defining all metrics with respect to these pseudo-perishing times, our results would continue to hold, thereby guaranteeing some notion of freshness in our end allocation.\footnote{If $\sigma$ is chosen according to the distribution of the pseudo-perishing time $\widetilde{T}_b$, the choice of $a_b$ may also have an impact on $\sigma$. For instance, if items are allocated in increasing order of $\mathbb{E}[\widetilde{T}_b]$, items for which $a_b$ is larger (i.e., those for which there is a quality premium) will be allocated earlier. If the decision-maker allocates items according to the {\it variance} of $\widetilde{T}_b$, on the other hand, given that $a_b$ is a constant, this choice would have no impact on $\sigma$.}

Alternatively, one could encode freshness considerations by modeling time-dependent utilities for each item. While this would be an interesting direction to explore theoretically, it is less practically relevant with respect to our motivating application of food banks, given that the {\it quantity} of goods allocated --- as opposed to an estimate of utility-based social welfare ---- remains the most important metric reported by these organizations~\citep{prendergast2022allocation,orgut2016achieving}.

\paragraph{On the allocation schedule $\sigma$.} In this work we assume the allocation schedule $\sigma$ is an input to the model, as opposed to a lever available to the decision-maker. While treating $\sigma$ as a decision variable would improve both the envy and inefficiency of our algorithm, the reason for this decoupling is two-fold. From a practical perspective, one can imagine that decision-makers have a preference for using interpretable allocation schedules, such as according to expected perishing time or ``best-by'' date. In contrast, having an allocation schedule depend on the specifics of the guardrail algorithm would be more challenging to communicate to stakeholders. This would be problematic in instances where spoilage occurred, and a unit with a later ``best-by'' date or expected perishing time was allocated {\it before} the unit that spoiled. Such bad events are harder to justify when $\sigma$ is an output of the decision-maker's optimization problem. By allowing $\sigma$ to be an input to the model, our framework can instead be used by the decision-maker in practice to test the performance of various interpretable allocation schedules, and heuristically determine the one that yields the strongest performance, as we do in our numerical experiments.

We moreover posit that treating $\sigma$ as a decision variable would be analytically challenging, as well as a computationally intractable task. {Namely, as we will see in \Cref{sec:construction_x_lower},} optimizing over $\sigma$ reduces to an optimization over $N!$ orderings, which cannot be done efficiently. We however make some progress toward characterizing ``high-quality'' allocation schedules by $(i)$ providing sufficient conditions on the relationship between $\sigma$ and the perishing distribution for our algorithm's spoilage to vanish asymptotically (\Cref{prop:gen_distribution_example}), and $(ii)$ numerically showing that allocating units in increasing order of a lower confidence bound on their perishing time is a robust choice for $\sigma$ (\Cref{sec:experiments_other}). We leave the derivation of structural results surrounding the optimal allocation schedule (which we conjecture will, in general, be highly distribution-dependent), as well as the design of approximately optimal allocation schedules, as interesting open questions. 

Finally, we note that the {\it ranking} of units is a sufficient statistic for the model we consider in this paper (as opposed to, say, granular forecasts of their perishing times). To see this, note that given a sequence of allocation quantities, the only question that remains for the decision-maker is {\it which} units she should allocate in each period. This reduces to the task of finding a ranking of units that minimizes spoilage, given the allocation quantities. One can easily see this when perishing times are known and deterministic; in this case, given an allocation quantity in each period, knowledge of the actual perishing times provides no additional advantage over the ranking that is induced by these true perishing times. This similarly holds for any forecast of units' perishing times treated as ground truth.

%% file: final_arxiv_parts/perishing_discussion.tex
\section{Fundamental Limits of Perishability: Lower Bounds}
\label{sec:construction_x_lower}

In the presence of perishable resources, a decision-maker must contend with two obstacles: $(i)$ the aggressiveness of the perishing process, and $(ii)$ errors in the perishing prediction $\sigma$. In this section we formalize the impact of these two challenges. Namely, we identify classes of perishing processes for which there is no hope of achieving the optimal fair allocation, $B/N$, and derive lower bounds on any algorithm's performance as a function of a salient quantity that captures both the aggressiveness of perishing, as well as the quality of the prediction $\sigma$. In the remainder of the section, we say that an online algorithm is {\it feasible} over a sample path if it does not run out of budget. Note that, if an algorithm is infeasible over a sample path, it necessarily achieves $\Denv = \Theta(1)$.

\subsection{On the Aggressiveness of the Perishing Process}\label{sec:offset-expiry}

We first argue that the proportional allocation $X^{opt}_{t,\theta} = B/N$ is unachievable unless one places restrictions on the rate at which items perish, even under full information over perishing times and demand.
%This would then imply optimality of $B/N$, since $B/N$ is optimal for the setting without perishable resources, an upper bound on our setting. 

To see this, consider an instance where all items perish at the end of the first round. There is no hope of achieving low envy in this setting, since there are no items left for arrivals from $t=2$ onwards. The following result establishes that the {\it only} perishing time realizations for which $B/N$ is achievable are ones in which the fraction of perished items ``lags behind'' the proportion of arrivals in each period. We formalize this via the notion of {\it offset-expiry}, defined below.

\begin{definition}[Offset-expiring process]\label{def:offset-expiry-sample-path}
    A perishing process $(T_b)_{b \in \mathcal{B}}$ is {\it offset-expiring} if:
    \begin{align}\label{eq:offset-expiry}
    \frac{P_{< t}}{B} \leq \frac{N_{< t}}{N} \quad \forall \ t \geq 2.
    \end{align}
\end{definition}
We provide an example below.
\begin{example}[Continued from \Cref{ex:sigma}]\label{ex:offset}
Consider a setting where  $T = B = 3$, with $N_1 = 1$, $N_2 = 2$, and $N_3 = 1$ (and $N = 4$ as a result). Suppose moreover $|\Theta| = 1$, with $w_{\theta} = 1$. For this sample path, we have:
\begin{align}
\frac{N_{<t}}{N} = \begin{cases}
\frac14 \quad &\text{for } t = 2\\
\frac34 &\text{for } t =3.
\end{cases}
\end{align}
Applying \Cref{eq:offset-expiry} to this sample path, a process is offset-expiring if and only if:
\begin{align}
P_{<t} \leq \begin{cases}
0 \quad &\text{for } t = 2\\
2 &\text{for } t =3.
\end{cases}
\end{align}
In words, for any offset-expiring process in this example, we must have $T_b \geq 2$ {almost surely} for all $b \in \mathcal{B}$, with no more than two units perishing at $t = 2$ in the worst case.
\end{example}

\cref{thm:offset-expiry-iff} establishes that {offset-expiry} exactly captures the trajectories whereby $B/N$ is a feasible allocation when units are allocated in increasing order of $T_b$. We refer to this latter ordering as the {\it hindsight optimal} ordering. We defer its proof to Appendix \ref{apx:proof-of-offset-expiry}.
\begin{theorem}\label{thm:offset-expiry-iff}
$X_{t,\theta} = B/N$ for all $t,\theta$ is feasible under the {hindsight optimal ordering} if and only if the perishing process is offset-expiring.
\end{theorem}

From a managerial perspective, the offset-expiry condition provides decision-makers with guidance on the selection of perishable goods to stock at the beginning of the horizon. Within the context of food banks, for instance, it highlights that rejecting perishable donations outright (as some food banks do) is too severe a policy; rather, the reasonable rule of thumb that says ``don't accept food if it will spoil faster than the rate at which you allocate it'' not only suffices, but is the {\it only} correct rule of thumb. It moreover underscores the importance of {\it jointly} considering the demand and perishing distributions in this selection. 
We provide necessary and sufficient conditions for offset-expiry to hold with high probability in \cref{sec:special-cases}.

\subsection{Unavoidable Loss Due to Perishing}
\label{ssec:unavoidable-loss}

The previous section established that, even under full information about the perishing process, there exist aggressive perishing distributions for which the proportional allocation $B/N$ is not achievable.
The following example illustrates how both the inherent aggressiveness of the exogenous perishing process, as well as the quality of the allocation schedule $\sigma$, play an important role in determining what an algorithm can achieve.

\begin{example}[Continued from \Cref{ex:sigma,ex:offset}]\label{ex:sigma-is-important}
Consider the following perishing distributions:
\begin{align}\label{eq:example-dists}
T_{\texttt{a}} =
\begin{cases}
1, & \text{w.p. } 1/2,\\
3, & \text{w.p. } 1/2,
\end{cases} \qquad 
T_{\texttt{b}} =
\begin{cases}
2, & \text{w.p. } 1/2,\\
4, & \text{w.p. } 1/2,
\end{cases} \qquad 
T_{\texttt{c}} =
\begin{cases}
1, & \text{w.p. } 1/2,\\
2, & \text{w.p. } 1/2.
\end{cases}
\end{align}
Recall, $N_1 = 1$, $N_2 = 2$, and $N_3 = 1$; moreover, $\sigma(\texttt{a}) = 2$, $\sigma(\texttt{b}) = 3$, and $\sigma(\texttt{c}) = 1$.  Note that $\sigma$ orders the units in increasing order of expected perishing time for this instance. 

By \Cref{ex:offset}, this process is not offset-expiring for any sample path in which $T_{\texttt{a}} = 1$ or $T_{\texttt{c}} = 1$. By \Cref{thm:offset-expiry-iff}, then, the proportional allocation of $B/N = 3/4$ is only feasible on sample paths for which \mbox{$T_{\texttt{a}} = 3$} and $T_{\texttt{c}} = 2$, under the hindsight optimal ordering. Even over such sample paths, however, the possibility remains that $B/N$ is infeasible under $\sigma$. To see this, consider the case where $T_{\texttt{a}} = 3$, $T_{\texttt{c}} = 2$, and $T_{\texttt{b}} = 2$. By allocating $3/4$ in each period according to ordering $\sigma$, $2.25$ units will have been allocated by the end of period $2$ (i.e., the entirety of $\texttt{c}, \texttt{a}$, and $0.25$ units of $\texttt{b}$). However, the remainder of $\texttt{b}$ perishes at the end of period $2$, leaving nothing for the individual who arrives in period $3$.
\end{example}

Motivated by this example, our key insight is that, for any static allocation $X$, there exists a high-probability worst-case loss, denoted by $\overline{\Delta}(X)$, that any algorithm incurs, due to both errors in $\sigma$, as well as the inherent aggressiveness of the perishing process. As a result, rather than having a budget of $B$ items, the algorithm has an {\it effective} budget of $B-\overline{\Delta}(X)$ items. Under this effective budget, any feasible static allocation must set $X$ such that $\overline{N}X \leq {B-\overline{\Delta}(X)}$, where $\overline{N}$ is a high-probability upper bound on $N$. Noting that $X = 0$ is always a feasible solution to this inequality, a natural choice is to set:
\begin{align}\label{eq:xlower}
\underline{X} = \sup\left\{X \ \mid \ X \leq \frac{B-\overline{\Delta}(X)}{\overline{N}}\right\},
\end{align}
if this supremum is achieved. When $\overline{\Delta}(X) = 0$ for all $X$ (e.g., either items don't perish or {the process is offset-expiring and predictions are perfect}), $\underline{X} = B/\overline{N}$, and we recover the conservative allocation of the no-perishing setting \citep{sinclair2021sequential}. 

By \Cref{eq:xlower}, given the worst-case perishing loss $\overline{\Delta}(X)$ for any allocation $X$, one can compute $\underline{X}$ via line search. Obtaining tight bounds on $\overline{\Delta}(X)$, however, is the challenging piece. To see this, note that for any algorithm, the quantity of goods that perished by the end of the horizon is:
\begin{align}\label{eq:perishing-expression}
\sum_{t\in[T]}\sum_{b \in \mathcal{B}_t^{alg}} (B_{t}^{alg}(b)-X_t^{alg}(b))^+\cdot \mathds{1}\{T_b = t\},
\end{align}
where $B_{t}^{alg}(b)$ is the {quantity} of item $b$ remaining at the beginning of period $t$, and $X_t^{alg}(b)$ is the quantity of item $b$ given out in period $t$. Since $X_t^{alg}(b)$ depends on the perishing realizations of previous rounds, computing this quantity requires the ability to simulate sufficiently many replications of the static allocation process under $X$, for all $X \in [0, B/\overline{N}]$, and for each of these replications to compute the number of unallocated goods that perished by the end of the horizon under this allocation, an approach which fails to scale.

To tackle this difficulty, it will be useful for us to consider a ``slow'' consumption process; at a high level, this process will allow us to bound the {\it latest possible time} any remaining unit $b$ could be allocated, with high probability. Under this process, $\underline{N}_{\leq t}$ --- a high-probability lower bound on $N_{\leq t}$ --- individuals arrive before $t+1$, $\underline{N}_{\leq t} {X}$ items are allocated up to period $t \in [T]$, and no items perish. For $b \in \mathcal{B}$, we let $\conftaub{1 \mid X, \sigma}$ be the period in which $b$ {\it would} have been entirely allocated under this slow consumption process. Formally,
\begin{align}
\conftaub{1 \mid X, \sigma} = \inf\left\{t \geq 1 : \underline{N}_{\leq t} X \geq \sigma(b)\right\}.
\end{align}
$\conftaub{1 \mid X,\ \sigma}$ represents an upper bound on the time an algorithm using static allocation $X$ would allocate $b$, since items ranked higher than $b$ may have perished, thus decreasing the time at which $b$ is allocated. We provide an example of the $\conftaub{1 \mid X, \sigma}$ construction below.

\begin{example}[continued from \Cref{ex:sigma,ex:offset,ex:sigma-is-important}]\label{ex:tau_b}
Fix $X = 3/4$.
Under deterministic demand, $\underline{N}_{\leq 1} = 1$, \mbox{$\underline{N}_{\leq 2} = 3$}, and $\underline{N}_{\leq 3}  = 4$. Then, the latest possible allocation times are given by: 
\begin{align*}
\begin{cases}
 \tau_{\texttt{a}}(1 \mid {X}, \sigma) = \inf\{t\geq 1: \frac34\underline{N}_{\leq t} \geq 2 \} = 2 \\
 \tau_{\texttt{b}}(1 \mid {X}, \sigma) = \inf\{t\geq 1: \frac34\underline{N}_{\leq t} \geq 3 \} = {3}\\
 \tau_{\texttt{c}}(1 \mid {X}, \sigma) = \inf\{t\geq 1: \frac34\underline{N}_{\leq t} \geq 1 \} = 2\\
\end{cases}
\end{align*}
{In words, under the slow consumption process, the latest possible allocation time for the entirety of unit $\texttt{c}$ is $t = 2$, since only $0.75$ units have been consumed in the first round. At $t = 2$, $2.25$ units have been consumed, meaning that $\texttt{c}$ and $\texttt{a}$ have both been consumed by the end of the second round. Finally, at $t = 3$, $3$ units have been consumed, meaning that $\texttt{b}$ must also have been consumed by the end of the third round.}
\end{example}

Note that a necessary condition for any unit $b$ to have spoiled under static allocation $X$ is for its perishing time $T_b$ to be such that $T_b < \min\{T, \tau_b(1 \mid {X},\sigma)\}$. This follows from the fact that, if $T_b$ exceeds the latest possible allocation time $\tau_b(1 \mid {X},\sigma)$, then $b$ must have been allocated before it perished, by definition. Accordingly, we define $\mu(X) = \sum_{b \in \mathcal{B}} \Pr\left(T_b < \min\left\{T, \conftaub{1 \mid X, \sigma}\right\}\right)$ to be an upper bound on the expected number of units that spoil under $X$, and use this quantity to define the worst-case loss $\overline{\Delta}(X)$. Formally:
\begin{align}
\label{eq:overline_delta}
\overline{\Delta}({X}) = \min\left\{B, \mu(X) + \conf^P_1(\mu(X))\right\},
\end{align}
where $\conf^P_1(\mu(X))$ is an appropriately chosen confidence bound, to be specified later. 

\begin{figure}[!t]
    \centering{
    \includegraphics[scale=0.4]{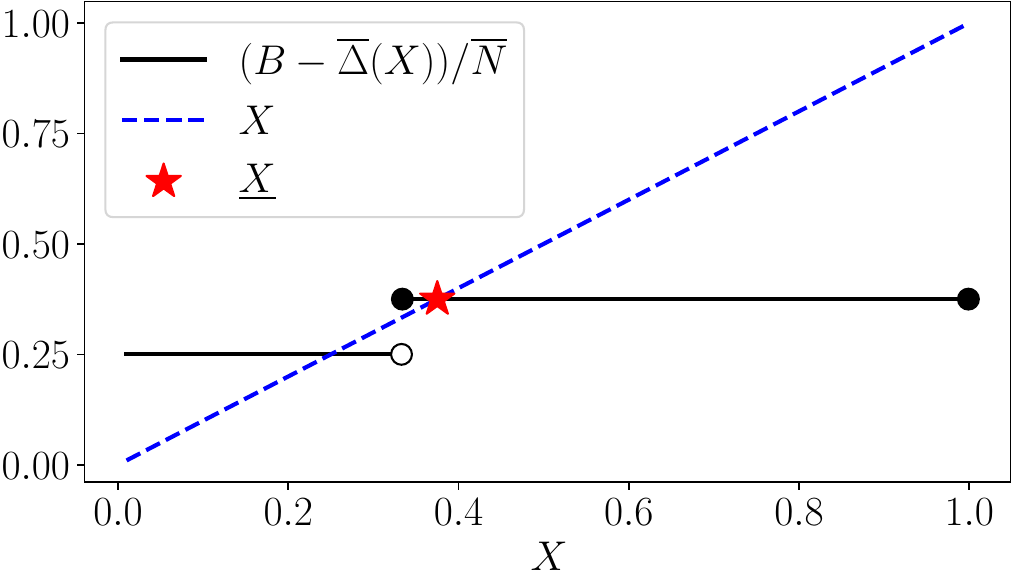}
    }
    \caption{ Illustrating the $\underline{X}$ construction for \cref{ex:xlower}. The dashed line corresponds to the line $Y = X$, and the solid line to $(B - \overline{\Delta}(X)) / \overline{N}$.  Here, $\underline{X}$ is represented by the star, the point at which the solid and dashed lines intersect.
 \label{fig:fixed-point} }
\end{figure}

With this machinery in hand, we illustrate the $\underline{X}$ construction for our running example.
\begin{example}[Continued from \Cref{ex:sigma,ex:offset,ex:sigma-is-important,ex:tau_b}]\label{ex:xlower}
Suppose the perishing distributions are as in \Cref{ex:sigma-is-important}, $\overline{N} = 4$, and $\conf_1^P(\mu(X)) = 0$ for all $X$. Let $X = B/\overline{N} = 3/4$. Continuing on from \Cref{ex:tau_b}, we have:
% \begin{align*}
% \mu(X) = \mathbb{P}\left(T_{\texttt{a}} < 2\right) + \mathbb{P}\left(T_{\texttt{b}} < 3\right) + \mathbb{P}\left(T_{\texttt{c}} < 2\right) = 3/2 \implies \overline{\Delta}(X) = 3/2.
% \end{align*}
\begin{align*}
\mu(X)
&= \mathbb{P}\left(T_{\texttt{a}} < 2\right)
 + \mathbb{P}\left(T_{\texttt{b}} < 3\right)
 + \mathbb{P}\left(T_{\texttt{c}} < 2\right) = \frac{3}{2} \implies \overline{\Delta}(X) = \frac{3}{2}.
\end{align*}
Then, $\frac{B-\overline{\Delta}(X)}{\overline{N}} = \frac{3-3/2}{4} = 3/8 < X$, meaning that setting an allocation of $X = 3/4$ is infeasible in this example. In this case, the solution to \Cref{eq:xlower} is $\underline{X} = 0.375$, as illustrated in 
\Cref{fig:fixed-point}.
\end{example}

The above example highlights that perishability forces the decision-maker to allocate less than the asymptotically optimal perishing-free allocation of $B/\overline{N}$, in order to guarantee envy-freeness. The notion of unavoidable {\it perishing-induced loss}, which we next define, uses the $\underline{X}$ construction to quantify how much less the decision-maker must allocate to account for spoilage that is out of her control.
\begin{definition}[Perishing-induced loss]\label{def:sigma-loss}
We define $\dperish = \frac{B}{\overline{N}} - \underline{X}$ to be any algorithm's \emph{perishing-induced loss}. We moreover term $T\dperish$ to be the {\it cumulative} perishing-induced loss.
\end{definition}
\cref{thm:lower_bound} establishes that any online algorithm's performance necessarily scales with $\dperish$.

\begin{theorem}
\label{thm:lower_bound}
    There exists an instance such that, for any online algorithm that is feasible with probability at least $\alpha$, the following holds with probability at least $\alpha$:
    \[
    \Deff \geq T \dperish \quad\quad \Denv \geq \dperish.
    \]
%    with probability at least $\alpha$.
\end{theorem}

\begin{rproof}
Consider an instance with $B=T$, $|\Theta| = 1$ and $N_t = 1$ for all $t$. Suppose moreover that resources have deterministic perishing times, with $T_b = b$ for all $b$, and $\sigma(b) = T+1-b$. Since the perishing and demand processes are deterministic, we let $\conf_1^P(\mu(X)) = 0$ for all $X$, and $\overline{N} = N$. For ease of notation, we omit the dependence of all quantities on $\theta$ in the remainder of the proof. The following lemma states that under this flipped ordering, any online algorithm is severely limited in its total allocation.

\begin{lemma}
\label{lem:lower_bound_total_alloc}
Any feasible algorithm must have $\sum_t X_t^{alg} \leq 1$.
\end{lemma}
\begin{rproof}
For any feasible algorithm, there must exist an available unit in period $T$. Since the only unit that has not perished by $t = T$ is $b = T$, it must be that this unit is available in period $T$. Thus, we must have $\sum_t X_t^{alg} \leq 1$ (else, $b = T$ will have been allocated before $T$).
\end{rproof}

\cref{lem:lower_bound_total_alloc} implies that any feasible {static} algorithm, which allocates a fixed amount $X_t^{alg} = X$ for all $t$, must have $X \leq \frac{1}{N}$. We use this fact to bound $\overline{\Delta}(X)$, for any feasible {static} allocation.

\begin{lemma}\label{lem:delta-bound}
For any $0 < X \leq \frac{1}{N}$, $\overline{\Delta}(X) \geq T - 1$.
\end{lemma}
\begin{rproof}
By definition: 
% \begin{align*}
% \conftaub{1 \mid X, \sigma} = \inf \{ t > 0 : N_{\leq t} X \geq \sigma(b)\} = \inf \{ t > 0 : t X \geq \sigma(b)\} = \lceil \frac{\sigma(b)}{X} \rceil =\lceil \frac{T+1-b}{X} \rceil \\ \geq T(T+1-b),
% \end{align*}
\begin{align*}
 \conftaub{1 \mid X, \sigma} = \inf \{ t > 0 : N_{\leq t} X \geq \sigma(b)\} 
= \inf \{ t > 0 : t X \geq \sigma(b)\}  = \lceil \frac{\sigma(b)}{X} \rceil  =\lceil \frac{T+1-b}{X} \rceil \\ \geq T(T+1-b),
\end{align*}

where the final inequality follows from the fact that $X \leq \frac1N = \frac1T$.
Hence, 
\begin{align*}
\overline{\Delta}(X) & = \sum_b \mathds{1}\left\{T_b < \min\{T,\conftaub{1 \mid X,\sigma}\}\right\} \geq \sum_b \mathds{1}\left\{b < \min\{T,T({T+1-b})\}\right\} = \sum_b \mathds{1}\left\{b < T\right\} = T-1,
\end{align*}
where the inequality uses the lower bound on $\conftaub{1 \mid X,\sigma}$ in addition to the assumption that~$T_b = b$.
\end{rproof}

Putting \cref{lem:lower_bound_total_alloc,lem:delta-bound} together, we have: 
\begin{align*}
\underline{X} & := \sup \{ X \mid X \leq \frac{B - \overline{\Delta}(X)}{N} \} \leq \sup \{ X \mid X \leq \frac{T - (T - 1)}{T} \}   = \frac{1}{T} \implies \dperish = 1-1/T,
\end{align*}
where the final equality follows from the fact that $\overline{\Delta}(X) = T-1$ for $X = 1/T$, and therefore $X \leq \frac{B - \overline{\Delta}(X)}{N}$ is tight at $X = 1/T$.

We now show the lower bounds on $\Denv$ and $\Deff$. By \cref{lem:lower_bound_total_alloc}, $\sum_t X_t^{alg} \leq 1$, which implies that $\min_t X_t^{alg} \leq \frac{1}{T}$.  Hence, $\Denv = \max_t |1 - X_t^{alg}| \geq 1 - \frac{1}{T} = \dperish$. Moreover, $\Deff = T - \sum_t X_t^{alg} \geq T - 1 = T \dperish$.
\end{rproof}

We conclude the section by noting the implicit dependence of $\dperish$ on the ordering $\sigma$. While $\sigma$ is an input to our model, \Cref{thm:lower_bound} implies that, if the decision-maker were to optimize over the allocation ordering, the correct choice would be to set $
    \sigma^* \in \arg\min_{\sigma} \dperish(\sigma),
$ where we emphasize the dependence of $\dperish$ on the order $\sigma$. Computing such a $\sigma^*$, however, is infeasible given the space of $N!$ orderings, as discussed in \Cref{ssec:discuss-model}. In \Cref{sec:special-cases} we identify sufficient conditions on the perishing process and allocation order that guarantee that $\dperish$ vanish asymptotically. We then use these conditions to identify practical and easily interpretable allocation schedules in our numerical experiments (see \cref{sec:experiments_other}).

\begin{remark}
We henceforth assume for simplicity that the supremum on the right-hand side of \Cref{eq:xlower} is attained. This is without loss to our results, since our bounds depend on the perishing-induced loss $\dperish$. Hence, if the supremum fails to be attained, one can define $\underline{X} = X^*-\epsilon$, where $X^*$ is the point of discontinuity and $\epsilon = o(1)$ guarantees feasibility of $\underline{X}$.
\end{remark}

\begin{remark}
Since $N_{\leq t} \geq t$ for all $t\in [T]$ one can similarly define $\conftaub{1 \mid X,\sigma} = \inf\{t \geq 1: t X \geq \sigma(b)\} = \lceil \frac{\sigma(b)}{X} \rceil$ as an upper bound on the latest possible allocation time. This quantity can then be interpreted as the ``effective rank'' of item $b$. This interpretable simplification comes at the cost of our algorithm's practical performance, but does not affect our subsequent theoretical bounds.
\end{remark}

%% file: final_arxiv_parts/results.tex
\section{Achieving the Optimal Trade-Off Via Perishing-Aware Guardrails}\label{sec:alg}

In this section we design an algorithm that achieves the lower bounds described above. We moreover instantiate our bounds for important special cases.

\subsection{Algorithm Description}

Our algorithm, \HopeGuardrailPerish, builds on the guardrail algorithm proposed by \citet{sinclair2021sequential} for the no-perishing setting. In particular, our algorithm always allocates the same amount to all individuals in a given round (independent of their type). As a result, the only source of envy is due to differences in allocations across rounds, as opposed to across individuals within the same round. Our algorithm controls this source of envy via an input parameter $L_T = o(1)$; it then tackles the problem of how much to allocate in each round in three main steps:
\begin{enumerate}
\item Constructing a {static} allocation (also referred to as {\it baseline allocation} or {\it lower guardrail}), $\underline{X} \in \mathbb{R}_+$.
\item Setting an {``aggressive''} allocation $\overline{X} = \underline{X} + L_T$ to improve efficiency.
\item Determining an appropriate {threshold condition} that indicates when to allocate aggressively (i.e., in which periods to allocate $\overline{X}$ instead of $\underline{X}$).
\end{enumerate}

At a high level, the lower guardrail $\underline{X}$ is constructed to be the {largest} feasible allocation for which the hindsight envy is zero (i.e., all individuals receive the same amount in each period), with high probability. For this to hold, it must be that the algorithm never runs out of budget under $\underline{X}$. This desideratum is precisely satisfied by the static allocation defined in \Cref{eq:xlower}, i.e., $$\underline{X} = \sup\left\{X \mid X \leq \frac{B-\overline{\Delta}(X)}{\overline{N}}\right\},$$
where we let $\overline{N} = \mathbb{E}[N] + \conf_{0,T}^N$ for an appropriately defined confidence term $\conf_{0,T}^N$.

While $\underline{X}$ is budget-respecting with high probability, allocating a single budget-respecting quantity in each period may result in significant losses to efficiency. Herein lies the motivation behind the use of the ``aggressive'' allocation $\overline{X} = \underline{X} + L_T$: in periods in which the algorithm ``feels'' as though it is holding onto more units than it will need throughout the remainder of the horizon, it allocates aggressively. More specifically, in each period, our algorithm verifies whether there is enough budget remaining to accommodate $(i)$ the aggressive allocation $\overline{X}$ in the current period, $(ii)$ the lower guardrail $\underline{X}$ in all future periods, under high demand, and $(iii)$ aggressive future perishing. If enough budget remains, it allocates $\overline{X}$; otherwise, it defaults to $\underline{X}$. Formally, in any given period $t$, given a pessimistic forecast of future perishing $\overline{P}_t$, our threshold condition is as follows:
\begin{align*}
B_t^{alg} \geq N_t\overline{X} + \overline{N}_{>t}\underline{X} + \overline{P}_t,
\end{align*}
where we let $\overline{N}_{>t} = \mathbb{E}[N_{>t}] + \conf_{t, T}^N$ for an appropriately defined confidence term $\conf_{t, T}^N$.

Our main challenge lies in the construction of $\overline{P}_t$. In particular, computing a tight upper bound on future perishing poses significant challenges that do not exist in the setting without perishable resources. In this latter setting, it is easy to obtain tight upper bounds on future demand since demand uncertainty is {\it exogenous}, i.e., it is invariant to the decisions made by the algorithm. On the other hand, uncertainty around future perishing is {\it endogenous}: {as discussed in \Cref{sec:construction_x_lower}}, though the distribution around perishing times is fixed, how many --- and which --- items perish depends heavily on the rate at which items are being allocated, which itself depends on the rate at which items perish. We sidestep this challenge by turning to the ``slow'' consumption process introduced in \Cref{sec:construction_x_lower}. Namely, fix a round $t \in [T]$, and let $\tau_b(t \mid \underline{X},\sigma)$ denote the latest possible allocation time of unit $b$, starting from time $t$. (This latest possible allocation time is a time-dependent generalization of the latest possible allocation time defined in \Cref{ssec:unavoidable-loss}; we formally define it below.) In line with the intuition provided earlier, a necessary condition for $b$ to have spoiled {after $t$} is for its perishing time $T_b$ to be such that {\mbox{$t \leq T_b < \min\{T, \tau_b(t \mid \underline{X},\sigma)\}$}}. We use this necessary condition to obtain a high-probability pessimistic estimate of the number of units that perish after period $t$, which we denote by $\overline{P}_t$. Crucially, $\overline{P}_t$ is independent of our algorithm's decisions, which will allow for its tractable analysis. We formalize the construction of this key quantity below. 
 
For all $t' \geq t$, let $\underline{N}_{[t,t']} = \Exp{N_{[t,t']}} - \conf_{t,t'}^N$ be a high-probability {lower} bound on the number of arrivals between $t$ and $t'$, for an appropriately defined confidence term $\conf_{t,t'}^N$. In line with the discussion in \Cref{ssec:unavoidable-loss}, with high probability, $\underline{N}_{< t} \underline{X} + \underline{N}_{[t,t']}\underline{X}$ corresponds to the smallest quantity of units that could have been consumed by $t'$, since any units perishing before $t'$ would only have increased consumption. Consider now any unit $b$ in the set of remaining items at the beginning of period $t$. Under the slow process defined above, if $\underline{N}_{<t}\underline{X}+ \underline{N}_{[t,t']}\underline{X} \geq \sigma(b)$, then $b$ must have been allocated by $t' \geq t$, according to the allocation schedule $\sigma$. We use this intuition to define $\tau_b(t \mid \underline{X}, \sigma)$, i.e.,
\begin{align}\label{eq:conftau}
\conftaub{t \mid \underline{X}, \sigma} = \inf\left\{t' \geq t: \underline{N}_{< t} \underline{X} + \underline{N}_{[t,t']}\underline{X} \geq \sigma(b)\right\}.
\end{align}

Let $\Bupper{t}$ denote the set of items remaining in period $t$ under this slow process. Via similar logic, \mbox{$\Bupper{t} = \{\sigma^{-1}(\lceil \underline{N}_{<t}\underline{X}\rceil),\ldots,\sigma^{-1}(B)\}$}. Then, an upper bound on the expected number of items that perish from $t$ onwards is given by:
\begin{align}\label{eq:eta}
\eta_t =  \sum_{b \in \Bupper{t}} \Pr\left(t\leq T_b< \min\{T, \conftaub{t \mid \underline{X}, \sigma}\} \right).
\end{align}
Here, we deliberately use $\Bupper{t}$ --- as opposed to $\Balg{t}$, which is endogenous to the allocation process --- so that our pessimistic perishing estimate is completely exogenous to our algorithm's decisions. This will be critical in our analysis, allowing us to sidestep a union bound over exponentially many subsets of remaining items. 

With this machinery in hand, we define our high-probability pessimistic forecast $\overline{P}_t$ as follows:
\begin{align}
\begin{cases}
\Pupper_t &= \min\left\{\Pupper_{t-1}, \eta_t + \conf_t^P(\eta_t)\right\} \quad \forall \ t \in [T],\\
\Pupper_0 &= B,
\end{cases}
\end{align}
for an appropriately defined confidence term $\conf_t^P(\cdot)$. {Note that $\Pupper_1 = \overline{\Delta}(\underline{X})$ by construction, a fact we will use in the analysis of our algorithm.}  

We present our algorithm, \HopeGuardrailPerish, in Algorithm \ref{alg:brief-perishing}. Recall, for $t\in[T]$, $\Balg{t}$ is used to denote the set of remaining resources at the beginning of time $t$, and $B_t^{alg}$ the {quantity} of remaining resources at the beginning of the period. Moreover, $\PUA_t^{alg}$ denotes the quantity of unallocated items that perished at the end of round $t$.

\begin{algorithm}[!t]
\DontPrintSemicolon % Some LaTeX compilers require you to use \dontprintsemicolon instead
\KwIn{Budget $B = B_1^{alg}$, allocation schedule $\sigma$, envy parameter $L_T$,  {arrival confidence terms $(\conf_{t,t'}^N)_{t,t'\in\{0,\ldots,T\}}$, {perishing confidence terms $(\conf_t^P(\cdot))_{t\in\{1,\ldots,T\}}$,} and perishing inputs $(\eta_t)_{t\in[T]}$ given by \Cref{eq:eta}}}
\KwOut{An allocation $X^{alg} \in \mathbb{R}^{T \times |\Theta|}$}
Compute $\underline{X} = \sup \left\{ X \mid X \leq \frac{B - \overline{\Delta}(X)}{\overline{N}} \right\}$ and set $\overline{X} = \underline{X} + L_T$.\;
\For{$t = 1, \ldots, T$}{ 
    Compute {{$\overline{P}_t = \min\left\{\Pupper_{t-1},\eta_t + \conf_t^P(\eta_t) \right\}$}  \tcp*[h]{Compute ``worst-case'' future perishing}}\;
    \uIf(\tcp*[h]{insufficient budget to allocate lower guardrail}){$B_{t}^{alg} < N_t\underline{X}$}{
        Set $X_{t, \theta}^{alg} = \frac{B_{t}^{alg}}{N_t}$ for each $\theta \in \Theta$. Allocate items $b\in\Balg{t}$ according to $\sigma$.\;}
    \uElseIf(\tcp*[h]{use upper guardrail}){$B_{t}^{alg} -   N_t \overline{X} \geq \underline{X}  (\Exp{N_{> t}} + \conf_{t, T}^N) + \overline{P}_t$}{
        Set $X_{t, \theta}^{alg} = \overline{X}$ for each $\theta \in \Theta$. Allocate items $b\in\Balg{t}$ according to $\sigma$.\;}
    \uElse(\tcp*[h]{use lower guardrail}){
        Set $X_{t, \theta}^{alg} = \underline{X}$ for each $\theta \in \Theta$. Allocate items $b\in\Balg{t}$ according to $\sigma$.}
        Update $B_{t+1}^{alg} = B_{t}^{alg}-N_{t}X_{t}^{alg} - \PUA^{alg}_t$\\ 
        }
\Return{$X^{alg}$}
	\caption{\HopeGuardrailPerish}
	\label{alg:brief-perishing}
\end{algorithm}

\subsection{Performance Guarantee}\label{sec:performance-guarantee}

\Cref{thm:upper_bound_perishing} below provides upper bounds on counterfactual envy and inefficiency, given hindsight envy parameter $L_T$.

\begin{theorem}
\label{thm:upper_bound_perishing}
For $t' > t$, define the confidence terms:
\begin{itemize}
    \item $\conf^N_{t,t'} = \sqrt{2(t'-t)|\Theta|\rho_{\max}^2 \log(2T^2 / \delta)}$
    \item {$\conf^P_t(\eta_t) = \frac{1}{2}\left(\log(3 t \log T / \delta) + \sqrt{\log^2(3 t \log T / \delta) + 8 \eta_t \log(3 t \log T / \delta)}\right)$}
\end{itemize}
Then, with probability at least $1-2\delta$, \Cref{alg:brief-perishing} achieves:
\begin{align*}
     \Denv  \lesssim \max \{L_T, \dperish + 1 / \sqrt{T} \} \qquad &
     \Deff  \lesssim \min\left\{\sqrt{T},L_T^{-1} + \sqrt{TL_T^{-1} \dperish} \right\}  + T \dperish  \\
   \Denvhind & \lesssim L_T 
\end{align*}
% \begin{align*}
%      \Denv \lesssim \max \{L_T, \dperish + 1 / \sqrt{T} \}&  & \Deff \lesssim \min\left\{\sqrt{T},L_T^{-1} + \sqrt{TL_T^{-1} \dperish} \right\} + T \dperish  \\
%    & \Denvhind \lesssim L_T &  %& \Dprop \lesssim \max \{L_T, \dperish + 1 / \sqrt{T}\} 
% \end{align*}
where $\lesssim$ drops poly-logarithmic factors of $T$, $\log(1/\delta)$, and absolute constants. 
\end{theorem}

In Appendix \ref{apx:tightness}, we show that these bounds are tight relative to the lower bounds in Theorems \ref{thm:lower_bounds} and \ref{thm:lower_bound}.
We dedicate the remainder of this section to further parsing the bounds on counterfactual envy and efficiency given the scalings of $L_T$ and $\dperish$. 
We prove \cref{thm:upper_bound_perishing} in \cref{sec:analysis}, and instantiate $\dperish$ for important special cases in \Cref{sec:special-cases}.

\begin{figure}[t]
\centering
\scalebox{.85}{
\tikzset{every picture/.style={line width=0.75pt}} %set default line width to 0.75pt 
\input{figures/intro_tikz/phase_transition}}
\caption{\centering
Graphical representation of the regimes described in \cref{cor:low-envy,cor:main_theorem_interpret,cor:high-envy-high-perish}. The curve $T\dperish L_T = 1$, which partitions the two cases of \Cref{cor:main_theorem_interpret}, is shown in gray.
\label{fig:phase_transition} }
\end{figure}

 {The following corollary simplifies our bounds in the ``low-envy'' setting where $L_T \lesssim 1/\sqrt{T}$.

 \begin{corollary}[Low-Envy]\label{cor:low-envy}
Suppose $L_T \lesssim 1/\sqrt{T}$. Then, \HopeGuardrailPerish achieves with probability at least $1-2\delta$:
\begin{align*}
&\Denv \lesssim \dperish + 1/\sqrt{T} & \Deff\lesssim \sqrt{T} + T\dperish.
\end{align*}
\end{corollary}
\cref{cor:low-envy} implies that there is no efficiency benefit to increasing $L_T$ as long as $L_T \lesssim 1/\sqrt{T}$. When $\dperish \lesssim 1/\sqrt{T}$, our algorithm incurs $\widetilde{O}(\sqrt{T})$ envy and inefficiency. In this case, these quantities are driven by the exogenous uncertainty in demand. When $\dperish \gtrsim 1/\sqrt{T}$, on the other hand, envy and inefficiency are driven by the unavoidable perishing-induced loss.

\cref{cor:high-envy-high-perish} next considers the ``high-envy, high-perishing'' setting, on the other extreme of the spectrum. 
\begin{corollary}[High-Envy, High-Perishing]\label{cor:high-envy-high-perish}
Suppose $L_T \gtrsim 1/\sqrt{T}$ and $\dperish \gtrsim 1/\sqrt{T}$. Then, \HopeGuardrailPerish achieves with probability at least $1-2\delta$:
\begin{align*}
&\Denv \lesssim \max\{L_T, \dperish\} &  \Deff \lesssim T\dperish.
\end{align*}
\end{corollary}
 Similarly in this regime, increasing $L_T$ doesn't guarantee arbitrary gains in efficiency; thus, setting \mbox{$L_T \sim \dperish$} is optimal. 
 We conclude by exploring the more nuanced ``high-envy, low-perishing'' regime. In this setting, our algorithm's guarantees depend on whether the efficiency gain from increasing envy, $L_T^{-1}$, exceeds $T\dperish$. We defer its proof to Appendix \ref{apx:main_theorem_interpret}.

\begin{corollary}[High-Envy, Low-Perishing]
\label{cor:main_theorem_interpret}
Suppose $L_T \gtrsim 1 / \sqrt{T}$, $\dperish \lesssim 1 / \sqrt{T}$, and \mbox{$T\dperish \lesssim L_T^{-1}$}. Then, \HopeGuardrailPerish achieves with probability at least $1 - 2\delta$:
\begin{align*}
    & \Denv \lesssim  L_T  & \Deff \lesssim L_T^{-1} .%\\
    % & \Denvhind \lesssim L_T & \Dprop \lesssim L_T 
\end{align*}
On the other hand, if {$T\dperish \gtrsim L_T^{-1}$}, \HopeGuardrailPerish achieves with probability at least $1 - 2 \delta$:
\begin{align*}
    & \Denv \lesssim  L_T  & \Deff \lesssim T \dperish. %\\
    % & \Denvhind \lesssim L_T & \Dprop \lesssim L_T 
\end{align*}
\end{corollary}
Thus, in ``high-envy, low-perishing'' regimes, if the cumulative {perishing-induced loss} is order-wise dominated by the efficiency gains from envy (otherwise phrased, our allocation schedule $\sigma$ is high-enough quality {and the perishing rate is mild}), increasing $L_T$ allows \HopeGuardrailPerish to achieve inversely proportional gains in efficiency, as in the setting without perishable resources. One can do this until moving into the regime where $L_T^{-1} \lesssim T\dperish$ (i.e., the cumulative perishing-induced loss dominates efficiency gains from envy). At this point, further increasing $L_T$ hurts envy, and has no order-wise impact on efficiency. We summarize these four regimes pictorially in \Cref{fig:phase_transition}}.

We conclude the section by noting that our result is a strict generalization of \citet{sinclair2021sequential}: in the simplest setting where \mbox{$\dperish = 0$} (i.e., no perishing, or the process is offset-expiring with perfect predictions), we recover the trade-off described in \Cref{thm:lower_bounds}.

\subsection{Analysis}\label{sec:analysis}

The proof of \cref{thm:upper_bound_perishing} is based on three main building blocks:
\begin{itemize}
    \item {\bf Defining and bounding the ``good event''} (\cref{sec:good_event}): We first show that, with probability at least $1-2\delta$, the realizations of future arrivals and perishing are no worse than our algorithm's pessimistic predictions. As a result, it suffices to condition over such ``good'' sample paths. 
    \item {\bf Establishing feasibility of $\underline{X}$} (\cref{sec:feasibility_allocations}): For this good event, we show that the static allocation $\underline{X}$ computed at the start of the horizon will never exhaust the budget, despite incurring perishing-induced loss.
    \item {\bf Improving efficiency via $\overline{X}$} (\cref{sec:t0-sec}): We next show that the threshold condition guarantees that the algorithm allocates aggressively enough throughout the horizon to ensure high efficiency. 
\end{itemize}
We use these building blocks for our final bounds on envy and efficiency in \cref{sec:proof_final_result}.

\subsubsection{Defining and bounding the ``good event.''}\label{sec:good_event}

We analyze the performance of our algorithm under a so-called ``good event'' $\mathcal{E}$, the intersection of the following two events:
\begin{enumerate}
\item  $\E_N = \left\{|N_{(t, t']} - \Exp{N_{(t,t']}}| \leq \conf^N_{t,t'}\ \forall \ t, t' > t  \right\},$
\item $\E_{\overline{P}} = \left\{\overline{P}_t \geq \PUA^{alg}_{\geq t} \, \forall \ t \in \{0,1,\ldots,T\} \right\}$, where $\PUA^{alg}_{\geq t}$ denotes the quantity of unallocated items that perished between the end of round $t$ and the end of round $T-1$.
\end{enumerate}
$\E$ represents the event that the arrival process falls close to its mean, and that $\overline{P}_t$ is a pessimistic estimate of the unallocated goods that perish in the future.

The following lemma provides the high-probability bound on $\Pr(\E_N)$. We defer its proof --- which follows from a standard application of Hoeffding's inequality --- to Appendix \ref{apx:hoeffding_app-perish}.

\begin{lemma}%[Concentration on Arrival Process]
\label{lem:hoeffding_app-perish}
$\E_N$ holds with probability at least $1-\delta$.
\end{lemma}

The main challenge in analyzing $\E$ lies in showing that $\Pupper_t$ is indeed a pessimistic estimate of our algorithm's future perishing. To upper bound the amount of unallocated resources that perished between $t$ and $T-1$, we must account for both the uncertainty in arrivals and the realized order in which resources perished, and relate these two sources of uncertainty to the time at which the algorithm intended to allocate these resources. Establishing this upper bound hinges upon the careful construction of the ``slow'' consumption process. In particular, leveraging $\tau_b(t \mid \underline{X}, \sigma)$ as an upper bound on the time that each unit is allocated by our algorithm, our bound on the number of units that perish from $t$ onwards is a sum of independent random variables that is decoupled from the algorithm's past (and future) decisions, thereby allowing us to tractably obtain tight bounds on the true perishing process. We formalize these ideas in \cref{lem:p_upper_bound}.

\begin{lemma}
\label{lem:p_upper_bound}
Given $\E_N$, $\E_{\overline{P}}$ holds with probability at least $1-\delta$.  %Moreover, $\Pupper_t$ is non-decreasing in $t$.%, with $\Pupper_1 = \overline{\Delta}(\underline{X})$.
\end{lemma}
\begin{rproofof}{\cref{lem:p_upper_bound}}
We prove the claim by induction. 

\noindent\textbf{Base case:} $t=0$. Since $\Pupper_0 = B$ by definition, we have $\Pupper_0 \geq \PUA^{alg}_{\geq 0}$ trivially.

\noindent\textbf{Inductive step:} $t-1\rightarrow t.$ %{We first upper bound $\PUA^{alg}_{\geq t}$ as a function of the worst-case perishing times $\conftaub{t \mid \underline{X}, \sigma}$.}
Let $\mathcal{A}_{\tau}$ be the set of items allocated in round $\tau$, for $\tau \geq t$. Since $\PUA^{alg}_{\geq t}$ represents the amount of unallocated goods that perished between $t$ and $T-1$, we have:
\begin{align}
\PUA^{alg}_{\geq t}  & {=} \sum_{\tau = t}^{T-1} \sum_{b \in \Balg{\tau}} \Ind{T_b = \tau, b \not\in \A_\tau} \notag \\ 
& = \sum_{b\in\Balg{t}}\Ind{T_b \geq t, b \text{ not allocated before }T_b} \notag \\
& \leq \sum_{b \in \Balg{t}} \Ind{t \leq T_b < \min\{T, \conftaub{t \mid \underline{X}, \sigma}\}} \label{eq:pua_ub_btalg}, 
\end{align}

% \begin{align}
% \PUA^{alg}_{\geq t}  {=} \sum_{\tau = t}^{T-1} \sum_{b \in \Balg{\tau}} \Ind{T_b = \tau, b \not\in \A_\tau} \\& = \sum_{b\in\Balg{t}}\Ind{T_b \geq t, b \text{ not allocated before }T_b} \notag \\
% & \leq \sum_{b \in \Balg{t}} \Ind{t \leq T_b < \min\{T, \conftaub{t \mid \underline{X}, \sigma}\}} \label{eq:pua_ub_btalg}, 
% \end{align}
where the inequality follows from the fact that {${t \leq T_b < \min\{T,\tau_b(t\mid \underline{X},\sigma)\}}$} is a necessary condition for item $b$ to have spoiled.

A similar argument gives us that $\Balg{t} \subseteq \Bupper{t}$. Plugging this into \eqref{eq:pua_ub_btalg}:
\begin{align}\label{eq:pua-ub-btalg2}
\PUA^{alg}_{\geq t} &\leq  \sum_{b \in \Bupper{t}} \Ind{t\leq T_b < \min\{T, \conftaub{t \mid \underline{X}, \sigma} \}}.
\end{align}

{Recall, $\eta_t = \sum_{b\in\Bupper{t}} \Pr(t \leq T_b < \min\{T, \conftaub{t \mid \underline{X}, \sigma}\})$. Applying a Chernoff bound (see \cref{cor:chernoff_bernoulli}) to the right-hand side of \eqref{eq:pua-ub-btalg2}, we obtain that, with probability at least $1-\delta/(3 t \log(T))$: 
\begin{align*}
 \PUA^{alg}_{\geq t}
&\leq \eta_t
   + \frac{1}{2}
     \Bigl( \log(3t\log(T)/\delta)  + \sqrt{
       \log^2(3t\log(T)/\delta)
       + 8\eta_t\log(3t\log(T)/\delta)
     } \Bigr) \\
&= \eta_t + \conf_t^P(\eta_t).
\end{align*}
% \begin{align*}
% \PUA^{alg}_{\geq t} 
% & \leq \eta_t + \frac{1}{2}\left( \log(3 t \log(T) / \delta) + \sqrt{\log^2(3 t \log(T) / \delta) + 8 \eta_t \log(3 t \log(T) / \delta)} \right) \\
% & = \eta_t + \conf_t^P(\eta_t).
% \end{align*}
Moreover, $\PUA^{alg}_{\geq t} \leq \PUA^{alg}_{\geq t-1} \leq \Pupper_{t-1}$, where the second inequality follows from the inductive hypothesis. Putting these two facts together, we obtain:
\begin{align*}
\PUA^{alg}_{\geq t} \leq \min\left\{ \Pupper_{t-1}, \eta_t + \conf_t^P(\eta_t) \right\} = \Pupper_t
\end{align*}
with probability at least $1-\delta/(3t \log(T))$.  A union bound over $t$ completes the proof of the result.}
\end{rproofof}

\cref{cor:4-delta-lem} follows from these high-probability bounds. We defer its proof to Appendix \ref{apx:4-delta-lem}.
\begin{lemma}\label{cor:4-delta-lem}
Let $\E = \E_N \cap \E_{\overline{P}}$. Then, $\Pr(\E) \geq 1 - 2\delta$.
\end{lemma}
In the remainder of the proof, it suffices to restrict our attention to $\E$.

\subsubsection{Feasibility of $\underline{X}$.}
\label{sec:feasibility_allocations}

We now show that, given $\E$, our algorithm never runs out of budget, and as a result always allocates $X_{t,\theta}^{alg} \in \{\underline{X},\overline{X}\}$. Since $X_{t,\theta}^{alg} = X_{t,\theta'}^{alg}$ for all $\theta,\theta'$, for ease of notation in the remainder of the proof we omit the dependence of $X_{t,\theta}^{alg}$ on $\theta$.
%We first start off by showing that given $\mathcal{E}$, the counterfactual algorithm which allocates the static allocation $\underline{X}$ in each round is feasible.  This will be used to highlight feasibility of the algorithm, since the hope is that the worst case is the algorithm allocates $\underline{X}$ at each time period.

\begin{lemma}\label{lem:enough-budget-perish-app}
Under event $\mathcal{E}$, $B_t^{alg} \geq N_{\geq t}\underline{X}$ for all $t \in [T]$.
\end{lemma}
\begin{rproofof}{\cref{lem:enough-budget-perish-app}}
By induction on $t$.

\noindent \textbf{Base Case}: $t = 1$. By definition:
\begin{align*}
&\underline{X} \leq \frac{B-\overline{\Delta}(\underline{X})}{\overline{N}} \implies B \geq \overline{N}\underline{X} + \overline{\Delta}(\underline{X}) \geq N\underline{X},
\end{align*}
where the final inequality follows from $\overline{N} \geq N$ under $\mathcal{E}$, and $\overline{\Delta}(\underline{X}) \geq 0$.

\noindent \textbf{Step Case}: $t - 1 \rightarrow t$. We condition our analysis on $(X_{\tau}^{alg})_{\tau < t}$, the algorithm's previous allocations.

\noindent \textit{Case 1}: $X_\tau^{alg} = \underline{X}$ for all $\tau < t$. By the recursive budget update, $B_t^{alg} = B - N_{< t} \underline{X} - \PUA^{alg}_{< t}$, %\\
where $\PUA^{alg}_{< t}$ denotes the quantity of unallocated goods that perished before the end of round $t$. To show that $B_t^{alg} \geq N_{\geq t} \underline{X}$, it then suffices to show that $B - \PUA^{alg}_{< t} \geq N \underline{X}$. We have:
\begin{align*}
\PUA^{alg}_{< t} \leq \PUA^{alg}_{\geq 1} \leq \Pupper_1 = \overline{\Delta}(\underline{X}),
\end{align*}
where the final inequality follows from \cref{lem:p_upper_bound}. Under $\E$, then, as in the base case:
\begin{align*}
B - \PUA^{alg}_{< t} \geq B-\overline{\Delta}(\underline{X}) \geq N\underline{X}.
\end{align*}

\noindent \textit{Case 2}: There exists $\tau < t$ such that $X_\tau^{alg} = \overline{X}$. Let $\tstar = \sup\{\tau < t : X_\tau^{alg} = \overline{X}\}$ {be the most recent time the algorithm allocated $\overline{X}$}. Again, by the recursive budget update:
    $$B_{t}^{alg} = B_{\tstar}^{alg} - N_{\tstar} \overline{X} - N_{(\tstar, t)} \underline{X} - \PUA^{alg}_{[\tstar, t)}.$$
Since $\overline{X}$ was allocated at $\tstar$, it must have been that $B_{\tstar}^{alg} \geq N_{\tstar}\overline{X} + \overline{N}_{>\tstar}\underline{X} + \overline{P}_{\tstar}$. Plugging this into the above and simplifying:
% \begin{align*}
%     B_t^{alg} & \geq N_{\tstar}\overline{X} + \overline{N}_{>\tstar}\underline{X} + \overline{P}_{\tstar} - N_{\tstar} \overline{X} - N_{(\tstar, t)} \underline{X} - \PUA^{alg}_{[\tstar, t)}\\
%     & = \overline{N}_{> \tstar} \underline{X}  - N_{(\tstar, t)} \underline{X} + \overline{P}_{\tstar} - \PUA^{alg}_{[\tstar, t)} \\
%     &\geq N_{\geq t} \underline{X} + \overline{P}_{\tstar} - \PUA^{alg}_{[\tstar, t)},
% \end{align*}
\begin{align*}
    B_t^{alg} & \geq N_{\tstar}\overline{X} + \overline{N}_{>\tstar}\underline{X} + \overline{P}_{\tstar} - N_{\tstar} \overline{X} - N_{(\tstar, t)} \underline{X} - \PUA^{alg}_{[\tstar, t)}\\
    & = \overline{N}_{> \tstar} \underline{X}  - N_{(\tstar, t)} \underline{X} + \overline{P}_{\tstar} - \PUA^{alg}_{[\tstar, t)} \\
    &\geq N_{\geq t} \underline{X} + \overline{P}_{\tstar} - \PUA^{alg}_{[\tstar, t)},
\end{align*}
where the second inequality follows from the fact that $\overline{N}_{> \tstar} \geq N_{> \tstar}$ under $\mathcal{E}$. Thus, it suffices to show that $\overline{P}_{t^*} \geq \PUA^{alg}_{[t^*,t)}$. This holds since $\PUA^{alg}_{[t^*, t)} \leq \PUA^{alg}_{\geq t^*} \leq \overline{P}_{t^*}$ by \cref{lem:p_upper_bound}.
\end{rproofof}

\subsubsection{Improving efficiency via $\overline{X}$.}\label{sec:t0-sec}

{Having established that the algorithm never runs out of budget, it remains to investigate the gains from allocating $\overline{X}$. By the threshold condition, whenever the algorithm allocates $\overline{X}$ it must be that there is enough budget remaining to allocate $\overline{X}$ in the current period, and $\underline{X}$ in all future periods, under high demand {\it and} high perishing. Thus, at a high level, $\overline{X}$ being allocated is an indication that the algorithm has been inefficient up until round $t$. The following lemma provides a lower bound on the last time the algorithm allocates $\underline{X}$. This lower bound will later on allow us to establish that, for most of the time horizon, the remaining budget is low relative to future demand, ensuring high efficiency.
%
% $B_{t}^{alg} < N_t \overline{X} + \underline{X}  (\Exp{N_{> t}} + \conf_{t, T}^N) + \overline{P}_t$. The following lemma provides a lower bound on the last time the algorithm allocates $\underline{X}$. This lower bound will later on allow us to establish that, for the majority of the time horizon, the remaining budget is low relative to future arrivals, thus ensuring high efficiency. %However, by showing that there is a time close to the end of the horizon such that we allocate $\underline{X}$, we are able to use the fact that the confidence radius becomes tighter, obtaining a better bound on efficiency.  We first address this with the following lemma, which states that there exists a final timestep $t_0$ close to the end of the time horizon for which the algorithm allocates $\underline{X}$ (and allocates $\overline{X}$ from then on).

\begin{lemma}
\label{lem:switching_point_last_perish}
Given $\mathcal{E}$, let $t_{0} = \sup \{t: X_{t}^{alg} = \underline{X} \}$ be the last time that $X_{t}^{alg} = \underline{X}$ (or else $0$ if the algorithm always allocates according to $\overline{X}$).  Then, for some $\tilde{c} = \widetilde{\Theta}(1)$, 
\[
t_{0} > T - \tilde{c} \left(\frac{1}{L_T} + \sqrt{\frac{T \dperish}{L_T}}\right)^2.
\]
\end{lemma}
We defer the proof of \cref{lem:switching_point_last_perish} to Appendix \ref{apx:switching_point_last_perish}. Observe that, as $\dperish$ increases, our algorithm stops allocating $\underline{X}$ earlier on. {We will next see how this loss propagates to our final efficiency bound.}
}

\subsubsection{Putting it all together.}
\label{sec:proof_final_result}

With these building blocks in hand, we prove our main result.
\begin{rproofof}{\cref{thm:upper_bound_perishing}}
By \cref{lem:enough-budget-perish-app}, the algorithm never runs out of budget under event $\mathcal{E}$, {which occurs with probability at least $1 - 2 \delta$}. As a result $X_{t,\theta}^{alg} \in \{\underline{X},\overline{X}\}$ for all $t\in[T]$, $\theta\in\Theta$. We use this to bound envy and efficiency.

\noindent \textbf{Counterfactual Envy}:  Recall, $\Denv = \max_{t, \theta} |w_\theta(X_{t, \theta}^{alg} - \frac{B}{N})| \leq w_{\max}\cdot |X_{t,\theta}^{alg}-\frac{B}{N}|$.  We consider two cases.

\textit{Case 1: $\underline{X} \leq \overline{X} \leq \frac{B}{N}$.}  By definition: 
\begin{align*}
    \frac{B}{N} - \underline{X} = \frac{B}{N} - \frac{B}{\overline{N}} + \frac{B}{\overline{N}} - \underline{X} \leq \frac{B}{\underline{N}} - \frac{B}{\overline{N}} + \dperish,
\end{align*}
where the inequality follows from the fact that $\underline{N} \leq N$ under $\E$, and $\dperish = B/\overline{N}-\underline{X}$ by definition.
We turn our attention to the first two terms:
\begin{align*}
   \frac{B}{\underline{N}} - \frac{B}{\overline{N}} &{=}  \frac{B}{\Exp{N} - \conf_{0,T}^N} - \frac{B}{\Exp{N} + \conf_{0,T}^N} \\
   &= \frac{B}{\Exp{N}}\left(\frac{1}{1-\frac{\conf_{0,T}^N}{\Exp{N}}} - \frac{1}{1 + \frac{\conf_{0,T}^N}{\Exp{N}}}\right) \\
   & = \Bav\left(\frac{1}{1-\frac{\conf_{0,T}^N}{\Exp{N}}} - \frac{1}{1 + \frac{\conf_{0,T}^N}{\Exp{N}}}\right).
\end{align*}
% \begin{align*}
%    \frac{B}{\underline{N}} - \frac{B}{\overline{N}} &{=}  \frac{B}{\Exp{N} - \conf_{0,T}^N} - \frac{B}{\Exp{N} + \conf_{0,T}^N} \\
%    &= \frac{B}{\Exp{N}}\left(\frac{1}{1-\frac{\conf_{0,T}^N}{\Exp{N}}} - \frac{1}{1 + \frac{\conf_{0,T}^N}{\Exp{N}}}\right)
%    = \Bav\left(\frac{1}{1-\frac{\conf_{0,T}^N}{\Exp{N}}} - \frac{1}{1 + \frac{\conf_{0,T}^N}{\Exp{N}}}\right).
% \end{align*}
Using the fact that $\conf_{0,T}^N =\sqrt{2T|\Theta|\rho_{\max}^2\log(2T^2/\delta)}$ and $\Exp{N} = \Theta(T)$, there exists $c_1, c_2 =\widetilde{\Theta}(1)$ such that, for large enough $T$, $\left(1-\frac{\conf_{0,T}^N}{\Exp{N}}\right)^{-1} \leq \left(1-c_1/\sqrt{T}\right)^{-1}\leq 1+2c_1/\sqrt{T}$  and $\left(1+\frac{\conf_{0,T}^N}{\Exp{N}}\right)^{-1} \geq \left(1+c_2/\sqrt{T}\right)^{-1} \geq 1-c_2/\sqrt{T}$.  Plugging this into the above:
\begin{align*}
    \frac{B}{\underline{N}} - \frac{B}{\overline{N}} &\leq \Bav\left(1+2c_1/\sqrt{T}-(1-c_2/\sqrt{T})\right) \\
    & \leq \Bav(2c_1+c_2)/\sqrt{T} \lesssim 1/\sqrt{T}.
\end{align*}
% \begin{align*}
%     \frac{B}{\underline{N}} - \frac{B}{\overline{N}} &\leq \Bav\left(1+2c_1/\sqrt{T}-(1-c_2/\sqrt{T})\right) \leq \Bav(2c_1+c_2)/\sqrt{T} \lesssim 1/\sqrt{T}.
% \end{align*}
Thus, we obtain $|X_{t, \theta}^{alg} - B/N| \lesssim 1 / \sqrt{T} + \dperish$.

\textit{Case 2: $\underline{X} \leq \frac{B}{N} \leq \overline{X}$.} We have:
%\begin{align*}
    $|X_{t, \theta}^{alg} - B/N|  = \max\left\{\frac{B}{N} - \underline{X}, \overline{X}-\frac{B}{N} \right\} \leq \overline{X} - \underline{X} = L_T.$
%\end{align*}

\noindent Combining these two cases, we obtain $\Denv \lesssim \max\{1 / \sqrt{T} + \dperish, L_T\}$.

\noindent \textbf{Hindsight Envy}: $\Denvhind$ is trivially bounded above by $w_{max}\cdot L_T \lesssim L_T$ since, for any $t, t'$:
$$w_\theta(X_{t', \theta'}^{alg} - X_{t, \theta}^{alg}) \leq w_{max}(\overline{X} - \underline{X}) = w_{max} L_T.$$

\noindent \textbf{Efficiency}:  
 Let $t_0=\sup\{t: X_{t}^{alg} = \underline{X}\}$. Then:
% \begin{align*}
%     \Deff & = B - \sum_{t, \theta} N_{t, \theta} X_{t, \theta} = B - \sum_{t} N_t X_{t}^{alg} \\
%     & = B_{t_0} + \sum_{t < t_0} N_t X_{t}^{alg} + \PUA^{alg}_{< t_0} - \sum_{t} N_{t} X_{t}^{alg} \\
%     &=  B_{t_0} - \sum_{t \geq t_0} N_t X_{t}^{alg} + \PUA^{alg}_{< t_0} \\
%     & < N_{t_0} \overline{X} + \overline{N}_{> t_0} \underline{X} + \overline{P}_{t_0} - N_{t_0} \underline{X} - N_{> t_0} \overline{X} + \PUA^{alg}_{< t_0}\\
%     & = \underline{X}(\overline{N}_{> t_0} - N_{> t_0}) - (\overline{X} - \underline{X})(N_{> t_0} - N_{t_0}) + \overline{P}_{t_0} + \PUA^{alg}_{< t_0},
% \end{align*}
\begin{align*}
    \Deff & = B - \sum_{t, \theta} N_{t, \theta} X_{t, \theta}^{alg} = B - \sum_{t} N_t X_{t}^{alg} \\
    & = B_{t_0} + \sum_{t < t_0} N_t X_{t}^{alg} + \PUA^{alg}_{< t_0} - \sum_{t} N_{t} X_{t}^{alg} \\
    &=  B_{t_0} - \sum_{t \geq t_0} N_t X_{t}^{alg} + \PUA^{alg}_{< t_0} \\
    & < N_{t_0} \overline{X} + \overline{N}_{> t_0} \underline{X} + \overline{P}_{t_0} - N_{t_0} \underline{X} - N_{> t_0} \overline{X} + \PUA^{alg}_{< t_0}\\
    & = \underline{X}(\overline{N}_{> t_0} - N_{> t_0}) - (\overline{X} - \underline{X})(N_{> t_0} - N_{t_0}) + \overline{P}_{t_0} + \PUA^{alg}_{< t_0},
\end{align*}
where the inequality follows from $X_{t_0}^{alg} = \underline{X}$ and the threshold condition for allocating $\overline{X}$.

Noting that $\underline{X} \leq \Bav$ and $\overline{N}_{> t_0} - N_{> t_0} \leq 2 \conf_{t_0, T}^N$, we have:
% \begin{align}\label{eq:ub-term-1}
% \underline{X}(\overline{N}_{> t_0} - N_{> t_0})&\leq \Bav\cdot 
% 2\confn{t_0}{T} \notag \\
% &\leq 2\Bav\sqrt{2\tilde{c}|\Theta|\rho_{\max}^2\log(2T^2/\delta)}\min\left\{\sqrt{T},L_T^{-1}+\sqrt{TL_T^{-1}\dperish} \right\},
% \end{align}
\begin{align}\label{eq:ub-term-1}
 \underline{X}(\overline{N}_{> t_0} - N_{> t_0}) &\leq \Bav\cdot 
2\confn{t_0}{T} \notag \\
&\leq 2\Bav\sqrt{2\tilde{c}|\Theta|\rho_{\max}^2\log(2T^2/\delta)}  \cdot \min\left\{\sqrt{T},L_T^{-1}+\sqrt{TL_T^{-1}\dperish} \right\},
\end{align}
where the second inequality follows from \cref{lem:switching_point_last_perish}.

We loosely upper bound the second term by:
\begin{align}\label{eq:ub-term-2}
& - (\overline{X} - \underline{X})(N_{> t_0} - N_{t_0}) \leq (\overline{X} - \underline{X}) N_{t_0}  \leq L_T |\Theta| (\mu_{max} + \rho_{\max}).
\end{align}

Finally, consider $\overline{P}_{t_0} + \PUA_{< t_0}^{alg}$. By construction, $\overline{P}_{t_0} \leq \overline{P}_1 = \overline{\Delta}(\underline{X})$. To upper bound $\PUA_{< t_0}^{alg}$, we consider the process that allocates $\underline{X}$ in each period to all arrivals. Let $B_t(\underline{X})$ denote the quantity of remaining items under this process, and $\setB{t}(\underline{X})$ the set of remaining items. %Note that $\setB{t} \subseteq \overline{\mathcal{B}}_t$, since $\overline{\mathcal{B}}_t$ includes items that may have perished under this slow process. 
We use $\PUA_t(\underline{X})$ to denote the quantity of unallocated items that perish at the end of period $t$ under this process, and $\PUA_{< t}(\underline{X})$ those that perished before the end of period $t$. The following lemma allows us to tractably bound $\PUA_{< t_0}^{alg}$ via this process. We defer its proof to Appendix \ref{apx:proof_lower_alg_comparison}.
\begin{lemma}
\label{lem:lower_alg_comparison}
    For all $t \in [T]$, 
    \begin{enumerate}
    \item $\Balg{t} \subseteq \setB{t}(\underline{X})$
    \item $\PUA_t^{alg} \leq \PUA_t(\underline{X})$.
    \end{enumerate}
\end{lemma}
Using these two facts, we have: $
\PUA_{< t_0}^{alg} \leq \PUA_{< t_0}(\underline{X}) \leq \PUA_{\geq 1}(\underline{X}) \leq \overline{P}_1 = \overline{\Delta}(\underline{X}).$  Hence,
\begin{align}\label{eq:ub-term-3}
& \overline{P}_{t_0} + \PUA_{< t_0}^{alg} \leq 2 \overline{\Delta}(\underline{X}) \leq 2 \overline{N} \dperish  \leq 2 \mu_{max} (1+\sqrt{2|\Theta|\rho_{\max}^2\log(2T^2/\delta)}) T \dperish,
\end{align}
where the second inequality follows from $\overline{\Delta}(\underline{X}) \leq B-\overline{N}\underline{X} = B-\overline{N}\left(\frac{B}{\overline{N}}-\dperish\right) = \overline{N}\dperish$, and the last inequality uses the definition of $\overline{N}$.
Putting bounds \eqref{eq:ub-term-1}, \eqref{eq:ub-term-2} and \eqref{eq:ub-term-3} together, we obtain:
% \begin{align*}
% \Deff \leq \Bav\cdot 
% 2\confn{t_0}{T} + L_T |\Theta| (\mu_{max} + \rho_{\max}) \\+ 2 \mu_{max} (1+\sqrt{2|\Theta|\rho_{\max}^2\log(2T^2/\delta)}) T \dperish.
% \end{align*}
\begin{align*}
 \Deff \leq &2\Bav\sqrt{2\tilde{c}|\Theta|\rho_{\max}^2\log(2T^2/\delta)}  \cdot \min\left\{\sqrt{T},L_T^{-1}+\sqrt{TL_T^{-1}\dperish} \right\}  + L_T|\Theta|(\mu_{max}+\rho_{\max})\\ &+ 2 \mu_{max} (1+\sqrt{2|\Theta|\rho_{\max}^2\log(2T^2/\delta)})T \dperish.
\end{align*}
Using the fact that $L_T = o(1)$, we obtain the final bound on efficiency.
\end{rproofof}

\begin{remark}
While \HopeGuardrailPerish achieves the optimal guarantees for arbitrary $\sigma$ and perishing distributions, the perishing-agnostic guardrail algorithm of \citet{sinclair2021sequential} is optimal in a special case of our setting. Namely, suppose the decision-maker knows the hindsight optimal ordering, and that there exist deterministic constants $\overline{N}$ and  $(\overline{P}_{<t}, \underline{N}_{<t}, t \geq 2)$, such that
% \begin{align}\label{eq:offset-relax}
% \mathbb{P}\left(P_{<t} \leq \overline{P}_{<t}, N_{<t} \geq \underline{N}_{<t}, N \leq \overline{N}\ \forall  \ t \geq 2\right) \geq 1-\delta \quad \text{and} \quad 
% \frac{\overline{P}_{<t}}{B} \leq \frac{\underline{N}_{<t}}{\overline{N}} \quad \forall \ t \ \geq 2,
% \end{align}
\begin{equation}\label{eq:offset-relax}
\begin{aligned}
&\mathbb{P}\left(
P_{<t} \leq \overline{P}_{<t},
N_{<t} \geq \underline{N}_{<t},
N \leq \overline{N}
\ \forall \ t \geq 2
\right)
\geq 1-\delta, 
&\text{ and } \quad \frac{\overline{P}_{<t}}{B}
\leq
\frac{\underline{N}_{<t}}{\overline{N}}
\quad \forall \ t \geq 2.
\end{aligned}
\end{equation}
for $\delta \in (0,1)$. Under the above condition, offset-expiry holds with probability at least $1-\delta$. Hence, \mbox{$\dperish = 0$} over the sample paths for which \eqref{eq:offset-relax} holds. As a result, the baseline allocation $\underline{X} = \frac{B}{\overline{N}}$ is feasible over these sample paths. Since allocating aggressively relative to $\underline{X}$ cannot increase the amount of goods that spoil, it follows that the algorithm of \citet{sinclair2021sequential} is optimal in this setting. While the assumption that the hindsight optimal ordering is known is unlikely to hold in practical settings where the perishing times are stochastic, this special case further highlights the two main features that render this setting challenging: $(i)$ the inherent aggressiveness of the perishing process, and $(ii)$ the unknown perishing order.
\end{remark}

%% file: figures/intro_tikz/phase_transition.tex
\tikzset{every picture/.style={line width=0.75pt}} %set default line width to 0.75pt        

\begin{tikzpicture}[x=0.75pt,y=0.75pt,yscale=-1,xscale=1]

%Straight Lines [id:da6027004850571227] 
\draw    (219,259) -- (219,7) ;
\draw [shift={(219,5)}, rotate = 90] [color={rgb, 255:red, 0; green, 0; blue, 0 }  ][line width=0.75]    (10.93,-4.9) .. controls (6.95,-2.3) and (3.31,-0.67) .. (0,0) .. controls (3.31,0.67) and (6.95,2.3) .. (10.93,4.9)   ;
%Straight Lines [id:da7927744524480898] 
\draw    (219,259) -- (479,259) ;
\draw [shift={(481,259)}, rotate = 180] [color={rgb, 255:red, 0; green, 0; blue, 0 }  ][line width=0.75]    (10.93,-4.9) .. controls (6.95,-2.3) and (3.31,-0.67) .. (0,0) .. controls (3.31,0.67) and (6.95,2.3) .. (10.93,4.9)   ;
% Horizontal guide line
\draw    (219,132) -- (470,132) ;

% Dashed curved line in gray
\draw[dashed, draw=gray, opacity=0.6] (272.62,12.7) .. controls (277.53,50.67) and (260.21,129.34) .. (354.5,132) ;
\draw[->, draw=gray, opacity=0.6, dashed] 
    (272.29,10.43) -- ++(0,0) node[anchor=south east] {} 
    coordinate[pos=0] (tmp) 
    -- ++(0,0); % arrowhead preserved

% Vertical line from curve end
\draw    (354.5,132) -- (354.5,7.43) ;

% Text Node
\draw (192.1,50) node [anchor=north west][inner sep=0.75pt]  [font=\large,rotate=-269.93] [align=left] {$L_T$};
% Text Node
\draw (425.55,260.06) node [anchor=north west][inner sep=0.75pt]  [font=\large] [align=left] {$\dperish$};
% Text Node
\draw (332.65,261.46) node [anchor=north west][inner sep=0.75pt]    {$T^{-1/2}$};
% Text Node
\draw (180.65,125.4) node [anchor=north west][inner sep=0.75pt]    {$T^{-1/2}$};
% Text Node
\draw (236.58,58.48) node [anchor=north west][inner sep=0.75pt]  [font=\normalsize] [align=left] {{\em High-Envy,} \\ {\em Low-Perishing} \\ (\cref{cor:main_theorem_interpret})};
% Text Node
\draw (337.58,186.48) node [anchor=north west][inner sep=0.75pt]  [font=\normalsize] [align=left] {{\em Low-Envy} \\ (\cref{cor:low-envy})};
% Text Node
\draw (387.58,58.48) node [anchor=north west][inner sep=0.75pt]  [font=\normalsize] [align=left] {{\em High-Envy,} \\ {\em High-Perishing} \\ (\cref{cor:high-envy-high-perish})};
% Text Node
% \draw (301.58,30.48) node [anchor=north west][inner sep=0.75pt]  [font=\normalsize] [align=left] {{\em High-Envy} \\ {\em Low-Perishing} \\ \cref{cor:main_theorem_interpret}};

\end{tikzpicture}

%% file: final_arxiv_parts/special-cases.tex
\subsection{Special Cases}\label{sec:special-cases}

We conclude the section by instantiating our bounds for important special cases. {Recall, \Cref{thm:offset-expiry-iff} established that the proportional allocation $B/N$ is achievable} if and only if the process is offset-expiring. As a result, for general, non-offset-expiring perishing processes, one should expect $\dperish = \Omega(1)$. In this case, \Cref{cor:low-envy,cor:high-envy-high-perish} imply that counterfactual envy is constant, and inefficiency grows linearly in $T$. {In this section we restrict our attention to a special class of distributions, parameterized by $\delta \in (0,1)$, which we term {\it $\delta$-offset-expiring}, for which $\dperish$ vanishes with respect to $T$ under the hindsight optimal ordering. We formally define this class of processes below.} 

\begin{definition}[$\delta$-offset-expiry.]\label{def:delta-offset-expiry}
Fix $\delta \in (0,1)$. A process is \emph{$\delta$-offset-expiring} if:
\begin{align}\label{eq:delta-offset}
{
\Pr\left(\frac{P_{<t}}{B} \leq \frac{N_{<t}}{N}  \  \forall \, t \geq 2\right) \geq 1-\delta.
}
\end{align}
\end{definition}

For simplicity, in this section we assume $B = N = T$, with $N_t = 1$ for all $t$, almost surely.\footnote{At the cost of cumbersome algebra, one can relax this assumption and derive entirely analogous results.} In this case, \Cref{eq:delta-offset} reduces to the condition that $\mathbb{P}(P_{<t} \leq t-1 \ \forall \ t\geq 2) \geq 1-\delta$.

\subsubsection{The i.i.d. setting.} 
We introduce some additional notation. For $t\in[T]$, we let \mbox{$\mathcal{B}^{rand}_{< t} = \{b: \Pr\left(T_b < t\right) \in (0,1)\}$} and $\mathcal{B}^{det}_{<t} = \{b: \Pr\left(T_b < t\right) = 1\}$. In words, $\mathcal{B}^{det}_{<t}$ represents the set of items that perish before period $t$, almost surely.

% In this case, the offset-expiry condition reduces to sample paths over which $P_{<t} \leq t-1$ for all $t\geq 2$. 

\Cref{prop:nec-cond-offset,prop:basic-conditions}  respectively establish necessary and sufficient conditions for $\delta$-offset-expiry to hold when the perishing distribution is i.i.d., for non-trivial values of $\delta$.

\begin{proposition}\label{prop:nec-cond-offset}
Suppose there exists $t \geq 2$ such that $\mathbb{E}[P_{<t}] > t - 1$. If $\mathcal{B}^{rand}_{<t} = \emptyset$, $\delta$-offset-expiry cannot be satisfied for any value of $\delta \in (0,1)$. Else, $\delta$-offset-expiry cannot be satisfied for $\delta < \frac{1}{2} -{\stdpt^{-3}} \cdot T$.  
\end{proposition}

\begin{proposition}\label{prop:basic-conditions}
Suppose $\mathbb{E}[P_{<t}] \leq t-1$ for all $t\geq 2$. Then, the perishing process is $\delta$-offset-expiring for any $\delta \geq \sum_{t=2}^T\min\left\{\left(\frac{\stdpt}{t-\mathbb{E}[P_{<t}]}\right)^2, \exp\left(-\frac{{2}(t-\mathbb{E}[P_{<t}])^2}{|\mathcal{B}^{rand}_{<t}|}\right) \right\}\mathds{1}\{|\mathcal{B}^{rand}_{<t}| > 0\}.$
\end{proposition}

\Cref{prop:nec-cond-offset,prop:basic-conditions} together illustrate how the mean and variability of the perishing distribution jointly impact offset-expiry.  In particular, \cref{prop:basic-conditions}  states that {\it either} the expected lag between demand and perishing, given by $t-\mathbb{E}[P_{<t}]$, must be large, {\it or} the coefficient of variation (CV) with respect to the random lag process $t-P_{<t}$, must be small. This then reduces to a bound on the variance of $P_{<t}$ in settings where perishing closely tracks demand.

Additionally, note that the necessary condition in \Cref{prop:nec-cond-offset} fails to hold for one of the most standard models of perishing: geometrically distributed perishing with parameter $1/T$ (that is, a constant fraction of items perish in each period). This highlights that one of the most popular models in the literature is, in a sense, far too pessimistic; for this setting, there is no hope of achieving low envy and efficiency with high probability. \Cref{thm:geometric_perishing} below shows how the perishing rate under this popular distribution affects offset-expiry and the perishing-induced loss $\dperish$.

\begin{proposition}
\label{thm:geometric_perishing}
Suppose $T_b \sim \emph{Geometric}(p)$ for all $b \in \mathcal{B}$, with $p =o(1/T)$.
Then, the perishing process is $\delta$-offset-expiring for any $\delta \geq 2 \log T \cdot \frac{Tp}{(1-Tp)^2}$. Moreover,
$\underline{X} \geq 1-3Tp-\frac{\log(3 \log(T)/\delta)}{T}$ for any ordering $\sigma$.  As a result,
\[
\dperish \leq 3Tp + \frac{\log(3 \log(T)/\delta)}{T}.
\]
\end{proposition}

\begin{figure}
\centering
    {
    \includegraphics[scale=0.4]{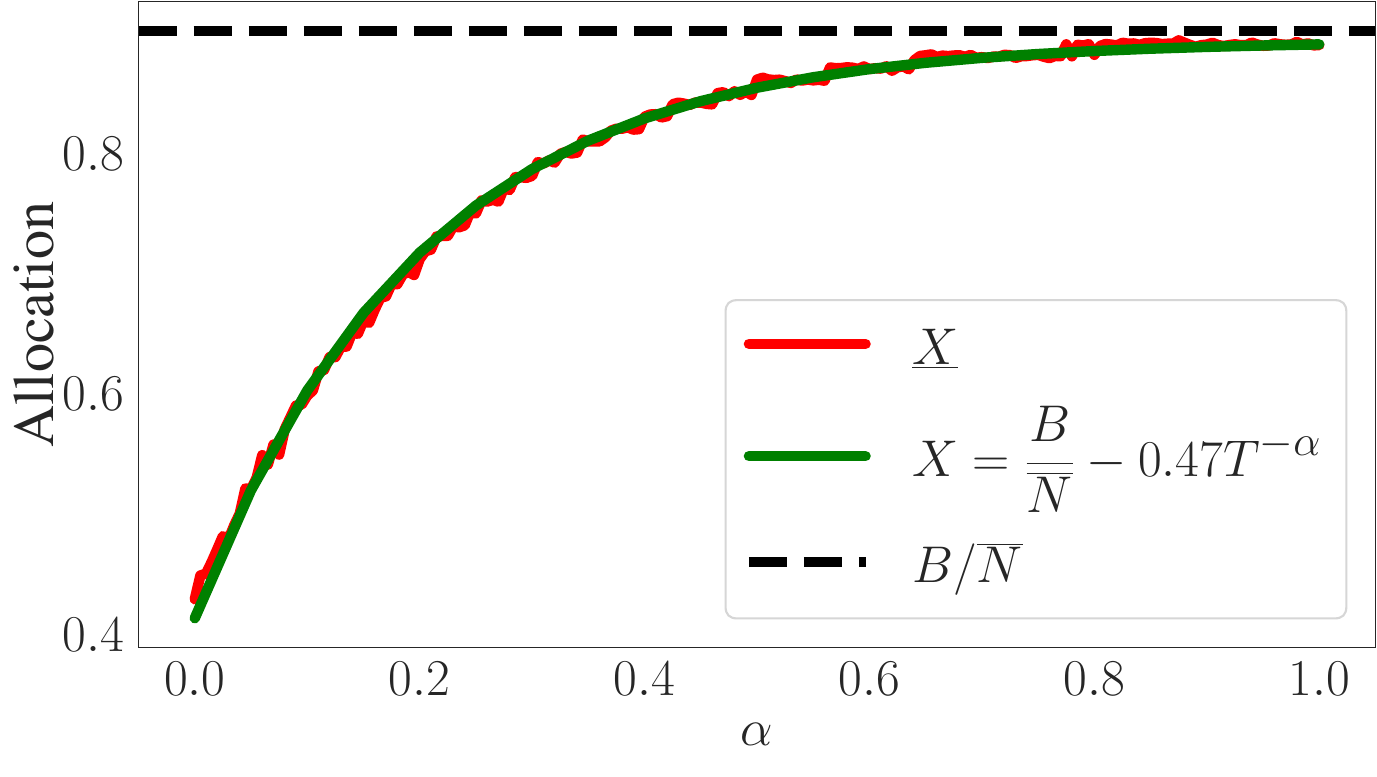}
    }
    \caption{
    Maximum feasible allocation $\underline{X}$ vs. $\alpha$, for $T_b \sim \emph{Geometric}(T^{-(1 + \alpha)})$, $B = 300$, $T = 150$, $N_t$ drawn from a truncated normal distribution $\mathcal{N}(2, 0.25)$, and $\delta = 1/T$. Here, $\underline{X}$ was calculated via line search, with Monte Carlo simulation used to estimate $\Delta(X)$ for each value of $X$. The dashed line represents the ``naive'' allocation $B / \overline{N} = 0.89$ which ignores possible perishing, and the second line is the curve of best fit to $\underline{X}$.
    \label{fig:geometric_dperish} }
\end{figure}

\cref{thm:geometric_perishing} establishes that, in the worst case,  $\underline{X}$ decays {\it linearly} in the rate $p$ at which goods spoil. This highlights the extent to which the commonly used (and practically pessimistic) geometric model limits both the kinds of perishable goods selected before allocation, as well as the rate at which a decision-maker can allocate goods. Letting $p = T^{-(1+\alpha)}$, $\alpha \in (0,1)$, \cref{thm:geometric_perishing} implies that $\dperish$ is on the order of $T^{-\alpha}$. Alternatively, if a decision-maker wants no more than $T^{-\alpha}$ loss relative to the proportional allocation, \cref{thm:geometric_perishing} provides an upper bound of $T^{-(1+\alpha)}$ on $p$.

We illustrate the scaling of $\underline{X}$ under geometric perishing numerically in \Cref{fig:geometric_dperish}, for an instance where the decision-maker also faces demand uncertainty. We observe that $\underline{X}$ is concave increasing in $\alpha$, and that the lower bound in \Cref{thm:geometric_perishing} provides a good fit for $\underline{X}$, even under demand uncertainty. For $\alpha$ close to 0  (i.e., $p \sim 1/T$), $\underline{X} \approx 0.4$, less than half of the ``naive'' no-perishing allocation, $B / \overline{N} =0.89$. For $\alpha = 1$, $\underline{X} \approx B / \overline{N}$. Note that, for $\alpha > 1$, $\underline{X}$ is limited by the confidence bound $\log(3\log T/\delta)/T$ in \cref{thm:geometric_perishing}. Plugging the lower bound on $\delta$ into this term, this implies that, even under deterministic demand, \mbox{$\dperish = \Omega(1/T)$}.

\subsubsection{The non-i.i.d. setting.} 
When the perishing distribution is not identical across different units, the quality of the allocation schedule $\sigma$ can have a significant impact on the perishing-induced loss $\dperish$. 
\cref{prop:gen_distribution_example} illustrates how the allocation schedule and perishing distribution come together in determining $\dperish$.  
\begin{proposition}
\label{prop:gen_distribution_example}
    Suppose there exists $\alpha \in (0,1)$ such that:
        \begin{enumerate}
        \item $\Exp{T_b} > \min\left\{T, \lceil \frac{\sigma(b)}{1 - T^{-\alpha}} \rceil\right\}$
        \item $\sum_b \frac{\Var{T_b}}{\left( \Exp{T_b} - \min\{T, \lceil\frac{\sigma(b)}{1-T^{-\alpha}}\rceil \} \right)^2} \leq \frac{1}{2}T^{1-\alpha}$.
    \end{enumerate}
    Then, for any $\delta \geq 3 \log(T) e^{-\frac18T^{1-\alpha}}$, the process is $\delta$-offset-expiring{, and $\underline{X} \geq 1 - T^{-\alpha}$.} As a result, \mbox{$\dperish  \leq T^{-\alpha}$}.
\end{proposition}
\cref{prop:gen_distribution_example} highlights the two key components that determine the baseline allocation: $(i)$ the lag between the expected perishing time $\Exp{T_b}$ and the expected allocation time $\min\{T, \lceil\frac{\sigma(b)}{1-T^{-\alpha}}\rceil \}$ (which we colloquially refer to as ``room to breathe''), and $(ii)$ the variance of the perishing time. Specifically, Condition 1 implies that, if $\Exp{T_b}$ is low, it must be that the item is allocated early on in the horizon (i.e., $\sigma(b)$ is low). This encodes the ``race against time'' intuition that is typically held around perishing. Condition 2 can be viewed as an upper bound on the cumulative adjusted {coefficient of variation} of the perishing process. High-variance perishing times and smaller ``room to breathe'' push $\alpha$ down, resulting in a lower maximum feasible allocation rate, and consequently higher perishing-induced loss. {In \cref{sec:experiments_other} we build on this intuition to propose heuristic schedules with low $\dperish$.}

%% file: final_arxiv_parts/experiments.tex
\section{Numerical Experiments}
\label{sec:experiments}

In this section we study the practical performance of \HopeGuardrailPerish via an extensive set of numerical experiments. We first consider one of the most popular (and aggressive) models of perishing: geometrically distributed perishing times. We explore the dependence of the envy-efficiency trade-off on the perishing rate for this setting, leveraging the empirical trade-off curves to provide guidance on the choice of envy parameter $L_T$. We moreover compare the performance of \HopeGuardrailPerish to its perishing-agnostic counterpart \citep{sinclair2021sequential} on both synthetic and real-world datasets~\citep{keskin2022data}. Lastly, for the non-i.i.d. setting we numerically study choices of allocation schedules that achieve the two-fold desideratum of being practical for decision-makers and achieving small perishing-induced loss.\footnote{See \url{https://github.com/seanrsinclair/perishing_fairness} for the code database.}

\subsection{Description of Metrics and Benchmark Policies}

We evaluate \HopeGuardrailPerish using $\delta = \frac{1}{T}$, benchmarking it against three algorithms from the literature:
\begin{itemize}
    \item \HopeGuardrail (\cite{sinclair2021sequential}): a guardrail-based algorithm designed for settings without perishable resources;
    \item \StaticXlower: the algorithm that allocates $X_{t, \theta} = \underline{X}$ for all $t,\theta$, until it runs out of resources;
    \item \StaticProp: the algorithm that allocates $X_{t, \theta} = \frac{B}{\overline{N}}$ for all $t,\theta$, until it runs out of resources.
\end{itemize}

We report the performance of these algorithms along five metrics: $(i)$ expected counterfactual envy $\Exp{\Denv}$, $(ii)$ expected hindsight envy $\Exp{\Denvhind}$, $(iii)$ expected inefficiency $\Exp{\Deff}$, $(iv)$ expected spoilage $\Exp{\Spoilage}$, and lastly $(v)$ stockout probability
$\Exp{\Stockout}$, i.e., the proportion of replications for which the algorithm runs out of resources before the end of the time horizon. 

\subsection{Geometric Perishing: Synthetic Experiments}
\label{sec:experiments_geometric}

We first consider the setting in which each unit perishes independently with probability $p$ in each period, i.e., $T_b \sim \text{Geometric}(p)$, for all $b \in \mathcal{B}$. Throughout the section, we assume $|\Theta| = 1$. Since perishing times are identically distributed, the allocation order $\sigma$ does not have any impact on the performance of the algorithm; hence, we assume $\sigma$ is the identity ordering. To study the impact of the perishing rate on algorithm performance, we let $p = T^{-(1+\alpha)}$, varying $\alpha \in  \{0.1, 0.2, 0.25, 0.3\}$ (see Appendix \ref{sec:experiment_details} for additional values of $\alpha$).
We let $T = 100$, $B = 2T$, and draw demand $N_t$ from a truncated normal distribution $\mathcal{N}(2, 0.25^2)$, for $t \in [T]$.   All results are averaged over 150 replications.

\subsubsection*{Empirical trade-off curves.}

\addtolength{\textfloatsep}{-0.2in}
\begin{figure}
    \begin{subfigure}{.45\textwidth}
        \centering
        \includegraphics[width=\textwidth]{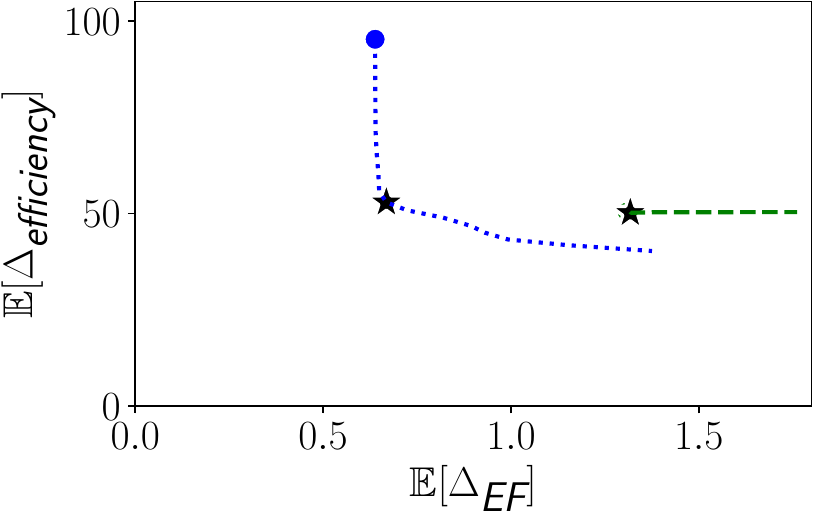}
        \caption{\centering $T_b \sim \emph{Geometric}(T^{-1.1}): \underline{X} = 0.70,$ $\dperish = 0.4, \Pr(\E_{OE}) = 0.89$.}\label{fig:tradeoff-1}
    \end{subfigure}
    \begin{subfigure}{.45\textwidth}
        \centering
        \includegraphics[width=\textwidth]{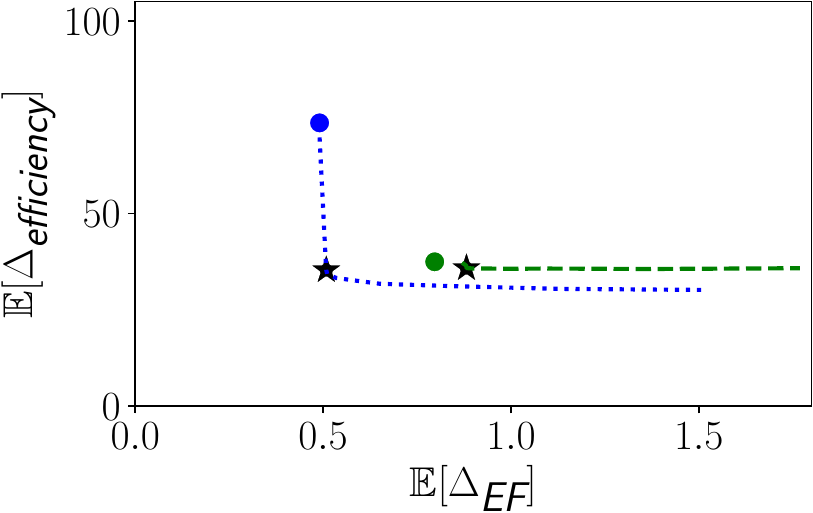}
        \caption{\centering $T_b \sim \emph{Geometric}(T^{-1.2}): \underline{X} = 0.84,$ $\dperish = 0.26, \Pr(\E_{OE}) = 0.97$}\label{fig:tradeoff-2}
    \end{subfigure}
        \begin{subfigure}{.45\textwidth}
        \centering
        \includegraphics[width=\textwidth]{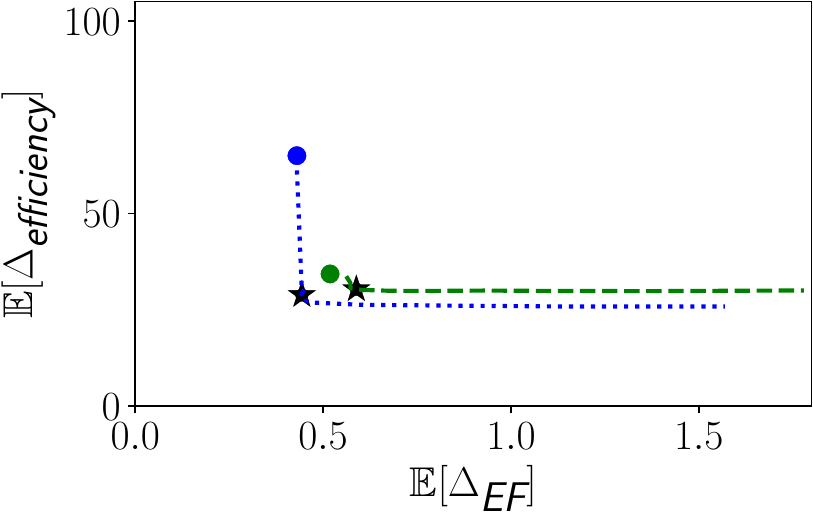}
        \caption{\centering $T_b \sim \emph{Geometric}(T^{-1.25}): \underline{X} = 0.89,$ $\dperish = 0.21, \Pr(\E_{OE}) = 0.99$}\label{fig:tradeoff-25}
    \end{subfigure}
    \begin{subfigure}{.45\textwidth}
        \centering
        \includegraphics[width=\textwidth]{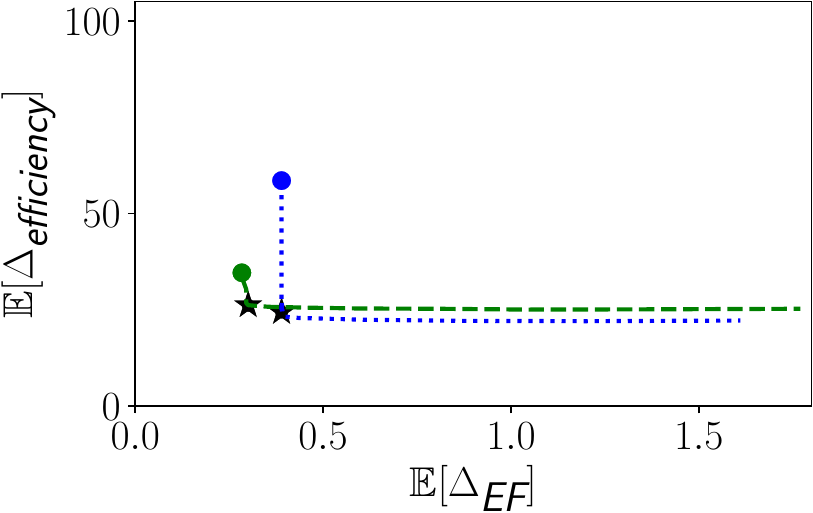}
        \caption{\centering $T_b \sim \emph{Geometric}(T^{-1.3}): \underline{X} = 0.93,$ $\dperish = 0.17, \Pr(\E_{OE}) = 0.99$}\label{fig:tradeoff-3}
    \end{subfigure}
    \begin{subfigure}{\textwidth}
        \centering
        \includegraphics[width = \textwidth]{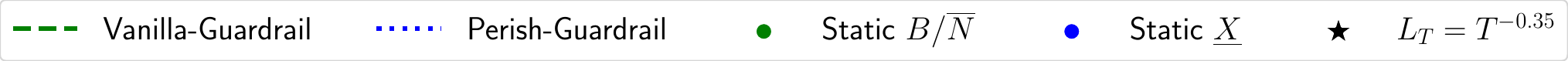}
    \end{subfigure}
    \caption{\centering Empirical trade-off between $\Exp{\Deff}$ and $\Exp{\Denv}$. The points on the trade-off curve correspond to increasing values of $L_T$, from left to right. \StaticProp and \StaticXlower respectively correspond to \HopeGuardrail and \HopeGuardrailPerish for $L_T = 0$.}
    \label{fig:geometric_tradeoff}
\end{figure}

We begin by numerically investigating the impact of the perishing rate on the envy-efficiency frontier, and use the empirical trade-off curves to provide guidance on how decision-makers should select $L_T$. For these instances, we compare \HopeGuardrailPerish and \HopeGuardrail \citep{sinclair2021sequential} with $L_T = T^{-\beta}$, $\beta \in \{0,0.05,0.1,\ldots,1\}$.

The empirical trade-off curves for each value of $\alpha$ can be found in \Cref{fig:geometric_tradeoff}. In each case, we report the probability that a sample path is offset-expiring, denoted by $\mathbb{P}(\mathcal{E}_{OE})$. For $\alpha \in \{0.1,0.2\}$ (\Cref{fig:tradeoff-1,fig:tradeoff-2}), \HopeGuardrail makes close to no gains in efficiency for any value of $L_T$. In these high-perishing settings, then, setting $L_T = 0$ is optimal (there is no trade-off), in stark contrast to the classic setting without perishability. As $\alpha$ increases (\Cref{fig:tradeoff-25,fig:tradeoff-3}), \HopeGuardrail attains small gains in efficiency, but plateaus very quickly. This yields the important insight that  \textbf{perishing-agnostic algorithms are {\it not} able to leverage unfairness to improve efficiency in the presence of perishable resources.} \HopeGuardrailPerish, on the other hand, sees extremely large gains in efficiency for a very small increase in envy, across all values of $\alpha$. For larger values of $L_T$, we observe marginal gains in efficiency; the horizontal asymptote observed across all values of $\alpha$ is precisely the cumulative perishing-induced loss, $T\dperish$ {(\cref{thm:lower_bound})}. Moreover, the vertical asymptote across all plots corresponds to the unavoidable loss due to demand uncertainty (\cref{thm:lower_bounds}).

Note that, for small values of $\alpha$, the \HopeGuardrailPerish empirical trade-off curve lies to the left of that of \HopeGuardrail, i.e., it achieves lower counterfactual envy across all values of $L_T$. As we will see below, this is due to the fact that \HopeGuardrail achieves extremely high stockout rates. This effect is diminished as $\alpha$ increases (i.e., the perishing rate decreases). As this happens, both curves move down and to the left (and closer) as they achieve lower counterfactual envy and inefficiency due to spoilage. When the perishing rate is negligible ($\alpha=0.3$), the empirical trade-off curve of \HopeGuardrail is slightly to the left of that of \HopeGuardrailPerish; this is due to the loss incurred by our modified $\underline{X}$ construction, which always allocates less than $B/\overline{N}$ as a baseline. However, even when perishing is negligible \HopeGuardrailPerish is slightly more efficient than \HopeGuardrail, despite its baseline allocation being lower. This runs counter to the intuition that \HopeGuardrail should be more efficient since it has a higher baseline allocation. The reason for this is the difference in the two algorithms' threshold allocation decisions. Our algorithm, \HopeGuardrailPerish, allocates $\overline{X} = \underline{X} + L_T$ if it forecasts that it has enough budget remaining to allocate $\overline{X}$ in period $t$, and $\underline{X}$ onwards. On the other hand, $\HopeGuardrail$ allocates $B/\overline{N} + L_T$ if it has enough budget remaining to allocate this high amount in period $t$, and $B/\overline{N}$ in all future periods. Since $B/\overline{N} > \underline{X}$, \HopeGuardrail depletes its budget faster than \HopeGuardrailPerish whenever they both allocate the lower guardrail. Hence, \HopeGuardrailPerish is able to allocate aggressively {more} frequently than \HopeGuardrail, which results in improved efficiency and decreased spoilage. 

These empirical trade-off curves moreover help to provide guidance on the choice of $L_T$. In particular, across all experiments, the cusp of the trade-off curve lies at $L_T \sim T^{-0.35}$, which is {\it larger} than no-perishing cusp of $L_T \sim T^{-1/2}$ \citep{sinclair2021sequential}. This is due to the fact that our baseline allocation is significantly lower to avoid perishing-induced stockouts; hence, in order to recover from this inefficiency, $L_T$ must be higher. We use this observation in the following experiments, comparing the performance of \HopeGuardrailPerish and \HopeGuardrail for $L_T \sim T^{-0.35}$.

\subsubsection*{Benchmark comparisons.}
\begin{figure}[t]
\begin{subfigure}[b]{\textwidth}
\centering
\includegraphics[width=.8\textwidth]{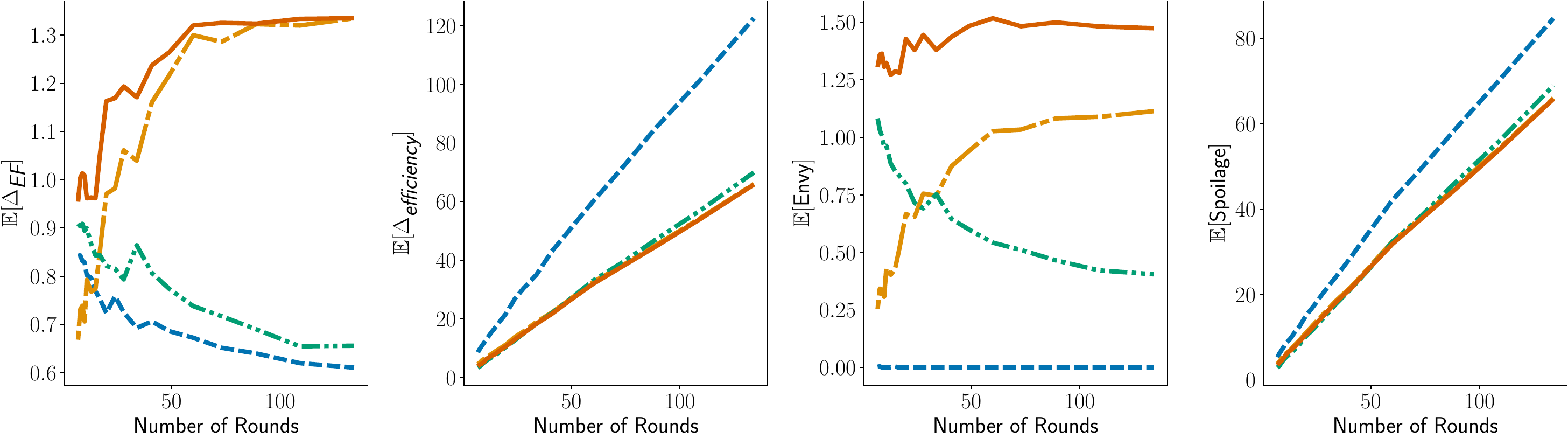}
\caption{\small $T_b \sim \emph{Geometric}(T^{-1.1})$}\label{fig:geom-0.1}
\end{subfigure}
\begin{subfigure}[b]{\textwidth}
\centering
\includegraphics[width=.8\textwidth]{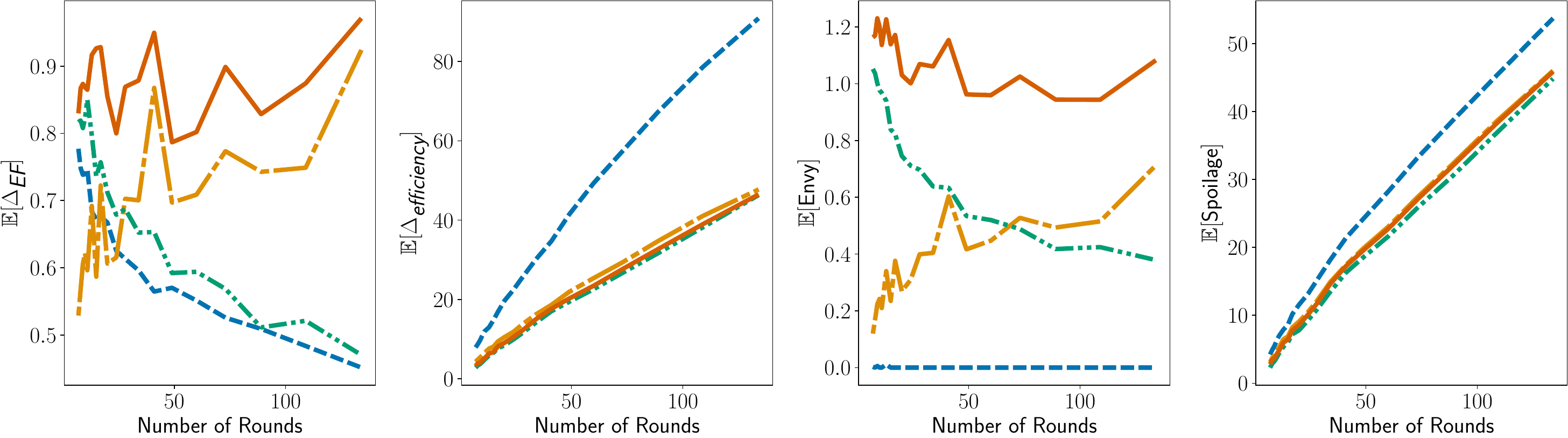}
\caption{\small $T_b \sim \emph{Geometric}(T^{-1.2})$}\label{fig:geom-0.2}
\end{subfigure}
\begin{subfigure}[b]{\textwidth}
\centering
\includegraphics[width=.8\textwidth]{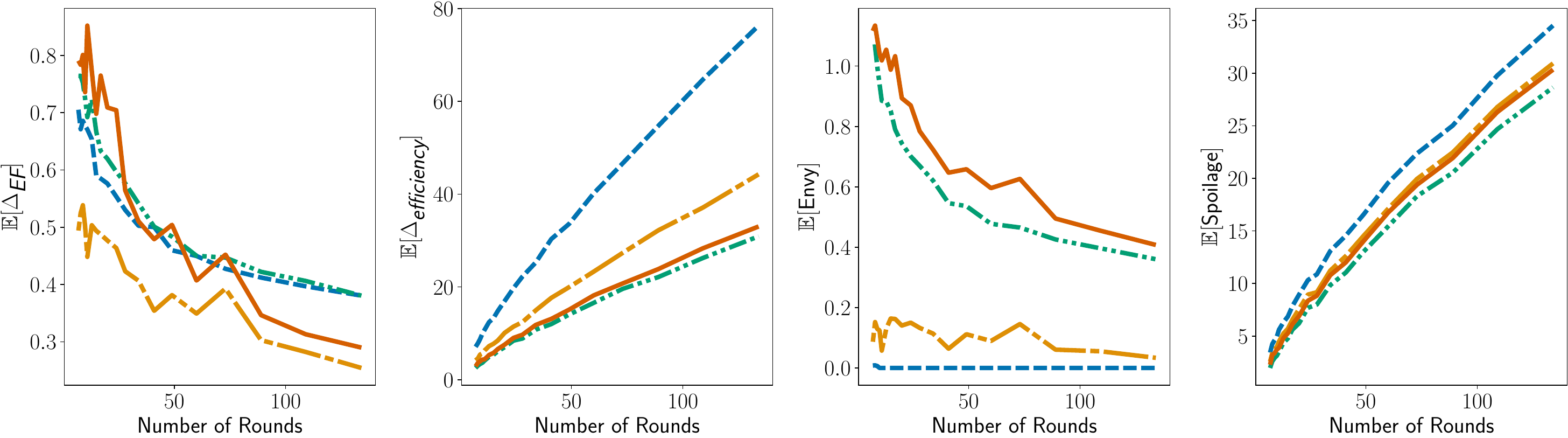}
\caption{\small $T_b \sim \emph{Geometric}(T^{-1.3})$}\label{fig:geom-0.3}
\end{subfigure}
\begin{subfigure}{\textwidth}
\centering
\includegraphics[width=.8\textwidth]{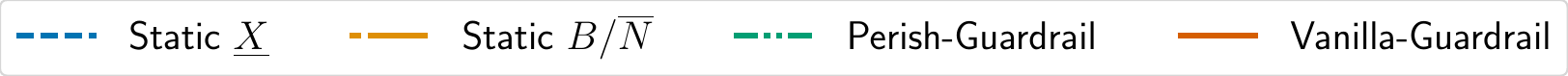}
\end{subfigure}
\caption{\centering Algorithm comparison across $\Exp{\Denv}, \Exp{\Deff}, \Exp{\Denvhind}$, and $\Exp{\Spoilage}$,  for $\alpha \in \{0.1,0.2,0.3\}$.}
\label{fig:geometric_simulations}
\end{figure}

Next, we compare the performance of all algorithms along our five metrics as the number of rounds $T$ grows large, for \mbox{$\alpha \in \{0.1,0.2,0.3\}$}. Based on our prior insights, we use \mbox{$L_T = T^{-0.35}$} for both \HopeGuardrailPerish and \HopeGuardrail. Our results can be found in \Cref{fig:geometric_simulations}.
We identify three regimes:
\begin{itemize}
    \item {\bf High Perishing} ($\alpha = 0.1$, \Cref{fig:geom-0.1}): While unfairness (as measured by $\Exp{\Denv}$ and $\Exp{\Denvhind}$) is decreasing in $T$ under \HopeGuardrailPerish and \StaticXlower, \HopeGuardrail and \StaticProp perform remarkably poorly along these two metrics. This is due to the fact that these latter algorithms fail to account for the unavoidable perishing loss, resulting in an extremely high stockout probability, as illustrated in \cref{tab:stockout_prob_table}, for $T = 150$. In contrast, the two perishing-aware algorithms rarely run out of resources. This underscores the importance of modifying the baseline guardrail $\underline{X}$, which was specifically constructed to avoid stockouts due to unavoidable perishing. 

    Comparing \StaticXlower to \HopeGuardrailPerish, our results also demonstrate that, in this high-perishing regime, the strategy of cautiously over-allocating by $L_T$ comes at a significant reduction in inefficiency $\Exp{\Deff}$, at close to no increase in counterfactual envy $\Exp{\Denv}$.
    
    \item {\bf Medium Perishing} ($\alpha = 0.2$, \Cref{fig:geom-0.2}): Though we observe similar trends as when $\alpha = 0.1$, all algorithms perform better across the board. \HopeGuardrailPerish exhibits slightly higher efficiency and lower spoilage than its perishing-agnostic counterpart in this regime, since it satisfies the threshold condition more frequently, as described above. Still, in \cref{tab:stockout_prob_table} we see that the perishing-agnostic algorithms run out of resources in over 50\% of replications. 
    
    \item {\bf Low Perishing} ($\alpha = 0.3$, \Cref{fig:geom-0.3}): For this smaller perishing rate, \HopeGuardrail stocks out significantly less frequently. Putting this together with the fact that $B/\overline{N} > \underline{X}$, this explains the fact that it has lower counterfactual envy than \HopeGuardrailPerish. However, along all other metrics \HopeGuardrailPerish improves upon \HopeGuardrail. The improvements in efficiency and spoilage are due to the same effects as described above; moreover, our algorithm improves upon \HopeGuardrail on $\Exp{\Denvhind}$ since it never stocks out.
\end{itemize}

\begin{table}[t]
\centering{
\begin{tabular}{c|rrrr}
\toprule
$\alpha$ & \StaticProp & \StaticXlower & \HopeGuardrail & \HopeGuardrailPerish \\
\midrule
$0.1$ & $0.99 \pm 0.02$ & $0.00 \pm 0.0$ & $1.00 \pm 0.0$ & $0.11 \pm 0.06$ \\
$0.2$ & $0.63 \pm 0.1$ & $0.00 \pm 0.0$ & $0.68 \pm 0.09$ & $0.03 \pm 0.03$ \\
$0.25$ & $0.17 \pm 0.07$ & $0.00 \pm 0.0$ & $0.22 \pm 0.08$ & $0.00 \pm 0.0$ \\
$0.3$ & $0.03 \pm 0.04$ & $0.00 \pm 0.0$ & $0.06 \pm 0.05$ & $0.00 \pm 0.0$ \\
\bottomrule
\end{tabular}
}
\caption{Comparison of stockout probabilities, for $T = 150$, $T_b \sim \text{Geometric}(T^{-(1+\alpha)})$. The second number in each cell corresponds to $95\%$ confidence intervals. 
\label{tab:stockout_prob_table} }
\end{table}

Overall, our results highlight the robustness of \HopeGuardrailPerish to perishability, as it is able to achieve similar if not improved performance compared to \HopeGuardrail in settings where there is limited perishing, with vastly superior performance in high-perishing settings.

\begin{table}
\centering{
\begin{tabular}{c|rrrrr}
\toprule
Algorithm & $\Exp{\Denv}$ & $\Exp{\Denvhind}$ & $\Exp{\Spoilage}$ & $\Exp{\Stockout}$ & $\Exp{\Deff}$\\
\midrule
\StaticProp & $1.18 \pm 0.01$ & $1.03 \pm 0.0$ & $346.4 \pm 2.6$ & $1.0 \pm 0.0$ & $444.6 \pm 2.7$ \\
\StaticXlower & $0.60 \pm 0.01$ & $0.0 \pm 0.0$ & $475.9 \pm 2.7$ & $0 \pm 0.0$ & $605.5 \pm 3.0$\\
\HopeGuardrail & $1.17 \pm 0.01$ & $1.44 \pm 0.0$ & $341.4 \pm 3.0$ & $1.0 \pm 0.0$ & $343.5 \pm 2.9$\\
\HopeGuardrailPerish & $0.78 \pm 0.04$ & $0.42 \pm 0.05$ & $372.2 \pm 3.2$ & $0.39 \pm 0.09$ & $372.7 \pm 3.2$ \\
\bottomrule
\end{tabular}
}
\caption{ Performance of the different algorithms (for $L_T = T^{-0.35}$) on the ``ginger'' dataset in \citet{keskin2022data}.  The second number in each cell corresponds to 95\% confidence intervals.
\label{tab:geometric_real_data} }
\end{table}

% \ifdefined\arxiv \newpage \else\fi

\subsection{Geometric Perishing: Real-World Instance}

We next investigate the performance of our algorithm using the real-life dataset provided by \citet{keskin2022data}, which tracks daily demand, replenishments, and perishing of ginger at a leading supermarket chain in China over $T = 365$ days. We assume each unit's perishing time is geometrically distributed with parameter $p$. Fitting this geometric perishing model to the data, we obtain an estimate of $p = 0.00224$. We similarly fit a truncated normal distribution to the demand data, obtaining $\mathcal{N}(3.2, 1.85^2)$. Finally, we set the total budget to $B = 365 \cdot 3.2 = 1168$, corresponding to the total expected demand over the horizon. (See Appendix \ref{sec:experiment_details} for additional details on the data calibration.)
For these inputs, $B / \overline{N} = 0.89$ and $\underline{X} = 0.46$. Under these parameters, the offset-expiry condition is only satisfied 65.2\% of the time.  Given this aggressive perishing, perishing-agnostic algorithms are expected to perform particularly poorly.

In \cref{tab:geometric_real_data} we compare the performance of the different algorithms.  We observe the following:
\begin{itemize}
    \item As conjectured, \StaticProp and \HopeGuardrail stock out on 100\% of replications since they fail to account for endogenous perishing. The high stockout probabilities of \StaticProp and \HopeGuardrail lead to high unfairness (vis-\`a-vis hindsight and counterfactual envy), since later arrivals receive allocations of zero. In contrast, \StaticXlower never stocks out. \HopeGuardrailPerish achieves a higher stockout rate of 40\%, due to its more aggressive allocations. Still, our algorithm's counterfactual envy and hindsight envy are over 30\% and 70\% lower, respectively, than that of \HopeGuardrail.
    \item \HopeGuardrailPerish allocates approximately 10\% fewer goods than \HopeGuardrail. It is notable, however, that it is more efficient than \StaticProp; this highlights that naively allocating more aggressively need not always generate {efficiency} gains.
\end{itemize}
Overall, we see that even when offset-expiry fails to hold over a significant fraction of sample paths, for small losses in efficiency our algorithm makes major gains in fairness relative to perishing-agnostic algorithms. This highlights the robustness of our algorithm to aggressive perishing.

\subsection{Non-i.i.d. Perishing}
\label{sec:experiments_other}

\Cref{thm:upper_bound_perishing} established that the theoretical performance of \HopeGuardrailPerish improves for orderings $\sigma$ under which $\dperish$ is lower. In practice, however, $\dperish$ has an extremely subtle impact on our algorithm's efficiency. Namely, while lower $\dperish$ implies a higher $\underline{X}$, this also results in a lower value of $\Pupper_t$ in each round $t$. The net effect on whether the threshold condition is satisfied in each period, then, remains unclear. In this section we study the impact of reasonable heuristic orderings on the practical performance of our algorithm.  Motivated by the insights of \cref{prop:gen_distribution_example}, we compare the following orderings:

\begin{itemize}
    \item {\bf Increasing Mean} ($\sMean$): Increasing order of $\Exp{T_b}$;
    \item {\bf Decreasing Coefficient of Variation} (CV) ($\sCV$): Decreasing order of $\STD{T_b} / \Exp{T_b}$. For fixed expected perishing time, this schedule allocates high-variance units earlier on. Conversely, for fixed variance, it allocates items according to the {\textbf{Increasing Mean}} schedule;
    \item {\bf Increasing Lower Confidence Bound} (LCB) ($\sLCB$): Increasing order of $\Exp{T_b} - 1.96 \STD{T_b}$. This ordering allocates items according to the lower bound of the 95\% confidence interval of the normal approximation to its perishing time. This lower bound is expected to be small if either the expected perishing time is small or the variance is large. 
\end{itemize}
We break ties randomly in all cases.

As in \cref{sec:experiments_geometric}, we draw the demands $N_t$ from a truncated normal distribution,   $\mathcal{N}(2, 0.25^2)$; we moreover let $T = 50$, $B = 100$, and $L_T = {T}^{-0.35}$. Finally, we consider two sets of perishing distributions:
\begin{itemize}
\item {\bf Instance 1: Front-loaded variability} \begin{align}
\label{eq:normal_example_A}
    T_b = \begin{cases}
        \text{Uniform}(T/2 - b/2, T/2 + b/2)  & \quad b \leq T\\
         T & \quad b > T
    \end{cases}
\end{align}
\item {\bf Instance 2: Back-loaded variability} \begin{align}
\label{eq:normal_example_B}
    T_b = \begin{cases}
        b+1  & \quad b \leq T \\
         \text{Uniform}(T, b+1) & \quad b > T
    \end{cases}
\end{align}
\end{itemize}
\cref{tab:uniform_example_A,tab:uniform_example_B} show the performance of our algorithm over these instances, for all three allocation schedules.  % 

\begin{table}
\centering
\begin{tabular}{c|rrrrrr}
\toprule
Order & $\Exp{\dperish}$ & $\Exp{\Denv}$ & $\Exp{\Denvhind}$ & $\Exp{\Spoilage}$ & $\Exp{\Stockout}$ & $\Exp{\Deff}$ \\ %& $\Exp{\Swap}$\\
\midrule
$\sMean$ & 
$0.10 \pm 0.004$ & $0.36 \pm 0.02$ & $0.51 \pm 0.02$ & $5.78 \pm 0.4$ & $0.04 \pm 0.04$ & $6.28 \pm 0.4$ \\ % & $6\% \pm 0.02\%$ \\
$\sCV$ / $\sLCB$ & $0.0 \pm 0.0$ & $0.44 \pm 0.02$ & $0.48 \pm 0.02$ & $1.24 \pm 0.1$ & $0.06 \pm 0.04$ & $1.79 \pm 0.1$ \\ % & $6\% \pm 0.02\%$\\
\bottomrule
\end{tabular}

\caption{\centering {Performance of \HopeGuardrailPerish for $L_T = T^{-0.35}$ on Instance 1 (\cref{eq:normal_example_A}).}
\label{tab:uniform_example_A} }
\end{table}

\begin{table}
\centering

\begin{tabular}{c|rrrrrr}
\toprule
Order & $\Exp{\dperish}$ & $\Exp{\Denv}$ & $\Exp{\Denvhind}$ & $\Exp{\Spoilage}$ & $\Exp{\Stockout}$ & $\Exp{\Deff}$ \\ %$\Exp{\Swap}$\\
\midrule
$\sMean$ / $\sLCB$ & $0.0 \pm 0.0$ & $0.41 \pm 0.02$ & $0.46 \pm 0.02$ & $0.0 \pm 0.0$ & $0.1 \pm 0.05$ & $0.56 \pm 0.07$ \\ % & $0\% \pm 0.0$\\
$\sCV$ & $0.47 \pm 0.0$ &  $0.51 \pm 0.05$ & $0.48 \pm 0.02$ & $48.3 \pm 0.08$ & $0.01 \pm 0.0$ & $48.7 \pm 0.07$ \\ % & $100\% \pm 0.0$\\
% Increasing LCB & $0.03 \pm 0.0$ & $0.43 \pm 0.05$ & $0.50 \pm 0.09$ & $0.0 \pm 0.0$ & $0.0 \pm 0.0$ & $2.3 \pm 0.1$\\
\bottomrule
\end{tabular}
\caption{\centering {Performance of \HopeGuardrailPerish for $L_T = T^{-0.35}$ on Instance 2 (\cref{eq:normal_example_B}).}
\label{tab:uniform_example_B} }
\end{table}

For the first instance, the {\bf Increasing Mean} schedule allocates the first $T$ items uniformly at random, ignoring the fact that, for $b \leq T$, as $b$ increases the item is more likely to perish earlier on in the horizon. The {\bf Decreasing CV} / \textbf{Increasing LCB} schedules, on the other hand, are identical: they allocate the first $T$ resources in decreasing order of $b$, and allocate the remaining uniformly at random.
Notably, the {\bf Decreasing CV} / \textbf{Increasing LCB} order achieves $\Exp{\dperish} = 0$, i.e., $\underline{X} = B/\overline{N}$, as in the no-perishing setting. (Note that $\dperish = 0$ implies that this is an optimal ordering.) Since its baseline allocation is higher it results in 78\% less spoilage than the {\textbf{Increasing Mean} order, and a 71\% decrease in inefficiency. However, this order performs slightly worse with respect to counterfactual envy and stockouts: this is again due to the more aggressive allocations.

For the second instance, the {\bf Increasing Mean} and {\bf Increasing LCB} schedules are identical: they allocate items lexicographically. The {\bf Decreasing CV} schedule, on the other hand, allocates the last $T$ items (in increasing order of $b$) before the first $T$ resources, since $\STD{T_b}=0$ for all $b \leq T$. In this setting, the first schedule is optimal with respect to perishing-induced loss, with $\Exp{\dperish} = 0$. This more aggressive allocation results in a 10\% stockout rate (versus 1\% for the {\bf Decreasing CV} schedule), but outperforms the {\bf Decreasing CV} order across all other metrics. This is intuitive as the number of errors in this latter, clearly bad order results in $\Exp{\dperish} = 0.47$, approximately $50\%$ of the baseline allocation $B / \overline{N}$. The algorithm then incurs both high inefficiency and spoilage.

These results indicate that the {\bf Increasing LCB} schedule is both a practical and robust candidate allocation order as it hedges against the inherent variability of the perishing process. 

%% file: final_arxiv_parts/conclusion.tex
\section{Conclusion}

This paper considers a practically motivated variant of the canonical problem of online fair allocation wherein a decision-maker has a budget of {\it perishable} resources to allocate {\it fairly} and {\it efficiently} over a fixed time horizon. Our main insight is that perishability fundamentally impacts the envy-efficiency trade-off derived for the no-perishing setting: while a decision-maker can arbitrarily sacrifice on envy in favor of efficiency in this latter setting, this is no longer the case when there is uncertainty around items' perishing times. We derive strong lower bounds to formalize this insight, which are a function of both the quality of the decision-maker's prediction over perishing times, as well as the inherent aggressiveness of the perishing process. We moreover design an algorithm that achieves these lower bounds; this algorithm relies on the construction of a baseline allocation that accounts for the unavoidable spoilage incurred by any online algorithm. From a technical perspective, the main challenge that the perishing setting presents is that the uncertainty around the quantity of resources that spoil in the future is {\it endogenous}, in contrast to the {\it exogenous} uncertainty on the number of arrivals in the classical setting. Deriving tight bounds on spoilage (both for our lower bounds as well as in the design of our algorithm) relied on the ``slow allocation'' construction, which rendered the highly coupled process amenable to tractable analysis. Our numerical experiments also demonstrate our algorithm's strong performance against state-of-the-art {\it perishing-agnostic} benchmarks.

From a practical perspective, our results provide concrete, intuitive guidance as to the composition of perishable donations that a food bank should accept in order to achieve a fair allocation (in the proportional sense). At a high level, the optimal mix of goods should ensure that cumulative perishing always lags behind cumulative demand. Moreover, we uncover the important --- and counterintuitive --- insight that guaranteeing a fair allocation in the presence of perishable goods may require the decision-maker to allocate significantly {\it less} than the static proportional benchmark. Finally, our numerical results suggest that an allocation schedule that ranks items according to high-probability {\it lower bounds} on their perishing times is a robust choice that appropriately hedges against the inherent variability of the perishing process. In particular, we observe that this allocation schedule outperforms natural alternatives such as ordering items according to expected perishing time.

In terms of future directions, our model assumes that the decision-maker allocates items according to a fixed allocation schedule. Though our results do not require that the perishing distribution be memoryless, allowing for time-varying / adaptive allocation schedules, though less practical, would improve our algorithm's performance in non-memoryless settings. This relates back to the question of deriving theoretical insights into the structure of optimal allocation schedules. Finally, though this paper considered exogenous {\it depletion} of the budget, a natural practical extension is one wherein $B$ evolves stochastically, accounting for external donations independent of the allocations made by the algorithm~\citep{onyeze2025sequential}.

%% file: final_arxiv_parts/appendix_perishing.tex
\section{\cref{sec:construction_x_lower} Omitted Proofs}\label{apx:proof-of-offset-expiry}

\subsection{Proof of \cref{thm:offset-expiry-iff}}

\begin{rproof}
We first argue that offset-expiry implies feasibility of $B/N$. Consider the allocation schedule which allocates goods in increasing order of perishing time (breaking ties arbitrarily), and is such that $X_{t,\theta} = B/N$ for all $t,\theta$, as long as there are resources remaining. {Noting that $(B/N) N_{< t}$ is precisely the cumulative allocation at the beginning of round $t$, this implies that we allocate (weakly) more than the number of goods with perishing time before round $t$ (i.e. $P_{< t}$).} Since we allocate goods in increasing order of perishing time, this also implies that no unit ever perishes under this sequence of allocations. Thus, the total allocation by the end of the horizon is $\frac{B}{N}\cdot N = B$, implying that $B/N$ is feasible.

We now argue that offset-expiry is {necessary} for $B/N$ to be feasible. To see this, consider the first period $t \geq 2$ for which $P_{< t}/B > N_{<t}/N$ (i.e., by the end of period $t-1$, there existed some unallocated goods that had perished). Then, the remaining budget at the start of period $t$ for any algorithm, denoted by $B_t^{alg}$, is:
\begin{align*}
B_t^{alg} \leq B-P_{<t} < B-N_{< t}\cdot \frac{B}{N} = N_{\geq t}\cdot \frac{B}{N},
\end{align*}
which implies that the remaining budget does not suffice to allocate $B/N$ to all arrivals from $t$ onwards. Hence, $B/N$ is not feasible.
\end{rproof}

\section{\Cref{sec:alg} Omitted Proofs}

\subsection{Tightness of bounds}
\label{apx:tightness}

Consider the random problem instance which achieves the lower bounds of \cref{thm:lower_bounds} with probability 1/2, and the lower bounds of \cref{thm:lower_bound} with probability 1/2. Putting these two bounds together, we have:
\begin{align*}
\Exp{\Denv} \gtrsim \dperish + 1/\sqrt{T}.
\end{align*}
By \cref{thm:upper_bound_perishing}, our algorithm achieves $\Exp{\Denv} \lesssim \max\{L_T, \dperish + 1/\sqrt{T}\}.$ Letting $L_T \lesssim \dperish + 1/\sqrt{T}$ then, our algorithm achieves this lower bound. We now argue that our algorithm is tight with respect to efficiency in this regime. Suppose $L_T = 0$. By \cref{thm:lower_bounds} and \cref{thm:lower_bound}, any online algorithm incurs:
\begin{align*}
\Exp{\Deff} \gtrsim T\dperish + \sqrt{T},
\end{align*}
which is achieved by our algorithm.

Consider now the regime in which $\Denv = L_T$, i.e., $L_T \gtrsim \dperish + 1/\sqrt{T}$. Again, randomizing between the two lower bounds, we have:
\begin{align}\label{eq:eff-lb}
\Exp{\Deff} \gtrsim T\dperish + \min\{\sqrt{T}, L_T^{-1}\}.
\end{align}

\noindent\textbf{Case 1:} $L_T^{-1} + \sqrt{T \dperish L_T^{-1}} \gtrsim \sqrt{T}$.
Here, our algorithm achieves $\Exp{\Deff} \lesssim \sqrt{T} + T\dperish.$ 

If $L_T^{-1} \gtrsim \sqrt{T}$, we achieve the bound in \eqref{eq:eff-lb}. Suppose now that $L_T^{-1} = o(\sqrt{T})$. Then, \eqref{eq:eff-lb} implies that $\Exp{\Deff} \gtrsim T\dperish + L_T^{-1}.$ We argue that, if $L_T^{-1} = o(\sqrt{T})$, then in this case $T\dperish \gtrsim \sqrt{T}$. $T\dperish$ then dominates both the lower bound in \eqref{eq:eff-lb}, as well as our upper bound, which gives us tightness. 

\medskip

\noindent\textbf{Case 2:} $L_T^{-1} + \sqrt{T \dperish L_T^{-1}} \lesssim \sqrt{T}$. Here, our algorithm achieves $\Exp{\Deff} \lesssim L_T^{-1} + \sqrt{T \dperish L_T^{-1}} + T\dperish.$ Since $L_T^{-1} \lesssim \sqrt{T}$, \eqref{eq:eff-lb} reduces to $\Exp{\Deff} \gtrsim T\dperish + L_T^{-1}.$ It is easy to check that $\sqrt{T\dperish L_T^{-1}} \lesssim \max\{L_T^{-1}, T\dperish\}$, which completes the tightness argument.

\subsection{\cref{sec:performance-guarantee} Omitted Proofs}\label{apx:performance-guarantee}
 
\subsubsection{Proof of \cref{cor:main_theorem_interpret}}\label{apx:main_theorem_interpret}

\begin{rproof}
Consider first the case where $T\dperish \lesssim L_T^{-1}$. 
Then:
\begin{align*}
L_T^{-1} + \sqrt{T\dperish L_T^{-1}}  \lesssim L_T^{-1} \lesssim \sqrt{T}, 
 \end{align*}
since $L_T \gtrsim 1/\sqrt{T}$ by assumption.
Thus, 
\begin{align*}
\Deff \lesssim \min\left\{\sqrt{T},L_T^{-1} + \sqrt{T\dperish L_T^{-1}} \right\} + T\dperish \lesssim L_T^{-1},
\end{align*}
where again we've used the assumption that $T\dperish \lesssim L_T^{-1}.$

% Under this we note that
% \[
%     L_T^{-1}(1 + \sqrt{T L_T \dperish}) \le 2 L_T^{-1} \leq \sqrt{T}
% \]
For the bound on $\Denv$, we use the facts that $L_T \gtrsim 1/\sqrt{T}$ and $\dperish \lesssim 1/\sqrt{T}$ to obtain: 
$$\Denv \lesssim \max\{L_T,\dperish + 1/\sqrt{T}\} \lesssim L_T.$$

\medskip

Suppose now $T\dperish \gtrsim L_T^{-1}$. In this case:
\begin{align*}
L_T^{-1} + \sqrt{T\dperish L_T^{-1}}  \lesssim T\dperish. 
\end{align*}
Using the fact that $\dperish \lesssim 1/\sqrt{T}$, we obtain:
\begin{align*}
\Deff \lesssim \min\left\{\sqrt{T},L_T^{-1} + \sqrt{T\dperish L_T^{-1}} \right\} + T\dperish \lesssim T\dperish.
\end{align*}

For the bound on $\Denv$, we similarly have $\Denv \lesssim L_T$, since $\dperish \lesssim 1/\sqrt{T}\lesssim L_T$, by assumption.
\end{rproof}

\subsection{\cref{sec:analysis} Omitted Proofs}\label{apx:analysis}

\subsubsection{Proof of \cref{lem:hoeffding_app-perish}}\label{apx:hoeffding_app-perish}
\begin{rproof}%{\cref{lem:hoeffding_app-perish}}
Fix $t' < t$. Recall, for all $\tau \in [T]$, $\rho_{\tau,\theta} \geq |N_{\tau,\theta}-\Exp{N_{\tau,\theta}}|$, which implies \[N_{\tau,\theta} \in  \left[\Exp{N_{\tau, \theta}} - \rho_{\tau, \theta}, \Exp{N_{\tau, \theta}} + \rho_{\tau, \theta}\right].\] Thus, from a simple application of Hoeffding's inequality (\cref{thm:hoeffding}):
\begin{align}\label{eq:hoeffding}
\Pr(|N_{(t,t']} - \Exp{N_{(t,t']}}| \geq \epsilon) \leq 2 \exp\left( - \frac{2 \epsilon^2}{\sum_\theta\sum_{\tau \in (t, t']} 4 \rho_{\tau, \theta}^2}\right)
\end{align}

We now consider our desired bound.
\begin{align*}
    & \Pr(|N_{(t,t']}-\Exp{N_{(t,t']}}| \leq \epsilon \, \forall\, t, t') \geq 1-\sum_{t,t'}\Pr(|N_{(t,t']}-\Exp{N_{(t,t']}}| \geq \epsilon) \\
    &\geq 1-\sum_{t,t'} 2\exp\left(-\frac{2\epsilon^2}{\sum_\theta\sum_{\tau\in(t,t']} 4\rho_{\tau,\theta}^2}\right) \geq 1-\sum_{t,t'} 2\exp\left(-\frac{\epsilon^2}{2|\Theta|\rho_{\max}^2(t'-t)}\right)
\end{align*}
where the first inequality follows from a union bound, the second inequality by plugging in Hoeffding's bound \eqref{eq:hoeffding}, and the third inequality by upper bounding $\rho_{\tau,\theta}$ by $\rho_{\max}$, for all $\tau \in (t,t']$.

\noindent Solving for $\epsilon$ such that $2\exp\left(-\frac{\epsilon^2}{2|\Theta|\rho_{\max}^2(t'-t)}\right) = \delta/T^2$, we obtain our result.
\end{rproof}

\subsubsection{Proof of \cref{cor:4-delta-lem}}\label{apx:4-delta-lem}

\begin{rproof}%{\cref{cor:4-delta-lem}}
The final high-probability bound follows from straightforward algebra, putting Lemmas \ref{lem:hoeffding_app-perish} and \ref{lem:p_upper_bound} together. Indeed, we have that:
\[
\Pr(\E) = 1 - \Pr(\E^c) = 1 - \Pr((\E_N \cap \E_{\Pupper})^c).
\]
Moreover:
\begin{align*}
\Pr(\E_N\cap\E_{\Pupper}) &= \Pr(\E_{\Pupper} \mid \E_N)\Pr(\E_N) \geq (1-\delta)^2 \geq 1-2\delta
\implies \Pr((\E_N \cap \E_{\Pupper})^c) \leq 2\delta.
\end{align*}
Plugging this in above, we obtain $\Pr(\E) \geq 1-2\delta$.
\end{rproof}

\subsubsection{Proof of \cref{lem:switching_point_last_perish}}\label{apx:switching_point_last_perish}
\begin{rproof}%{}
Consider first the case in which the algorithm always allocates $\underline{X}$ (i.e., $t_0 = T$). Then, the inequality is trivially satisfied and it suffices to prove the lower bound for $t_0 < T$. We have: 
\begin{align}\label{eq:budget}
B_T^{alg} &= B_{t_0}^{alg} - N_{t_0} \underline{X} - \PUA^{alg}_{t_0} - N_{(t_0, T)} \overline{X} - \PUA^{alg}_{> t_0} \notag \\
&<N_{t_0} \overline{X} + \overline{N}_{> t_0} \underline{X} + \overline{P}_{t_0}- N_{t_0} \underline{X} - \PUA^{alg}_{t_0} - N_{(t_0, T)} \overline{X} - \PUA^{alg}_{> t_0} \notag \\
&=N_{t_0}L_T + (\overline{N}_{> t_0}-N_{>t_0}) \underline{X} + {N_T\underline{X}} - N_{(t_0, T)} L_T + \overline{P}_{t_0}  - \PUA^{alg}_{\geq t_0},
\end{align}
where the first inequality follows from the fact that $B_{t_0}^{alg} < N_{t_0} \overline{X} + \overline{N}_{> t_0} \underline{X} + \overline{P}_{t_0}$ since $X_{t_0}^{alg} = \underline{X}$, and the second inequality uses $\overline{X} = \underline{X} +L_T$ and re-arranges terms. Since $\overline{X}$ was allocated at $T$, $B_T^{alg} -N_T\overline{X} \geq 0$, which then implies that $B_T^{alg}-N_T\underline{X} \geq N_TL_T$. Plugging this fact into \eqref{eq:budget} and re-arranging, we obtain:
{\begin{align}\label{eq:t_0_part}
%&N_{t_0}L_T + (\overline{N}_{> t_0}-N_{>t_0}) \underline{X} + N_T\underline{X} - N_{(t_0, T)} L_T + \overline{P}_{t_0}  - \PUA^{alg}_{\geq t_0} - N_T\overline{X} > B_T^{alg}-N_T\overline{X} \notag \\
%\implies &
N_{t_0}L_T + (\overline{N}_{> t_0}-N_{>t_0}) \underline{X} - N_{>t_0} L_T + \overline{P}_{t_0} - \PUA^{alg}_{\geq t_0} > 0.
\end{align}
}

We now upper bound the left-hand side of \eqref{eq:t_0_part}. Using the facts that $\PUA^{alg}_{\geq t_0} \geq 0$, $\underline{X} \leq \Bav$ by construction, and $\overline{P}_{t_0} \leq \overline{P}_1 = \overline{\Delta}(\underline{X})$, {we have, for $C = \sqrt{2|\Theta|\rho_{max}^2 \log(2T^2 / \pfail)}$}:
\begin{align}\label{eq:ugly-quad}
0 < N_{t_0}L_T + (\overline{N}_{> t_0}-N_{>t_0}) \underline{X} - N_{>t_0} L_T + \overline{P}_{t_0} - \PUA^{alg}_{\geq t_0} \leq \rho_{\max}L_T + 2C\Bav\sqrt{T-t_0} \notag \\ -L_T(T-t_0)+ \overline{\Delta}(\underline{X}).
\end{align}

Consider now the quadratic function $f(x) = -L_Tx^2 + 2C\Bav x + \rho_{\max}L_T + \overline{\Delta}(\underline{X})$, which has a positive root at:
\begin{align*}
x^+ &= \frac{2C\Bav + \sqrt{4C^2\Bav^2+4L_T\left(\rho_{\max}L_T + \overline{\Delta}(\underline{X})\right)}}{2L_T}\\
&=\frac{C\Bav}{L_T} + \sqrt{\frac{C^2\Bav^2}{L_T^2} + \frac{\rho_{\max}L_T + \overline{\Delta}(\underline{X})}{L_T}}\\
&\leq 2\frac{C\Bav}{L_T} + \sqrt{\frac{\overline{\Delta}(\underline{X})}{L_T}} + \sqrt{\rho_{\max}}\\
&< {c}\left(\frac{1}{L_T}+\sqrt{\frac{\overline{\Delta}(\underline{X})}{L_T}}\right),
\end{align*}
for some ${c} \in \widetilde{\Theta}(1)$. Thus, for all $x \geq {c}\left(\frac{1}{L_T}+\sqrt{\frac{\overline{\Delta}(\underline{X})}{L_T}}\right)$, $f(x) \leq 0$. Letting $x = \sqrt{T-t_0}$, we obtain that the right-hand side of \eqref{eq:ugly-quad} is non-positive for all $t_0$ such that $T-t_0 \geq {c}^2\left(\frac{1}{L_T}+\sqrt{\frac{\overline{\Delta}(\underline{X})}{L_T}}\right)^2 \iff t_0 \leq T-{c}^2\left(\frac{1}{L_T}+\sqrt{\frac{\overline{\Delta}(\underline{X})}{L_T}}\right)^2$, which would lead to a contradiction. 

Concluding, we have $t_0 >  T-{c}^2\left(\frac{1}{L_T}+\sqrt{\frac{\overline{\Delta}(\underline{X})}{L_T}}\right)^2 {\geq T-\tilde{c}^2\left(\frac{1}{L_T}+\sqrt{\frac{T\dperish}{L_T}}\right)^2}$, where the final inequality follows from the fact that $\overline{\Delta}(\underline{X}) \leq B-\overline{N}\underline{X} = \overline{N}\dperish \lesssim T\dperish$.%, where $\tilde{c} = \mu_{max} (1 + \sqrt{2|\Theta|\rho_{max}^2 \log(T^2 / \delta)})$.
\end{rproof}

\subsubsection{Proof of \cref{lem:lower_alg_comparison}}
\label{apx:proof_lower_alg_comparison}

\begin{rproof} We show the two properties by induction on $t$.

\begin{figure}
\centering     %%% not \center
\ifdefined\opre

 \else

 \fi

\begin{subfigure}[b]{0.3\paperwidth}
\scalebox{1}{
\tikzset{every picture/.style={line width=0.75pt}} %set default line width to 0.75pt        
\input{figures/perishing_proof_tikz/case_one}
}
\caption{Case I}
\end{subfigure} \quad
\begin{subfigure}[b]{0.3\paperwidth}
\scalebox{1}{
\tikzset{every picture/.style={line width=0.75pt}} %set default line width to 0.75pt        
\input{figures/perishing_proof_tikz/case_two}
}
\caption{Case II}
\end{subfigure} \quad
\begin{subfigure}[b]{0.3\paperwidth}
\scalebox{1 }{
\tikzset{every picture/.style={line width=0.75pt}} %set default line width to 0.75pt        
\input{figures/perishing_proof_tikz/case_three}
}
\caption{Case III}
\end{subfigure}
\ifdefined\opre \caption{Illustration of the three cases in the induction step of the proof of the fact that $\PUA_{t+1}^{alg} \leq \PUA_{t+1}(\underline{X})$ (\cref{lem:lower_alg_comparison}). Here, we assume $\setB{t+1}^{alg} \subseteq \setB{t+1}(\underline{X})$. The squares across all three plots show the resources in $\setB{t+1}(\underline{X})$, ordered left to right according to $\sigma$; the dashed squares correspond to $\setB{t+1}^{alg}$. The gray-shaded region corresponds to resource fractionally allocated by the $\underline{X}$ process at the beginning of $t+1$, and the cross-hatched region to the set of resources allocated by our algorithm. Finally, the red arrow corresponds to the resource $b$ considered in each case.
}
\else
\caption{Illustration of the three cases in the induction step of the proof of the fact that $\PUA_{t+1}^{alg} \leq \PUA_{t+1}(\underline{X})$ (\cref{lem:lower_alg_comparison}). Here, we assume $\setB{t+1}^{alg} \subseteq \setB{t+1}(\underline{X})$. The squares across all three plots show the resources in $\setB{t+1}(\underline{X})$, ordered left to right according to $\sigma$; the dashed squares correspond to $\setB{t+1}^{alg}$. The gray-shaded region corresponds to resource fractionally allocated by the $\underline{X}$ process at the beginning of $t+1$, and the cross-hatched region to the set of resources allocated by our algorithm. Finally, the red arrow corresponds to the resource $b$ considered in each case.
}
\fi
 \label{fig:nested_pua_proof}
\end{figure}

\noindent \textbf{Base Case $t = 1$}.  By definition, $\Balg{1} = \mathcal{B} = \setB{t}(\underline{X})$.  We now argue that $\PUA_1^{alg} \leq \PUA_1(\underline{X})$.  Suppose there exists a resource $b$ which perished at the end of $t=1$.  Then, either:
\begin{enumerate}
    \item $b$ was neither allocated by our algorithm, nor under the $\underline{X}$ allocation. Hence, it perished unallocated under both allocations.
    \item $b$ was allocated by our algorithm but not by the $\underline{X}$ allocation. Hence, it perished unallocated under $\underline{X}$ but not our algorithm.
    \item $b$ was allocated under the $\underline{X}$ allocation but not by our algorithm. This could never hold, since both algorithms begin with the same set of resources and our algorithm allocated (weakly) more than $\underline{X}$, and under the same ordering $\sigma$.
\end{enumerate}

\noindent \textbf{Step case $t \rightarrow t+1$}.  We first show that $\Balg{t+1} \subseteq \setB{t+1}(\underline{X})$.  {Let $b \in \Balg{t+1}$.  Then $b$ did not perish, and moreover $b \in \Balg{t}$.} Then, by the inductive  hypothesis, $b \in \setB{t}(\underline{X})$. Consider the following cases:
\begin{enumerate}
    \item $b$ was not allocated under the $\underline{X}$ process. Then, $b \in \setB{t+1}(\underline{X})$, since it did not perish.
    \item $b$ was allocated under the $\underline{X}$ process. In this case, since the algorithm allocated (weakly) more than $\underline{X}$ according to the same ordering, this resource {\it must} have been available to both the algorithm and the $\underline{X}$ process. This then contradicts that $b$ was not allocated by the algorithm.
\end{enumerate}

We now argue that $\PUA_{t+1}^{alg} \leq \PUA_{t+1}(\underline{X})$. Suppose there exists a resource $b$ that perished at time $t+1$. We consider the following cases (see Figure \ref{fig:nested_pua_proof} for an illustration):
\begin{enumerate}
    \item $b$ was neither allocated by our algorithm, nor under the $\underline{X}$ allocation. Hence, it perished unallocated under both allocations.
    \item $b$ was allocated by our algorithm but not by the $\underline{X}$ allocation. Hence, it perished unallocated under $\underline{X}$ but not our algorithm.
    \item $b$ was allocated under the $\underline{X}$ allocation but not by our algorithm. Then, $b$ {\it must} have either perished or been allocated before $t+1$, since the set of remaining resources under our algorithm is (weakly) nested in the set of remaining resources under the $\underline{X}$ for all $t' \leq t$, by the inductive hypothesis. Thus, $b$ could not have perished at the end of $t+1$ under our algorithm's sample path.
\end{enumerate}
\end{rproof}

%% file: figures/perishing_proof_tikz/case_one.tex
% Pattern Info
 
\tikzset{
pattern size/.store in=\mcSize, 
pattern size = 5pt,
pattern thickness/.store in=\mcThickness, 
pattern thickness = 0.3pt,
pattern radius/.store in=\mcRadius, 
pattern radius = 1pt}
\makeatletter
\pgfutil@ifundefined{pgf@pattern@name@_n33m1xbw1}{
\pgfdeclarepatternformonly[\mcThickness,\mcSize]{_n33m1xbw1}
{\pgfqpoint{0pt}{0pt}}
{\pgfpoint{\mcSize}{\mcSize}}
{\pgfpoint{\mcSize}{\mcSize}}
{
\pgfsetcolor{\tikz@pattern@color}
\pgfsetlinewidth{\mcThickness}
\pgfpathmoveto{\pgfqpoint{0pt}{\mcSize}}
\pgfpathlineto{\pgfpoint{\mcSize+\mcThickness}{-\mcThickness}}
\pgfpathmoveto{\pgfqpoint{0pt}{0pt}}
\pgfpathlineto{\pgfpoint{\mcSize+\mcThickness}{\mcSize+\mcThickness}}
\pgfusepath{stroke}
}}
\makeatother
\tikzset{every picture/.style={line width=0.75pt}} %set default line width to 0.75pt        

\begin{tikzpicture}[x=0.75pt,y=0.75pt,yscale=-1,xscale=1]

%Shape: Square [id:dp11213243805688111] 
\draw   (79,51) -- (109,51) -- (109,81) -- (79,81) -- cycle ;
%Shape: Square [id:dp25825906771660456] 
\draw   (120,51) -- (150,51) -- (150,81) -- (120,81) -- cycle ;
%Shape: Square [id:dp8921841060799822] 
\draw  [pattern=_n33m1xbw1,pattern size=6pt,pattern thickness=0.75pt,pattern radius=0pt, pattern color={rgb, 255:red, 0; green, 0; blue, 0}][dash pattern={on 4.5pt off 4.5pt}] (160,51) -- (190,51) -- (190,81) -- (160,81) -- cycle ;
%Shape: Square [id:dp585806274576244] 
\draw  [dash pattern={on 4.5pt off 4.5pt}] (199,51) -- (229,51) -- (229,81) -- (199,81) -- cycle ;
%Shape: Square [id:dp011586004588473697] 
\draw   (240,51) -- (270,51) -- (270,81) -- (240,81) -- cycle ;
%Shape: Square [id:dp38249086910680963] 
\draw   (280,51) -- (310,51) -- (310,81) -- (280,81) -- cycle ;
%Shape: Rectangle [id:dp8614288268637429] 
\draw  [draw opacity=0][fill={rgb, 255:red, 144; green, 144; blue, 144 }  ,fill opacity=1 ] (79,51) -- (101,51) -- (101,81) -- (79,81) -- cycle ;
%Straight Lines [id:da8707389102866245] 
\draw [color={rgb, 255:red, 255; green, 0; blue, 0 }  ,draw opacity=1 ][line width=1.5]    (135,111) -- (135,87) ;
\draw [shift={(135,84)}, rotate = 90] [color={rgb, 255:red, 255; green, 0; blue, 0 }  ,draw opacity=1 ][line width=1.5]    (14.21,-6.37) .. controls (9.04,-2.99) and (4.3,-0.87) .. (0,0) .. controls (4.3,0.87) and (9.04,2.99) .. (14.21,6.37)   ;

\end{tikzpicture}

%% file: figures/perishing_proof_tikz/case_two.tex
% Pattern Info
 
\tikzset{
pattern size/.store in=\mcSize, 
pattern size = 5pt,
pattern thickness/.store in=\mcThickness, 
pattern thickness = 0.3pt,
pattern radius/.store in=\mcRadius, 
pattern radius = 1pt}
\makeatletter
\pgfutil@ifundefined{pgf@pattern@name@_5e6r2x0kt}{
\pgfdeclarepatternformonly[\mcThickness,\mcSize]{_5e6r2x0kt}
{\pgfqpoint{0pt}{0pt}}
{\pgfpoint{\mcSize}{\mcSize}}
{\pgfpoint{\mcSize}{\mcSize}}
{
\pgfsetcolor{\tikz@pattern@color}
\pgfsetlinewidth{\mcThickness}
\pgfpathmoveto{\pgfqpoint{0pt}{\mcSize}}
\pgfpathlineto{\pgfpoint{\mcSize+\mcThickness}{-\mcThickness}}
\pgfpathmoveto{\pgfqpoint{0pt}{0pt}}
\pgfpathlineto{\pgfpoint{\mcSize+\mcThickness}{\mcSize+\mcThickness}}
\pgfusepath{stroke}
}}
\makeatother
\tikzset{every picture/.style={line width=0.75pt}} %set default line width to 0.75pt        

\begin{tikzpicture}[x=0.75pt,y=0.75pt,yscale=-1,xscale=1]

%Shape: Square [id:dp09357775651496847] 
\draw   (79,51) -- (109,51) -- (109,81) -- (79,81) -- cycle ;
%Shape: Square [id:dp15662402472798687] 
\draw   (120,51) -- (150,51) -- (150,81) -- (120,81) -- cycle ;
%Shape: Square [id:dp22199261524776293] 
\draw  [pattern=_5e6r2x0kt,pattern size=6pt,pattern thickness=0.75pt,pattern radius=0pt, pattern color={rgb, 255:red, 0; green, 0; blue, 0}][dash pattern={on 4.5pt off 4.5pt}] (160,51) -- (190,51) -- (190,81) -- (160,81) -- cycle ;
%Shape: Square [id:dp028470169530659062] 
\draw  [dash pattern={on 4.5pt off 4.5pt}] (199,51) -- (229,51) -- (229,81) -- (199,81) -- cycle ;
%Shape: Square [id:dp9218290841121348] 
\draw   (240,51) -- (270,51) -- (270,81) -- (240,81) -- cycle ;
%Shape: Square [id:dp7291273279648698] 
\draw   (280,51) -- (310,51) -- (310,81) -- (280,81) -- cycle ;
%Shape: Rectangle [id:dp482412119475383] 
\draw  [draw opacity=0][fill={rgb, 255:red, 144; green, 144; blue, 144 }  ,fill opacity=1 ] (79,51) -- (101,51) -- (101,81) -- (79,81) -- cycle ;
%Straight Lines [id:da9103633355602507] 
\draw [color={rgb, 255:red, 255; green, 0; blue, 0 }  ,draw opacity=1 ][line width=1.5]    (175,112) -- (175,88) ;
\draw [shift={(175,85)}, rotate = 90] [color={rgb, 255:red, 255; green, 0; blue, 0 }  ,draw opacity=1 ][line width=1.5]    (14.21,-6.37) .. controls (9.04,-2.99) and (4.3,-0.87) .. (0,0) .. controls (4.3,0.87) and (9.04,2.99) .. (14.21,6.37)   ;

\end{tikzpicture}

%% file: figures/perishing_proof_tikz/case_three.tex
% Pattern Info
 
\tikzset{
pattern size/.store in=\mcSize, 
pattern size = 5pt,
pattern thickness/.store in=\mcThickness, 
pattern thickness = 0.3pt,
pattern radius/.store in=\mcRadius, 
pattern radius = 1pt}
\makeatletter
\pgfutil@ifundefined{pgf@pattern@name@_h1knmyavb}{
\pgfdeclarepatternformonly[\mcThickness,\mcSize]{_h1knmyavb}
{\pgfqpoint{0pt}{0pt}}
{\pgfpoint{\mcSize}{\mcSize}}
{\pgfpoint{\mcSize}{\mcSize}}
{
\pgfsetcolor{\tikz@pattern@color}
\pgfsetlinewidth{\mcThickness}
\pgfpathmoveto{\pgfqpoint{0pt}{\mcSize}}
\pgfpathlineto{\pgfpoint{\mcSize+\mcThickness}{-\mcThickness}}
\pgfpathmoveto{\pgfqpoint{0pt}{0pt}}
\pgfpathlineto{\pgfpoint{\mcSize+\mcThickness}{\mcSize+\mcThickness}}
\pgfusepath{stroke}
}}
\makeatother
\tikzset{every picture/.style={line width=0.75pt}} %set default line width to 0.75pt        

\begin{tikzpicture}[x=0.75pt,y=0.75pt,yscale=-1,xscale=1]

%Shape: Square [id:dp43616324533464046] 
\draw   (79,51) -- (109,51) -- (109,81) -- (79,81) -- cycle ;
%Shape: Square [id:dp25230886622787185] 
\draw   (120,51) -- (150,51) -- (150,81) -- (120,81) -- cycle ;
%Shape: Square [id:dp553985313946842] 
\draw  [pattern=_h1knmyavb,pattern size=6pt,pattern thickness=0.75pt,pattern radius=0pt, pattern color={rgb, 255:red, 0; green, 0; blue, 0}][dash pattern={on 4.5pt off 4.5pt}] (160,51) -- (190,51) -- (190,81) -- (160,81) -- cycle ;
%Shape: Square [id:dp706192907673109] 
\draw  [dash pattern={on 4.5pt off 4.5pt}] (199,51) -- (229,51) -- (229,81) -- (199,81) -- cycle ;
%Shape: Square [id:dp18015802978742967] 
\draw   (240,51) -- (270,51) -- (270,81) -- (240,81) -- cycle ;
%Shape: Square [id:dp650850085675371] 
\draw   (280,51) -- (310,51) -- (310,81) -- (280,81) -- cycle ;
%Shape: Rectangle [id:dp7555306302001175] 
\draw  [draw opacity=0][fill={rgb, 255:red, 144; green, 144; blue, 144 }  ,fill opacity=1 ] (79,51) -- (101,51) -- (101,81) -- (79,81) -- cycle ;
%Straight Lines [id:da17528430147505736] 
\draw [color={rgb, 255:red, 255; green, 0; blue, 0 }  ,draw opacity=1 ][line width=1.5]    (95,112) -- (95,88) ;
\draw [shift={(95,85)}, rotate = 90] [color={rgb, 255:red, 255; green, 0; blue, 0 }  ,draw opacity=1 ][line width=1.5]    (14.21,-6.37) .. controls (9.04,-2.99) and (4.3,-0.87) .. (0,0) .. controls (4.3,0.87) and (9.04,2.99) .. (14.21,6.37)   ;

\end{tikzpicture}

%% file: final_arxiv_parts/special-cases-proofs.tex
\subsection{\cref{sec:special-cases} Omitted Proofs}\label{apx:discussion}

For ease of notation, we let $\nu_t = \mathbb{E}[P_{<t}]$ for all $t\in\{2,\ldots,T\}$.

\subsubsection{Proof of \cref{prop:nec-cond-offset}}

\begin{rproof}
{
Let $t \in [T]$ be such that $\mathbb{E}[P_{<t}] > t-1$. Then:
\begin{align}\label{eq:def}
\Pr(P_{<t} \leq t-1 \ \forall \ t \geq 2) &\leq \Pr(P_{<t} \leq t-1)
= \Pr\left(\sum_{b\in \mathcal{B}} \mathds{1}\{T_b < t\} \leq t-1 \right),
\end{align}
where the second equality is by definition.

Consider first the case where $\mathcal{B}^{rand}_{<t} = \emptyset$. In this case, if $b$ perishes before $t$ with strictly positive probability, it must be that $b \in \mathcal{B}^{det}_{<t}$. Then:
\begin{align}\label{eq:det}
\Pr\left(\sum_{b\in \mathcal{B}} \mathds{1}\{T_b < t\} \leq t-1 \right) &= \Pr\left(\sum_{b\in \mathcal{B}^{det}_{<t}} \mathds{1}\{T_b < t\} \leq t-1 \right) = \Pr\left(|\mathcal{B}^{det}_{<t}| \leq t-1 \right),
\end{align}
where the second equality follows from the fact that items in $\mathcal{B}^{det}_{<t}$ perish before $t$ with probability 1. By the same reasoning:
\begin{align*}
&t-1 < \mathbb{E}[P_{<t}] = \sum_{b\in\mathcal{B}}\Pr(T_b<t) = \sum_{b \in\mathcal{B}^{det}_{<t}} \Pr(T_b < t) = |\mathcal{B}^{det}_{<t}|
\implies \Pr(|\mathcal{B}^{det}_{<t}| \leq t-1) = 0.
\end{align*}
Plugging this back into \eqref{eq:def}, we obtain $\Pr(P_{<t} \leq t-1 \ \forall \ t\geq 2) = 0$.

\medskip 

Consider now the case where $\mathcal{B}^{rand}_{<t} \neq \emptyset$. The goal is to show the existence of $\epsilon$ such that \mbox{$\Pr(P_{<t} \leq t-1 \ \forall t \geq 2) \leq \epsilon$}.  Define the random variable: 
$$Y_b = \mathds{1}\left\{T_b < t\right\}-\Pr(T_b < t), \quad  b \in \mathcal{B}^{rand}_{< t}.$$ By construction, $\mathbb{E}[Y_b] = 0, 0 < \mathbb{E}[Y_b^2] \leq 1$, and $\mathbb{E}[|Y_b|^3] \leq 1$. We have:
\begin{align*}
\Pr(P_{<t} \leq t-1) &= \Pr\left(\sum_{b\in \mathcal{B}^{rand}_{< t}} \mathds{1}\{T_b < t\} \leq t-1 -|\mathcal{B}^{det}_{<t}|\right)\\&=\Pr\left(\sum_{b\in \mathcal{B}^{rand}_{< t}} Y_b \leq t-1 -|\mathcal{B}^{det}_{<t}| - \sum_{b\in\mathcal{B}^{rand}_{<t}}\Pr(T_b < t)\right).
\end{align*}
By assumption, $\mathbb{E}[P_{<t}] = |\mathcal{B}^{det}_{<t}| + \sum_{b\in\mathcal{B}^{rand}_{<t}}\Pr(T_b < t) > t-1$. Hence,
\begin{align*}
\Pr(P_{<t} \leq t-1) 
&\leq \Pr\left(\sum_{b\in \mathcal{B}^{rand}_{< t}} Y_b \leq 0\right) = \Pr\left(\frac{\sum_{b\in \mathcal{B}^{rand}_{< t}} Y_b}{\sqrt{\sum_{b\in \mathcal{B}^{rand}_{< t}} \mathbb{E}[Y_b^2]}} \leq 0\right).
\end{align*}
Let $\Phi(\cdot)$ denote the cdf of the standard normal distribution. By the Berry-Esseen Theorem,
\begin{align*}
\Pr\left(\frac{\sum_{b\in \mathcal{B}^{rand}_{< t}} Y_b}{\sqrt{\sum_{b\in \mathcal{B}^{rand}_{< t}} \mathbb{E}[Y_b^2]}} \leq 0\right)&\leq \Phi(0)+\left(\sum_{b\in\mathcal{B}^{rand}_{< t}}\mathbb{E}[Y_b^2]\right)^{-3/2}\cdot\sum_{b\in\mathcal{B}^{rand}_{< t}}\mathbb{E}[|Y_b|^3] \\
&= \frac12 + \left(\sum_{b\in\mathcal{B}^{rand}_{< t}}\text{Var}[\mathds{1}\{T_b < t\}]\right)^{-3/2}\cdot\sum_{b\in\mathcal{B}^{rand}_{< t}}\mathbb{E}[|Y_b|^3] \\
&= \frac12 + \left(\text{Var}\left[\sum_{b\in\mathcal{B}^{rand}_{< t}}\mathds{1}\{T_b < t\}\right]\right)^{-3/2}\cdot\sum_{b\in\mathcal{B}^{rand}_{< t}}\mathbb{E}[|Y_b|^3] \\
&\leq \frac12 + \text{Var}[P_{<t}]^{-3/2} \cdot T \\
&= \frac12 + {\stdpt^{-3}} \cdot T.
\end{align*}
Putting this all together, we obtain:
\begin{align*}
\Pr(P_{<t} \leq t-1 \ \forall t \geq 2) \leq \frac12 +{\stdpt^{-3}} \cdot T.
\end{align*}
As a result, the perishing process cannot be $\delta$-offset-expiring for any $\delta < \frac12-{\stdpt^{-3}} \cdot T$.
}
\end{rproof}

\subsubsection{Proof of \cref{prop:basic-conditions}}
\begin{rproof}
Recall, $P_{<t} = |\mathcal{B}^{det}_{<t}| + \sum_{b \in \mathcal{B}^{rand}_{<t}} \mathds{1}\{T_b < t\}$. By Chebyshev's inequality:
\begin{align*}
\mathbb{P}(P_{<t} > t-1) &= \Pr\left(P_{<t}-\nu_t \geq t-\nu_t\right) \leq \left(\frac{\stdpt}{t-\nu_t}\right)^2.
\end{align*}
%where the final inequality follows from the upper bound on $\nu_t$.

Similarly, by Hoeffding's inequality we have:
\begin{align*}
\mathbb{P}(P_{<t} > t-1) &= \Pr\left(\sum_{b\in\mathcal{B}^{rand}_{<t}}\mathds{1}\{T_b < t\}-\Pr(T_b < t) {\geq} t - \nu_t\right)\leq \exp\left(-\frac{{2}\left(t-\nu_t\right)^2}{|\mathcal{B}^{rand}_{<t}|}\right).
\end{align*}

Then, via a straightforward union bound we have:
\begin{align*}
\mathbb{P}(P_{<t} \leq t-1 \ \forall \, t \geq 2) \geq 1-\sum_{t=2}^T \Pr(P_{<t} > t-1) \geq  1-\sum_{t=2}^T \min\left\{\left(\frac{\stdpt}{t-\nu_t}\right)^2, \exp\left(-\frac{{2}(t-\nu_t)^2}{|\mathcal{B}^{rand}_{<t}|}\right)\right\}.
\end{align*}
\end{rproof}

\subsubsection{Proof of \cref{thm:geometric_perishing}}
\begin{rproof}

For all $b \in \mathcal{B}$, we have:
\begin{align}\label{eq:bounding-std-and-nu}
\mathbb{P}(T_b \leq t-1) = 1-(1-p)^{t-1}
\implies 
&\begin{cases}\nu_t &= T\left(1-(1-p)^{t-1}\right) \\ \text{Var}[P_{<t}] &= T\left(1-(1-p)^{t-1}\right)(1-p)^{t-1}.
\end{cases}
\end{align}

By \cref{prop:basic-conditions}, $\delta$-offset-expiry holds for all $\delta \geq \sum_{t=2}^T\left(\frac{\stdpt}{t-\nu_t}\right)^2$. By \eqref{eq:bounding-std-and-nu}:
\begin{align*}
\left(\frac{\stdpt}{t-\nu_t}\right)^2 &=
\frac{T\left(1-(1-p)^{t-1}\right)(1-p)^{t-1}}{\left(t-T\left(1-(1-p)^{t-1}\right)\right)^2}
\end{align*}
Using the fact that $(1-p)^{t-1} \geq 1-(t-1)p$, we have:
\begin{align*}
t-T(1-(1-p)^{t-1}) \geq t-T(t-1)p > 0,
\end{align*}
where the final inequality follows from the assumption that $p \leq 1/T$. Taking derivatives, it is easy to show that the function $f(x) = \frac{(1-x)x}{(t-T(1-x))^2}$ is decreasing for $t-T(1-x) > 0$. Hence, leveraging the same lower bound on $(1-p)^{t-1}$ we have:
\begin{align*}
\left(\frac{\stdpt}{t-\nu_t}\right)^2 &\leq
\frac{T(t-1)p(1-(t-1)p)}{\left(t-T(t-1)p\right)^2} \leq \frac{Ttp(1-(t-1)p)}{t^2\left(1-Tp\right)^2} \leq \frac{Tp(1-(t-1)p)}{(t-1)(1-Tp)^2} \\
\implies \sum_{t=2}^T \left(\frac{\stdpt}{t-\nu_t}\right)^2 & \leq \frac{Tp}{(1-Tp)^2}\sum_{t=2}^T\left(\frac{1}{t-1}-p\right) \\ &=  \frac{Tp}{(1-Tp)^2}\left(\sum_{t=1}^{T-1}\frac{1}{t}-p(T-1)\right) \\
&\leq \frac{Tp}{(1-Tp)^2}\cdot 2\log T.
\end{align*}
Hence, the perishing process is $\delta$-offset-expiring for all $\delta \geq 2\log T \cdot \frac{Tp}{(1-Tp)^2}$.

We conclude by showing the lower bound on $\underline{X}$. 
By definition,  
\begin{align*}
\overline{\Delta}(X) &= \mu(X) + \frac12\left(\log(3 \log(T)/\delta) + \sqrt{\log^2(3 \log(T)/\delta)+8\mu(X)\log(3 \log(T)/\delta)}\right)\\
&\leq \mu(X) + \frac12\left(2\log(3 \log(T)/\delta) +\frac{1}{2\log(3 \log(T)/\delta)}\cdot 8\mu(X)\log(3 \log(T)/\delta)\right)\\
&=3\mu(X) + \log(3 \log(T)/\delta),
\end{align*}
where the inequality follows from concavity. Moreover:
\begin{align*}
 \mu(X) = \sum_b \Pr(T_b < \min\{T, \taub{1 \mid X,\sigma}\}) & \leq \sum_b \Pr(T_b < T) = T(1-(1-p)^{T-1})\\
 \implies \overline{\Delta}(X) &\leq 3 T(1-(1-p)^{T-1}) + \log(3 \log(T)/\delta)\\
 &\leq 3T(T-1)p+\log(3 \log(T)/\delta).
\end{align*}
Since any feasible stationary allocation $X$ must satisfy $X \leq \frac{T-\overline{\Delta}(X)}{T}$, it suffices to have:
\begin{align}\label{eq:geom-rhs}
X &\leq \frac{T-3 T(T-1)p-\log(3 \log(T)/\delta)}{T}=1-3(T-1)p-\frac{\log(3 \log(T)/\delta)}{T}.
\end{align}
Noting that the right-hand side of \eqref{eq:geom-rhs} is non-negative for $\delta \geq 3 \log T \cdot \exp(-(T - 3T^2 p))$, and that $2\log T \cdot \frac{Tp}{(1-Tp)^2}\geq 3 \log T \cdot \exp(-(T - 3T^2 p))$ for $p = o(1)$, we obtain the result.
\end{rproof}

\subsubsection{Proof of \cref{prop:gen_distribution_example}}
\label{app:prop_gen_distr_proof}

\begin{rproof}
For ease of notation, we let $\mu_b = \Exp{T_b}$.
For a stationary allocation $X = 1-T^{-\alpha}$, let $\mu(X) = \sum_b \Pr(T_b < \min\{T, \lceil \frac{\sigma(b)}{X}\rceil\})$, and $\mu_b = \mathbb{E}[T_b]$. By Chebyshev's inequality, we have:
\begin{align}\label{eq:mu-gen-ex}
\mu(X) &\leq \sum_b \Pr \left( T_b-\mu_b \leq \min\left\{T, \lceil\frac{\sigma(b)}{X}\rceil \right\} -\mu_b \right) \leq \sum_b \Pr\left(|T_b - \mu_b| \geq \mu_b - \min\left\{T, \lceil\frac{\sigma(b)}{X}\rceil \right\}\right) \notag \\
& \leq \sum_b \frac{\Var{T_b}}{\left( \mu_b - \min\{T, \lceil\frac{\sigma(b)}{X}\rceil \} \right)^2},
\end{align}
where we used the assumption that $\mu_b > \min\{T, \lceil\frac{\sigma(b)}{X}\rceil \}$ for all $b \in \mathcal{B}$.
However, by assumption we have
\[ \mu(X) \leq \sum_b \frac{\Var{T_b}}{\left( \mu_b - \min\{T, \lceil\frac{\sigma(b)}{X}\rceil \} \right)^2} \leq \frac{1}{2}T^{1-\alpha}\]

Moreover, by definition:
\begin{align*}
\overline{\Delta}(X) &= \mu(X) + \frac{1}{2}\left(\log(3 \log(T) / \delta) + \sqrt{\log^2(3 \log(T) / \delta) + 8 \mu(X) \log(3 \log(T) / \delta)}\right) \\
& \leq \mu(X) + \log(3 \log(T)/\delta) + \sqrt{2\mu(X)\log(3 \log(T)/\delta)}.
\end{align*}

Since any feasible $X$ must satisfy ${X} \leq 1-\overline{\Delta}({X})/T$, we have that $\overline{\Delta}(X)\leq T^{1-\alpha}$ for $X = 1-T^{-\alpha}$. Thus, it suffices for $\mu(X)$ to satisfy
\begin{align*}
    & \mu(X) +  \sqrt{2\mu(X) \log(3 \log(T)/\delta)} + \log(3 \log(T)/\delta) \leq T^{1-\alpha} %w\\
    % & \iff \mu(X) \leq T^{1-\alpha}-\sqrt{\left(2T^{1-\alpha}-\log(3 \log(T)/\delta)\right)\log(3 \log(T)/\delta)}
\end{align*}
%for $\delta \geq e^{-2T^{1-\alpha}}$. 
We have that $\mu(X) \leq \frac12T^{1-\alpha}$ satisfies this inequality for all $\delta \geq 3\log(T)e^{-\frac18T^{1-\alpha}}$.
\end{rproof}

%% file: final_arxiv_parts/experiments_full.tex
\section{Simulation details}
\label{sec:experiment_details}

\paragraph{Computing Infrastructure.} The experiments were conducted on a personal computer with an Apple M2, 8-core processor and 16.0GB of RAM.

\begin{figure}
\centering
\includegraphics[scale=0.5]{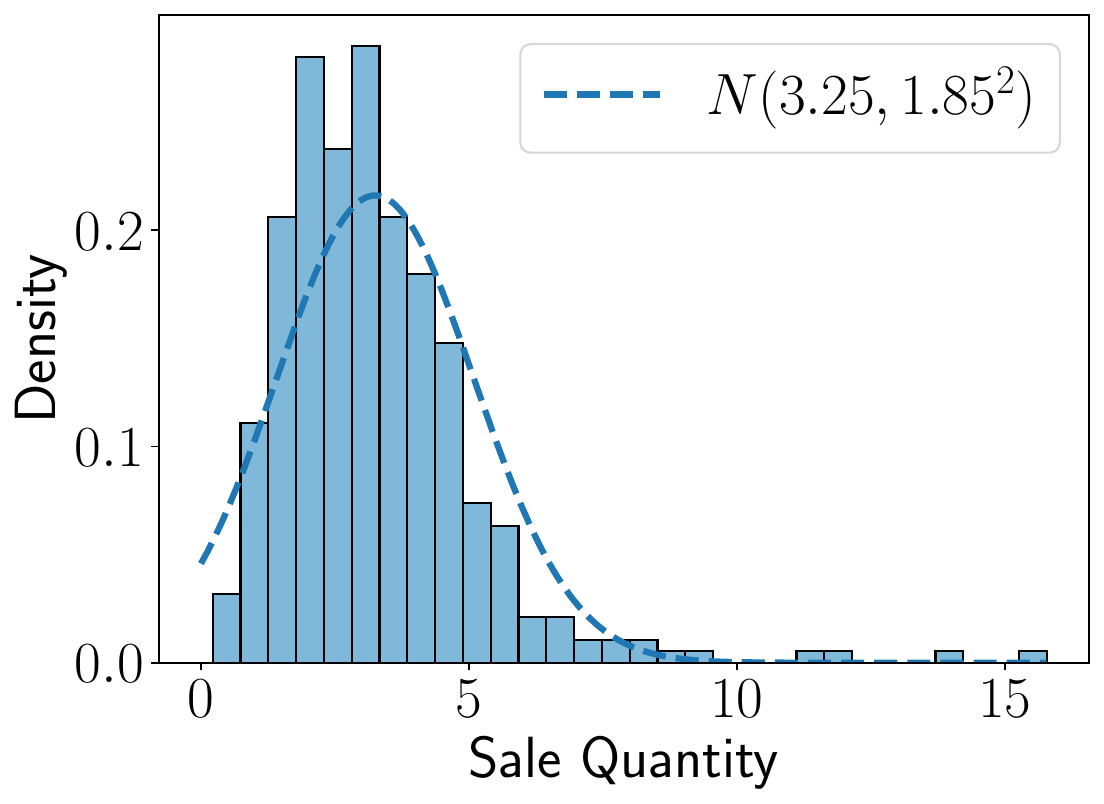}
\caption{Empirical distribution of daily ginger sales with fitted normal distribution using the sample mean and standard deviation.}
\label{fig:ginger_plot}
\end{figure}

\paragraph{Real-World Simulations.} We use the dataset of \citet{keskin2022data}, which contains detailed daily information on orders, inventory, sales, prices, costs, and holiday indicators for a variety of perishable fruits and vegetables across multiple stores, for a leading supermarket chain in  China. Following \citet{keskin2022data}, we focus on ``ginger'' in store ``A,'' as this product experiences no stock-outs over the sample period. The data span the full calendar year of 2013 and, for each day, include beginning inventory, replenishment quantity, realized sales, ending inventory, price, and a holiday indicator.  All inventory and sales quantities are normalized in the dataset, so there is no explicit notion of physical ``units.'' Throughout, we treat these normalized quantities as inventory units for modeling purposes.

As in \citet{keskin2022data}, we assume demand is normally distributed and fit the mean and standard deviation using the full year of data. Specifically, we treat observed daily sales as a proxy for demand and fit a normal distribution to this data, obtaining mean $\mu = 3.2$ and standard deviation $\sigma = 1.85$.\footnote{This specification ignores potential effects from price changes, seasonality, or holiday-related demand shifts, as these effects are not the focus of this work.}  See \Cref{fig:ginger_plot} for a histogram of the true sales data against the normal approximation.

To estimate the perishability rate, we fit a maximum likelihood estimator under the assumption that each unit independently perishes each day with probability $p$.  We compute the number of units that perish on each day $t$ as follows:
\[
\text{Perish}_t = \text{BeginStockQty}_t + \text{Restock}_t - \text{SalesQty}_t - \text{EndStockQty}_t,
\]
where $\text{Restock}_t$ is used to denote the replenishment quantity on day $t$. (The dataset reports the replenishment quantity separately from the initial inventory level for each day $t$.)
Under the assumption that each unit perishes i.i.d. with probability $p$ on each day, the number of perished units on day $t$ is binomially distributed with $\text{BeginStockQty}_t$ trials and parameter $p$. Hence, the maximum likelihood estimate of $p$ is given by:
\[
\hat{p} = \frac{\sum_t \text{Perish}_t}{\sum_t \text{BeginStockQty}_t}.
\]
On the ginger dataset, this results in $\hat{p} = 0.0024.$

\newpage

\subsection{Additional Results}
\label{sec:additional_numerical_results}

\begin{figure}[!h]
\ifdefined\opre
\centering
\else \fi
    \begin{subfigure}{.45\textwidth}
        \centering
        \includegraphics[width=\textwidth]{figures/geometric_figures/geometric_perishing_tradeoff_0-3.pdf}
        \caption{\centering $T_b \sim \emph{Geometric}(T^{-1.3}): \underline{X} = 0.93,$ $\dperish = 0.17, \Pr(\E_{OE}) = 0.99$}
    \end{subfigure}
    \begin{subfigure}{.45\textwidth}
        \centering
        \includegraphics[width=\textwidth]{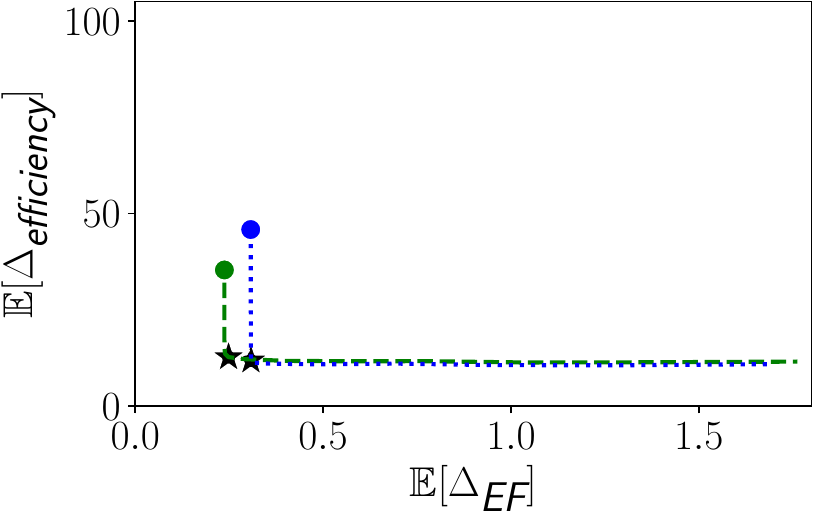}
        \caption{\centering $T_b \sim \emph{Geometric}(T^{-1.5}): \underline{X} = 1.03,$ $\dperish = 0.07, \Pr(\E_{OE}) = 1$}
    \end{subfigure}
    \begin{subfigure}{.45\textwidth}
        \centering
        \includegraphics[width=\textwidth]{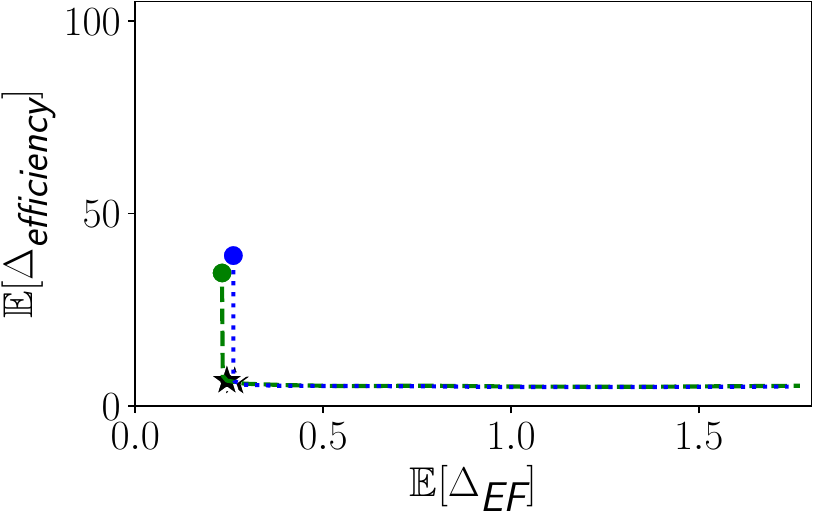}
        \caption{\centering $T_b \sim \emph{Geometric}(T^{-1.7}): \underline{X} = 1.06,$ $\dperish = 0.03, \Pr(\E_{OE}) = 1$}
    \end{subfigure}
        \begin{subfigure}{.45\textwidth}
        \centering
        \includegraphics[width=\textwidth]{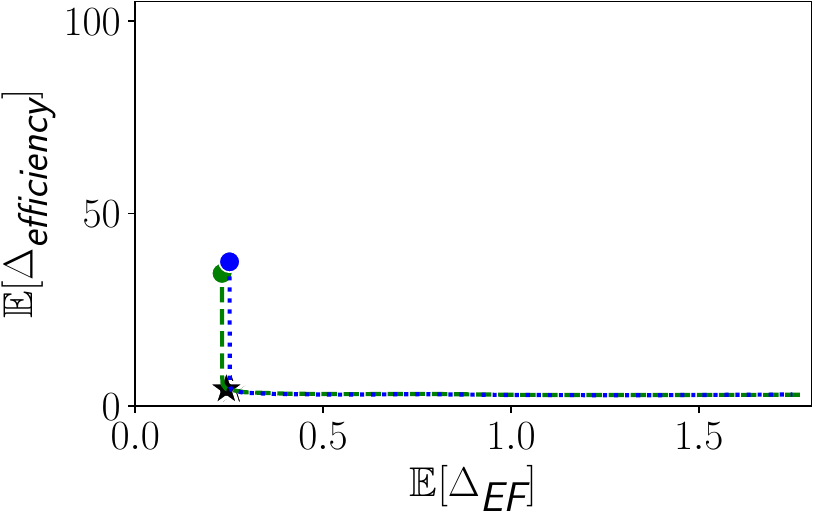}
        \caption{\centering $T_b \sim \emph{Geometric}(T^{-1.9}): \underline{X} = 1.07,$ $\dperish = 0.02, \Pr(\E_{OE}) = 1$}
    \end{subfigure}
    \begin{subfigure}{\textwidth}
        \centering
        \includegraphics[width = \textwidth]{figures/geometric_figures/tradeoff_legend.pdf}
    \end{subfigure}
    \ifdefined\opre     \caption{\centering Empirical trade-off between $\Deff$ and $\Denv$ for the different algorithms under various values of $L_T$.}
 \else
    \caption{Empirical trade-off between $\Deff$ and $\Denv$ for the different algorithms under various values of $L_T$.} \fi
    \label{fig:geometric_tradeoff_app}
\end{figure}

\begin{figure}
\ifdefined\opre 
\else \fi
\centering
% \begin{subfigure}[b]{\textwidth}
% \centering
% \includegraphics[width=.8\textwidth]{figures/geometric_figures/geometric_perishing_0-1.pdf}
% \caption{$T_b \sim \emph{Geometric}(T^{-(1.1)})$}
% \end{subfigure}
\begin{subfigure}[b]{\textwidth}
\centering
\includegraphics[width=.8\textwidth]{figures/geometric_figures/geometric_perishing_2_0-3.pdf}
\caption{$T_b \sim \emph{Geometric}(T^{-1.3})$}
\end{subfigure}
\begin{subfigure}[b]{\textwidth}
\centering
\includegraphics[width=.8\textwidth]{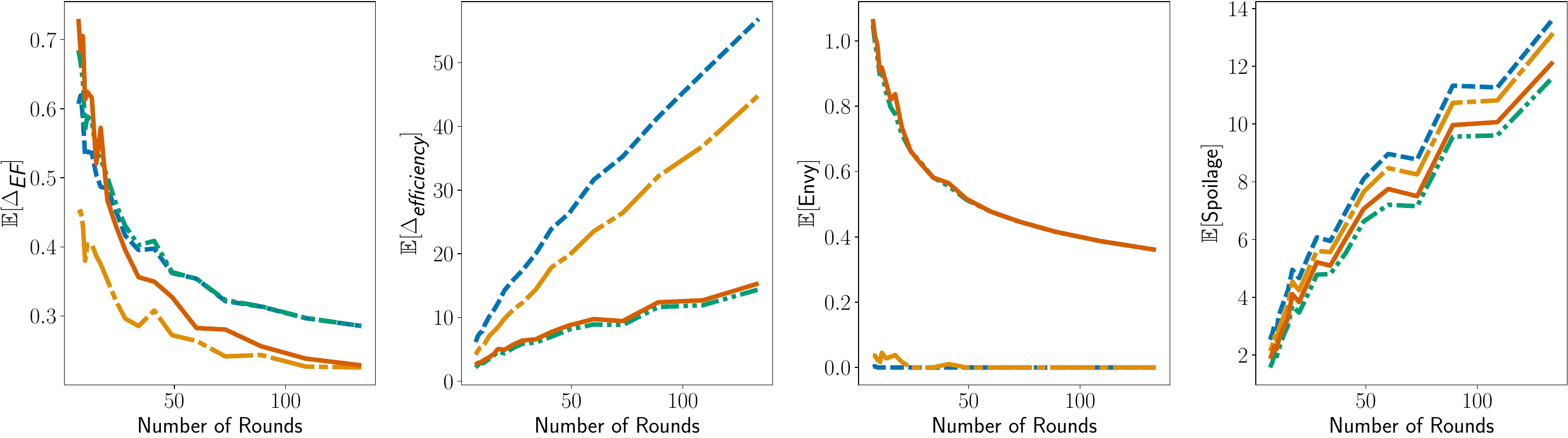}
\caption{$T_b \sim \emph{Geometric}(T^{-1.5})$}
\end{subfigure}
\begin{subfigure}[b]{\textwidth}
\centering
\includegraphics[width=.8\textwidth]{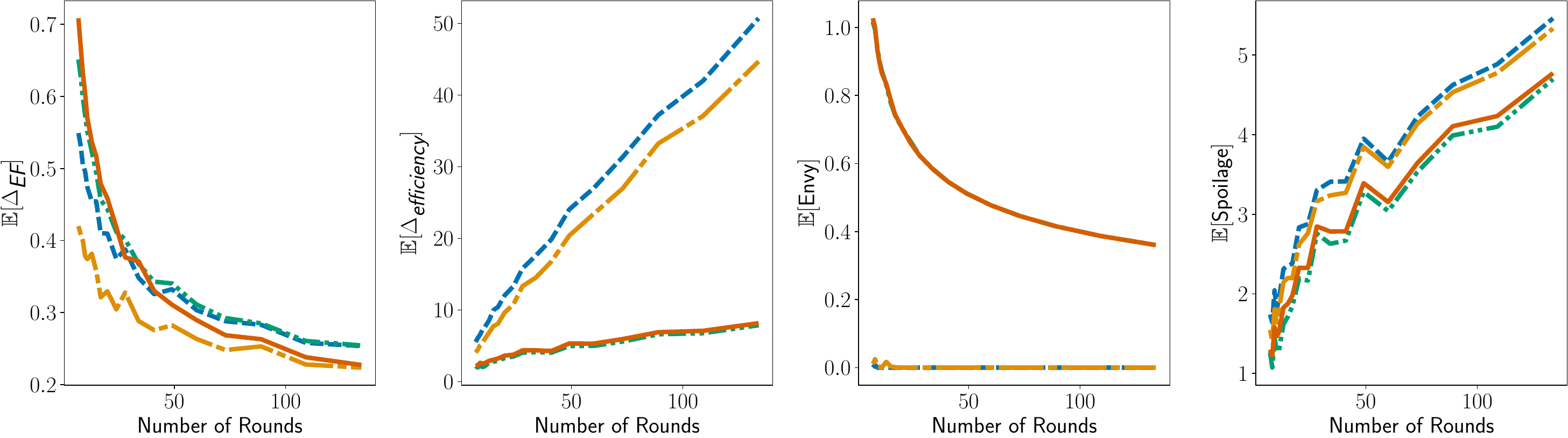}
\caption{$T_b \sim \emph{Geometric}(T^{-1.7})$}
\end{subfigure}
\begin{subfigure}[b]{\textwidth}
\centering
\includegraphics[width=.8\textwidth]{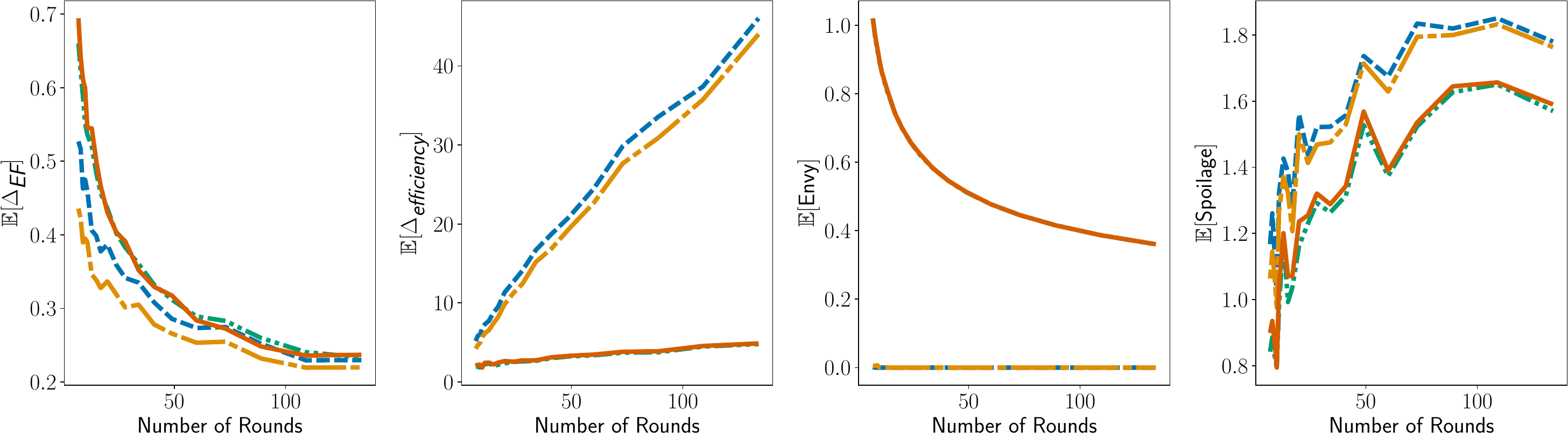}
\caption{$T_b \sim \emph{Geometric}(T^{-1.9})$}
\end{subfigure}
\begin{subfigure}{\textwidth}
\centering
\includegraphics[width=.8\paperwidth]{figures/geometric_figures/t_scaling_legend.pdf}
\end{subfigure}
\ifdefined\opre \caption{\centering Numerical results for $T_b \sim \emph{Geometric}(p)$ as described in \cref{sec:experiments_geometric}.}
 \else
\caption{Numerical results for $T_b \sim \emph{Geometric}(p)$ as described in \cref{sec:experiments_geometric}.} \fi
\label{fig:geometric_simulations_app}
\end{figure}

\clearpage

%% file: final_arxiv_parts/useful_lemmas.tex
\section{Useful lemmas}
\label{sec:lemmas}

We use the following standard theorems throughout the proof.  See, e.g. \citet{vershynin2018high} for proofs and further discussion.

\begin{lemma}[Hoeffding's Inequality]
\label{thm:hoeffding}
Let $X_1, \ldots, X_n$ be independent random variables such that $a_i \leq X_i \leq b_i$ almost surely, with $S_n = \sum_i X_i$.  Then, for all $t > 0$:
\begin{align*}
    \Pr(|S_n - \Exp{S_n}| \geq t) \leq 2 \exp\left(-\frac{2t^2}{\sum_i (b_i - a_i)^2}\right).
\end{align*}
\end{lemma}

The next is a Chernoff bound for Bernoulli random variables. See \citet{mitzenmacher2017probability}.
\begin{lemma}[Chernoff Bound for Sum of Bernoulli Random Variables]
\label{thm:chernoff_bernoulli}
Consider a sequence of Bernoulli random variables $(X_i)_{i\in[N]}$, independently distributed with probability of success $p_i \in (0,1)$. Let $X = \sum_i X_i$, and let $\mu = \Exp{X} = \sum_i p_i$.  Then, for all $\epsilon > 0$:
\[
    \Pr(X \geq (1+\epsilon)\mu) \leq \exp\left(-\frac{\epsilon^2}{2+\epsilon}\mu\right).
\]
\end{lemma}
{
\begin{corollary}
\label{cor:chernoff_bernoulli}
Consider a sequence of Bernoulli random variables $(X_i)_{i\in[N]}$, independently distributed with probability of success $p_i \in (0,1)$. Let $X = \sum_i X_i$, and let $\mu = \Exp{X} = \sum_i p_i$.  Then, for any $\delta > 0$, with probability at least $1 - \delta$ we have:
\[
    X \leq \mu + \frac{1}{2} \left(\log(1 / \delta) + \sqrt{\log^2(1 / \delta) + 8 \mu \log(1 / \delta)} \right).
\]
\end{corollary}
\begin{rproof}
    Setting the right hand side equal to $\delta$ in \cref{thm:chernoff_bernoulli} and solving for $\epsilon$, we have:
    \begin{align*}
    &\frac{\epsilon^2}{2+\epsilon}\mu = \log(1/\delta)
    \iff \epsilon^2\mu-\epsilon\log(1/\delta)-2\log(1/\delta) = 0\\
    \iff &\epsilon = \frac{\log(1/\delta) + \sqrt{\log^2(1/\delta) + 8\mu\log(1/\delta)}}{2\mu}.
    \end{align*}
    Plugging this value of $\epsilon$ into $(1 + \epsilon) \mu$ we have the result.
\end{rproof}
}